\documentclass[11pt]{amsart}
\usepackage{palatino, mathpazo}
\usepackage[mathscr]{eucal}
\usepackage{amssymb}
\usepackage[all]{xy}
\usepackage{graphicx}
\usepackage{texdraw}
\usepackage{enumitem}
\usepackage{color}

\newtheorem{theorem}{Theorem}[section]

\newtheorem{corollary}[theorem]{Corollary}

\theoremstyle{remark}
\newtheorem{remark}[theorem]{Remark}
\newtheorem{example}[theorem]{Example}

\title{Generalized Lam\'e equations with finite monodromy}
\author{You-Cheng Chou}
\address{Department of Mathematics, National Taiwan University, Taipei, Taiwan}
\email{b99201040@ntu.edu.tw}
\date{August 22, 2014}

\begin{document}
\maketitle
\begin{abstract}

In this paper, we study the algebraic form of the symmetric generalized Lam\'e equations which have finite projective monodromy groups. In particular, we consider equations with $3$ regular singular points on a flat torus $T$ which takes the form
\begin{equation*}
\begin{split}
\frac{d^2 y}{dz^2}-&\left[n_1(n_1+1)(\wp(z+a) +\wp(z-a))\right.\\
&\qquad\qquad\left. +A_1(\zeta(z+a) - \zeta(z-a))+n_0(n_0+1) \wp(z)+B\right]y=0,
\end{split}
\end{equation*}
where $n_1, n_0 \in \Bbb R$, $A_1,B \in \Bbb C$, and $\wp$ is the Weierstrass elliptic function. 

We give a complete list of all the group types that occur as the finite projective monodromy groups on the algebraic form and give the corresponding parameters $n_0$ and $n_1$. For equations with only $1$ or $2$ regular singular points, we further determine their monodromy group types. 

The main tool used is the Grothendieck correspondence which gives a bijection between Belyi pairs and dessin d'enfants. By Klein's theorem, we may regard generalized Lam\'e equations with finite monodromy as pullbacks of the hypergeometric equations. In this paper, we will restrict our cases with the assumption that the pullback maps are Belyi functions. Under this setup, our main results consist of a systematic construction on the required dessin. In particular a gluing procedure of dessin will be developed to enable inductive constructions of the dessin. 
\end{abstract}

\section{Introduction}

The generalized Lam\'e equation on a flat torus $T = \Bbb C/\Lambda$ with multiple regular singular points takes the following form
\begin{equation} \label{Wei-gene-Lame}
\frac{d^2 y}{dz^2}-\left[ \sum_{i=1}^{s} n_i(n_i+1)\wp(z-a_i)+\sum_{i=1}^{s}A_i\zeta(z-a_i)+B\right]y=0,
\end{equation}
where $\wp$ is the Weierstrass $\wp$ function associated to the lattice $\Lambda$, $\zeta$ is the Weierstrass $\zeta$ function, $n_i\in\mathbb{R}$, $A_i, B\in\mathbb{C}$ ($1 \le i \le s$) and $\sum_{i = 1}^s A_i=0$ for the periodic reason  (c.f.~\cite{Boss3}).

We will focus on the special case of equation (\ref{Wei-gene-Lame}) by assuming that the singular points have a center and choosing $A_i$ such that the equation is invariant under the transformation $z\rightarrow -z$. We get a symmetric form
\begin{equation} \label{C-Wei-Lame}
\begin{split}
&\frac{d^2 y}{dz^2}-\left[ \sum_{i=1}^{s}n_i(n_i+1)\left(\wp(z+a_i)+\wp(z-a_i)\right)\right.
\\
&\qquad\qquad\qquad
+\sum_{i=1}^{s}A_i(\zeta(z+a_i)-\zeta(z-a_i))+n_0(n_0+1)\wp(z)+B\Bigg]y=0. 
\end{split}
\end{equation}
In this special symmetric form, by changing the variable $x=\wp(z)$ which induces the double covering $T \to \Bbb P^1$, replacing the term $\wp(z\pm a)$ by the addition law for $\wp(z)$:
\[
\wp(z\pm a)=\frac{1}{4}\left\lbrace\frac{\wp'(z)\pm \wp'(a)}{\wp(z)-\wp(a)}\right\rbrace^{2}-\wp(z)-\wp(a),\]
and using the formula 
\[
\zeta(z-a) - \zeta(z+a) + 2\zeta(a)=\frac{\wp'(a)}{\wp(z) - \wp(a)},
\]
we get its algebraic form on $\Bbb P^1$,
\begin{equation} \label{C-Alg-Lame}
\begin{split}
 p(x)\frac{d^2 y}{dx^2} + \frac{1}{2}  p'(x)  \frac{dy}{dx} 
- \left\lbrace\sum_{i=1}^{s} n_i (n_i+1) \left(\frac{p(x) + \wp'(a_i)^{2}}{2(x - \wp(a_i))^2} - 2x - 2\wp(a_i)\right)\right. &
\\
 - \sum_{i=1}^{s} A_i\left.\left( \frac{\wp'(a_i)}{x - \wp(a_i)} - 2\zeta(a_i)\right) + n_0(n_0 + 1)x + B\right\rbrace y=0, &
\end{split}
\end{equation}
where $e_1,e_2,e_3\in \mathbb{C}$, $p(x)\equiv 4\prod{(x-e_i)}=4x^3-g_2 x-g_3$, and $e_i$'s are distinct. It is a second order Fuchsian equation with regular singular points at $\lbrace e_1,e_2,e_3,\wp(a_i),\infty\rbrace$. The caracteristic exponents are $\lbrace 0,1/2\rbrace$ at points $e_1,e_2,e_3$, $\lbrace -n_i,n_i+1\rbrace$ at points $\wp(a_i)$, and $\lbrace -n_0/2,(n_0+1)/2\rbrace$ at point $\infty$. For convenience, we say that equation (\ref{C-Alg-Lame}) has local exponent
$$\left\lbrace
\begin{array}{ccccccc}
e_1 & e_2 & e_3 & \wp(a_1)& \cdots &\wp(a_s) & \infty\\
0 & 0 & 0 & -n_1 &\cdots & -n_s & -n_0/2 \\
\frac{1}{2} & \frac{1}{2} & \frac{1}{2} & n_1+1 & \cdots & n_s+1 & (n_0+1)/2
\end{array}
\right\rbrace.$$

The main goal of this paper is to determine whether equation (\ref{C-Alg-Lame}) has finite (projective) monodromy group and to determine its group type. To study this problem, we first review the basic Schwarz list and Klein's theorem (\cite{BD}, theorem 3.4). We start by recalling the normal form of the hypergeometric differential equations: For $\lambda,\mu,\nu\in\mathbb{C}$, let\[ H_{\lambda,\mu,\nu}y \equiv \frac{d^2 y}{dx^2} + \left(\frac{1-\lambda^2}{4x^2}+\frac{1-\mu^2}{4(x-1)^2}+\frac{\lambda^2+\mu^2-\nu^2-1}{4x(x-1)}\right)y = 0.\] It is the unique second order normal differential equation with regular singular points only at $0,1,\infty$, with constant Wronskian, and with exponent difference $\lambda,\mu,\nu$ at $0,1,\infty$ respectively.\\
    
\begin{theorem}[Klein, c.f.~\cite{BD}]\label{Klein} If $L$ is a second order Fuchsian equation with finite but not cyclic projective monodromy group $G$, then it is the pullback by a rational map of one of the hypergeometric operators $H_{1/\lambda,1/\mu,1/\nu}$, where $\lambda, \mu, \nu$ are contained in the basic Schwarz list defined below. 

Conversely, such a pullback has a finite projective monodromy group. Furthermore, in this case, $G$ is the subgroup of $\widetilde{G}$, the projective monodromy group of the corresponding hypergeometric operators $H_{1/\lambda,1/\mu,1/\nu}$.
\end{theorem}

Here is the basic Schwarz list:\\\\
\begin{tabular}{c|ccc}
&$\displaystyle \frac{1}{\lambda}\quad\frac{1}{\mu}\quad\frac{1}{\nu}$ & group type & order of $\widetilde{G}$\\\\ \hline \\\uppercase\expandafter{\romannumeral 1} & $\displaystyle\frac{1}{2}\quad\frac{1}{2}\quad\frac{1}{n}$ & dihedral ($D_{n}$) & $2n$\\\\ \hline\\\uppercase\expandafter{\romannumeral 2} & $\displaystyle\frac{1}{2}\quad\frac{1}{3}\quad\frac{1}{3}$ & tetrahedral ($A_4$) & 12\\\\\hline\\\uppercase\expandafter{\romannumeral 3} & $\displaystyle\frac{1}{2}\quad \frac{1}{3}\quad \frac{1}{4}$ & octahedral ($S_4$) & 24\\\\\hline\\\uppercase\expandafter{\romannumeral 4} & $\displaystyle\frac{1}{2}\quad\frac{1}{3}\quad\frac{1}{5}$ & icosahedral ($A_5$) & 60\\
\end{tabular}\\\\

\begin{remark}
Let $Ly=0$ and $L'y=0$ be two second order Fuchsian differential equations on $\mathbb{P}^1$, and let $F:\mathbb{P}^1\rightarrow \mathbb{P}^1$ be a meromorphic function. We say that $Ly=0$ is a (weak) pullback of $L'y=0$ by $F$ if there exist independent solutions $y_1,y_2$ of $Ly=0$ and $y_1',y_2'$ of $L'y=0$ such that $y_1/y_2=(y_1'/y_2')\circ F$. Hence $\hat{L}$ is a (weak) pullback of $Ly$ by $F$ if and only if it has the same normal form as $L'$. In this paper, following the terminology in \cite{BD}, we use the term \emph{pullback} instead of \emph{weak pullback}.
\end{remark}

In this paper, we give all the possible $n_0$ and $n_1$ of equation (\ref{C-Alg-Lame}), with $s=1$, such that it has finite projective monodromy group. Conversely, for any given $n_0$ and $n_1$ in the list, we classify all the possible finite projective monodromy group that may occur in equation (\ref{C-Alg-Lame}). We divide this problem into three parts. The two special cases $n_0=0$ or $n_1=0$ will be firstly discussed in section 2 and section 3 respectively. The remaining third case has $n_0$ and $n_1$ being both nonzero. It is the most complicated case and will be solved in section 4. 

The main tool of the proof is the Grothendieck correspondence between Belyi pairs and dessin d'enfants (see section 2).  The idea to apply the Grothendieck correspondence to study the finite monodromy problem of ordinary differential equations was first introduced by Litcanu in \cite{Lit}. It describes how to  translate the problem on existence of the pullback functions in Theorem \ref{Klein} to the problem on existence of the dessin, if the pullback function is a Belyi function. This is the case that the pullback function maps all the regular-singular points of second order Fuchian equation $L$ in Theorem \ref{Klein} to $0,1$ or $\infty$, the regular-singular points of hypergeometric equation. Since, in this case, the pullback function only ramified at three points $0,1,\infty$ and then it is a Belyi function. This idea had caught more and more attention since then.  However, the construction of the dessin is still a tricky problem and so far not many systematic methods had been developed. Our paper is an attempt to contribute such a systematic construction in the context of generalized Lam\'e equations. 

It is important to note that the number of Belyi pairs up to equivalence is countable, since the number of dessin up to equivalence is. Furthermore, there is a bijection between the pullback function up to equivalence and equation (\ref{C-Alg-Lame}) up to scaling $(x\rightarrow ax)$. We will show that if $n_0 \notin\lbrace1/2\rbrace+\mathbb{Z}$ and $n_i\notin\mathbb{Z}/2$, $1\leq i\leq s$, the pullback function must be a Beyli function and then the number of equation (\ref{C-Alg-Lame}) with finite monodromy is countable up to scaling. In other words, if we fix $B=0$ or $1$ after scaling the equation, there are only countable many tori such that equation (\ref{C-Alg-Lame}) has finite monodromy. For the case, $n_0\in\lbrace1/2\rbrace+\mathbb{Z}$ or $n_i\in\mathbb{Z}/2$ for some $1\leq i\leq s$, we can not have the countable conclusion with previous argument, since the exponent difference at some regular-singular points of equation (\ref{C-Alg-Lame}) is integer and then the pullback function may not be a Belyi function.

We start with the classical case: For the special case $n_1, A_1=0$, (\ref{C-Alg-Lame}) reduces to the original Lam\'e equation
\begin{equation}
\frac{d^2y}{dz^2}-[n_0(n_0+1)\wp(z)+B]y
=0.
\end{equation}
Its algebraic form is given by
\begin{equation} \label{Alg-Lame}
L_{n_0,B,g_{2},g_{3}}y\equiv\frac{d^{2}y}{dx^{2}}+\frac{1}{2}\sum_{i=1}^{3}\frac{1}{x-e_{i}}\frac{dy}{dx}-\frac{n_0(n_0+1)x+B}{p(x)}y=0,
\end{equation} 
which has local exponent 
$$\left\lbrace
\begin{array}{cccc}
e_1 & e_2 & e_3  & \infty\\
0 & 0 & 0 &  -n_0/2 \\
\frac{1}{2} & \frac{1}{2} & \frac{1}{2} & (n_0+1)/2
\end{array}
\right\rbrace.$$

There is a well-known result that the Lam\'e equation never has cyclic or $A_4$ as its finite projective monodromy group. We refer to (\cite{BD}, Theorem 4.1) or (\cite{B2}, Proposition 1.1, Proposition 3.1) for the details.

We can further consider the monodromy group $M\subset $GL$(2,\mathbb{C})$ of Lam\'e equation $L_{n_0,B,g_2,g_3}y=0$. Let $\gamma_i$ be the local monodromy matrices at $e_i$, and $\gamma_{\infty}$ be the local monodromy matrix at $\infty$. We have the equation $\gamma_1^2=\gamma_2^2=\gamma_3^2={I}_2$ and $\gamma_1\gamma_2\gamma_3\gamma_{\infty}={I}_2$. Hence, we get the conclusion that $M$ is generated by the reflections $\gamma_1,\gamma_2$ and $\gamma_3$ and therefore is a \emph{complex reflection group}. Using the Shephard-Todd classification we can list all possible finite complex reflection groups that could be the monodromy group of the Lam\'e equation. We get further restrictions on the cases by considering rank $2$ groups generated by order $2$ reflections. These are 5 possible types:

\begin{enumerate}

\item 
$G(N,N/2,2)$ $(N\in2\mathbb{Z}_{\geq 0})$ is the group of order $4N$ generated by
\[
\left(\begin{array}{cc}e^{2\pi i/N} & 0 \\0 & e^{-2\pi i/N}\end{array}\right), \quad\left(\begin{array}{cc} -1 & 0 \\0 & 1 \end{array}\right) , \quad \left(\begin{array}{cc} 0 & 1 \\1 & 0\end{array}\right).
\]
Its projective group is $D_{N/2} $ if $4\mid N$, and is $D_{N}$ otherwise.

\item 
$G(N,N,2)$ $(N\geq 3)$ is the group of order $2N$ generated by 
\[
\left(\begin{array}{cc}e^{2\pi i/N} & 0 \\0 & e^{-2\pi i/N}\end{array}\right) , \quad\left(\begin{array}{cc} 0 & 1 \\1 & 0\end{array}\right).
\]Its projective group is $D_{N/2}$ if $N$ is even and $D_{N}$ otherwise.

\item 
$G_{12}$ is the group of order $48$ generated by 
\[
\frac{1}{\sqrt{2}}\left(\begin{array}{cc} 0 & 1+i \\1-i & 0\end{array}\right), \quad \frac{1}{\sqrt{2}}\left(\begin{array}{cc}1 & 1 \\1 & -1\end{array}\right), \quad \frac{1}{\sqrt{2}}\left(\begin{array}{cc} 1 & i \\-i & -i\end{array}\right). 
\]Its projective group is $S_{4}$.

\item 
$G_{13}$ is the group of order $96$. It is generated by the elements of $G_{12}$ and $i\cdot {I}_2$. Its projective group is $S_{4}$.

\item 
$G_{22}$ is the group of order 120 generated by 
\[
\left(\begin{array}{cc} i & 0 \\0 & i\end{array}\right), \quad 
\frac{1}{\sqrt{5}}\left(\begin{array}{cc} \zeta_5-\zeta_5^4 & \zeta_5^2-\zeta_5^3 \\\zeta_{5}^2-\zeta_5^3 & \zeta_{5}^{4}-\zeta_{5}\end{array}\right), \quad 
\frac{1}{\sqrt{5}}\left(\begin{array}{cc} \zeta_5^3-\zeta_5 & 1-\zeta_5 \\\zeta_5^4-1 & \zeta_5^2-\zeta_5^4\end{array}\right).
\] Its projective group is $A_{5}$.

\end{enumerate}

Beukers and Waall had listed all the possible $n_0$ in $L_{n_0,B,g_2,g_3}y=0$ with given finite monodromy group in each of the above $5$ types: 

\begin{theorem}[\cite{BA}, Theorem 4.4]\label{BD,Intro}

Let $M$ be the monodromy group of the Lam\'e equation $L_{n_0,B,g_2,g_3}y=0$.
\begin{enumerate}
\item 
$M\cong G(N,N/2,2) \Rightarrow n_0\in\lbrace 1/2\rbrace+\mathbb{Z}$ and $N=4$;
\item 
$M\cong G(N,N,2),\ N\geq 3\Rightarrow n_0\in\mathbb{Z}$ and $N\neq 4$;
\item 
$M\cong G_{12} \Rightarrow n_0\in \lbrace \pm1/4 \rbrace+\mathbb{Z}$;
\item 
$M\cong G_{13}\Rightarrow n_0 \in \lbrace \pm1/6 \rbrace+\mathbb{Z}$;
\item 
$M\cong G_{22}\Rightarrow n_0\in \lbrace \pm1/10,\pm3/10,\pm1/6\rbrace+\mathbb{Z}$.
\end{enumerate}
Moreover, if $M$ is finite and $n\geq -1/2$, then $n_0>0$.

\end{theorem}

There are some partial results for the converse of Theorem \ref{BD,Intro}:

\begin{enumerate}
\item 
For $n_0 = k+\frac{1}{2},k\in\mathbb{Z}_{\geq 0}$, there exists a polynomial $p_{n_0}(B,g_2,g_3)$ of degree $k+1$ in $B$ such that the Lam\'e equation $L_{n_0,B,g_2,g_3}y=0$ has finite monodromy if and only if $p_{n_0}(B,g_2,g_3)=0$. $p_{n_0}$ is known as the Brioschi-Halphen determinant. We refer to (\cite{Po}, p.164) for a proof.

\item 
For $n_0\in\mathbb{N}$, this question is still unsolved. It is related to the modular form problem. See \cite{BD}, \cite{D}, \cite{D1} and \cite{Boss2} for recent progresses.
\end{enumerate}

In section 2, we give the converse side of Theorem \ref{BD,Intro} for all the remaining cases. That is, for $M\in\lbrace G_{12},G_{13},G_{22}\rbrace$ with the corresponding value $n$ in Beukers and Wall's list, there exists $g_2,g_3,B$ such that $L_{n_0,B,g_2,g_3}y=0$ has $M$ as its monodromy group. We first prove the theorem on the level of projective monodromy group.
\begin{theorem}[$=$ Theorem \ref{BM}, c.f.~\cite{M}, \cite{N}]\label{BM,Intro}
Let $PM$ be the projective monodromy group of $L_{n,B,g_2,g_3}$. Then, for each of the following pairs of PM, n there exists a Lam\'e equation:
\begin{enumerate}
\item $PM\cong S_4$ with $n\in \lbrace \pm1/4,\pm1/6\rbrace+\mathbb{Z}$;
\item $PM\cong A_5$ with $ n\in \lbrace \pm1/10,\pm3/10,\pm1/6\rbrace+\mathbb{Z}$.
\end{enumerate}
\end{theorem}

 After I finished the proof, I found that Maier and Nakanishi had achieved the same conclusions. Maier \cite{M} wrote down the pullback function explicitly for some special types of tori. Nakanishi \cite{N} also used the Grothendieck correspondence and constructed the corresponding dessin. Indeed, we prove the same conclusions but with different constructions of the dessin. To make the construction easier (more canonical) and \emph{suitable for generalizations to equations with more regular singular points}, we leave the point with large ramification index in the outer face (c.f.~section 2 for the terminologies). Our method will also be used in later sections. Moreover, we can further determine its monodromy group by combining with the fact that if a second order differential equation has algebraic Wronskian then it has finite projective monodroy group if and only if it has finite monodromy group. This leads to

\begin{corollary} [$=$ Corollary \ref{C1}] \label{C1,Intro}
Let M be the monodromy group of $L_{n,B,g_2,g_3}y = 0$. Then, for each of the following pairs of M, n there exists a Lam\'e equation:
\begin{enumerate}
\item 
$M\cong G_{12}$ with $n\in \lbrace \pm1/4\rbrace+\mathbb{Z}$;
\item 
$M\cong G_{13}$ with $n\in \lbrace \pm1/6\rbrace+\mathbb{Z}$;
\item 
$M\cong G_{22}$ with $n\in \lbrace \pm1/10,\pm 3/10,\pm1/6\rbrace+\mathbb{Z}$.
\end{enumerate}
\end{corollary}

Now we describe our new results for the generalized Lam\'e equations: For the case $n_0=0$, the equation is given by
\begin{equation} \label{Wei-2sing-Lame}
\begin{split}
\frac{d^{2}y}{dz^{2}}-[n_1(n_1+1)(\wp(z+a)+&\wp(z-a))\\
&+A_1(\zeta(z+a)-\zeta(z-a))+B]y=0.
\end{split}
\end{equation} 
It is the simplest generalization for the original Lam\'e equation. Its algebraic form is 
\begin{equation} \label{Alg-2sing-Lame}
\begin{split}
 p(x)\frac{d^2 y}{dx^2} + \frac{1}{2}p'(x)\frac{dy}{dx} 
 - \left\lbrace n_1(n_1+1) \left(\frac{p(x)+\wp'(a)^{2}}{2(x-\wp(a))^2}-2x-2\wp(a)\right) \right. &
 \\ 
\left. - A_1\left( \frac{\wp'(a)}{x - \wp(a)} - 2\zeta(a)\right) + B\right\rbrace y=0, &
\end{split}
\end{equation}
with local exponent
$$\left\lbrace
\begin{array}{ccccc}
e_1 & e_2 & e_3 & \wp(a)  & \infty\\
0 & 0 & 0 & -n_1 & 0 \\
\frac{1}{2} & \frac{1}{2} & \frac{1}{2} & n_1+1 & \frac{1}{2}
\end{array}
\right\rbrace.$$
Since four of its singularities have exponent difference equals to $1/2$, we also conclude that the monodromy group of equation (\ref{Alg-2sing-Lame}) is a complex reflection group. Hence, the Shepperd-Todd classification can also be applied to list all the possible finite monodromy group in this case. It turns out that the resulting cases are the same as the original Lam\'e equation. We also give all the possible $ n_1\notin\mathbb{Z}/2$ for the given monodromy group. 

\begin{theorem} [$=$ Theorem \ref{BM1}]\label{BM1,Intro}
Let M be the monodromy group of equation (\ref{Alg-2sing-Lame}) and $n \notin\mathbb{N}/2$. Then 
\begin{enumerate}
\item 
$M \cong G(N,N/2,2)\Rightarrow n\in\lbrace\pm 1/4\rbrace+\mathbb{Z}$, $N$ is even and $N\geq 4$;
\item
 $M \cong G_{12} \Rightarrow n\in\lbrace \pm 1/6,\pm 1/4,\pm1/3\rbrace+\mathbb{Z}$;
\item 
$M\cong G_{13}\Rightarrow n\in\lbrace\pm 1/8,\pm3/8,\pm 1/6,\pm 1/4,\pm 1/3\rbrace +\mathbb{Z}$ ;
\item 
$M\cong G_{22}\Rightarrow n\in(\mathbb{Z}/6)\cup(\mathbb{Z}/10) \setminus(\mathbb{Z}/2) $.
\end{enumerate}
Moreover, if $n\geq -1/2$, then $n\geq 1/6$.
\end{theorem}

Furthermore, we prove the converse side of it. Namely for given $n_1$, $M$ in Theorem \ref{BM1,Intro}, there exist $g_2,g_3,\wp(a),A_1,B$ such that equation (\ref{Alg-2sing-Lame}) has $M$ as its monodromy group.

\begin{theorem} [$=$ Theorem \ref{BM2}] \label{BM2,Intro}
Let M be the projective monodromy group of the generalized Lam\'e equation. Then, for each of the following pairs of M, n there exists a generalized Lam\'e equation (\ref{Alg-2sing-Lame}):
\begin{enumerate}
\item 
$M\cong G_{12}$ with $n\in\lbrace \pm1/6,\pm1/4,\pm1/3\rbrace+\mathbb{Z}$;
\item 
$M\cong G_{13}$ with $n\in\lbrace\pm1/8,\pm3/8,\pm1/6,\pm1/4,\pm1/3\rbrace +\mathbb{Z}$;
\item 
$M\cong G_{22}$ with $n\in(\mathbb{Z}/6)\cup(\mathbb{Z}/10) \setminus(\mathbb{Z}/2) $.
\end{enumerate}
\end{theorem}

We remark that, as in the case of classical Lam\'e equations, the study on the exceptional case $n \in \Bbb N/2$ can also be linked to certain modular form problems. This will be discussed in a forthcoming work \cite{CW}.

For the case that $n_0$ and $n_1$ are both non-zero, the equation is given by 
\begin{equation} \label{C-ale-3pts}
\begin{split}
p(x)\frac{d^2 y}{dx^2}+\frac{1}{2}p'(x){\frac{dy}{dx}}
-\left\lbrace n_1(n_1+1)\left(\frac{p(x)+\wp'(a)^{2}}{2(x-\wp(a))^{2}}-2x-2\wp(a)\right) \right. &
\\
\left. - A_1\left( \frac{\wp'(a)}{x - \wp(a)} - 2\zeta(a)\right) + n_0(n_0+1)x + B\right\rbrace y=0, &
\end{split}
\end{equation}
which has local exponent
$$\left\lbrace
\begin{array}{ccccc}
e_1 & e_2 & e_3 & \wp(a)  & \infty\\
0 & 0 & 0 & -n_1 & -n_0/2 \\
\frac{1}{2} & \frac{1}{2} & \frac{1}{2} & n_1+1 & (n_0+1)/2
\end{array}
\right\rbrace.$$
Notice that the local monodromies around $\wp(a),\infty$ are in general not reflection elements, hence the monodromy group of equation (\ref{C-Alg-Lame}) may not be a complex reflection group. Nevertheless, we can still list all the possible group types that may occur as the projective monodromy group $PM$. 

All the group types given in Klein's theorem turn out are all possible. Furthermore, we can generalize the method used in original Lam\'e equation to exclude the cyclic case (c.f.~Theorem \ref{A-Lame-Cyc}). Then we will give all the possible $n_0$ and $n_1$ for the projective monodromy group being in the basic Schwarz list (c.f. Theorem \ref{n0n1} and Theorem \ref{n0n1-3}). Furthermore, we prove the converse side of it: For any given $n_0, n_1$ and $PM$ as described above, there exists equation (\ref{C-ale-3pts}) with $PM$ as its projective monodromy group (Theorem \ref{n0n1} and Theorem \ref{3-pts M}). 

To conclude the introduction, we remark that to construct the dessin we make use of the "gluing method" extensively. Our gluing method gives a systematic way to construct the dessin with large degree points and faces (see the proofs of case (2), case (7) and case (14) of Theorem \ref{3-pts M}). It describes how to construct the white points, faces and black points of the dessin with large degrees respectively. This gluing method can also be generalized to construct the dessin which corresponds to the generalized Lam\'e equation (\ref{C-Alg-Lame}) with $2n+1$ regular singular points for $n \ge 2$. However, so far we obtain only partial results in this way and we plan to give a complete treatment of it in a later work (see Remark \ref{R. Gene}).

\section{Existence of Lam\'e equations with $G_{12}$, $G_{13}$, $G_{22}$ as the monodromy group}

In this section, we only consider the case that the Lam\'e equation has $S_4$ or $A_5$ as its projective monodromy group. The following theorem due to Beukers and Waall gives strong restrictions to the cases of $n$.
\begin{theorem}[\cite{BA}, Theorem 4.4]\label{BD}

Let $M$ be the monodromy group of the Lam\'e equation $L_{n,B,g_2,g_3}y=0$.
\begin{enumerate}
\item 
$M\cong G(N,N/2,2) \Rightarrow n\in\lbrace 1/2\rbrace+\mathbb{Z}$ and $N=4$;
\item 
$M\cong G(N,N,2),\ N\geq 3\Rightarrow n\in\mathbb{Z}$ and $N\neq 4$;
\item 
$M\cong G_{12} \Rightarrow n\in \lbrace \pm1/4 \rbrace+\mathbb{Z}$;
\item 
$M\cong G_{13}\Rightarrow n \in \lbrace \pm1/6 \rbrace+\mathbb{Z}$;
\item 
$M\cong G_{22}\Rightarrow n\in \lbrace \pm1/10,\pm3/10,\pm1/6\rbrace+\mathbb{Z}$.
\end{enumerate}
Moreover, if $M$ is finite and $n\geq -1/2$, then $n>0$.

\end{theorem}

The main work of this section gives the converse side of Theorem \ref{BD} for the case $M\in\lbrace G_{12},G_{13},G_{22}\rbrace$. In other words, for given $n$ in the Theorem \ref{BD} corresponding to $M$ there exists $B,g_2,g_3$ such that the equation $L_{n,B,g_2,g_3}y=0$ has $M$ as its monodromy group. We will highly use the Grothendieck correspondence of the Belyi pair and dessin d'enfants, or simply dessin. 

Here we recall the contents of the Grothendieck correspondence and explain how we use it. If $X$ is a compact Riemann surface and $F:X\rightarrow \mathbb{P}^1$ is a covering unramified over $\mathbb{P}^1 \setminus\lbrace 0,1,\infty\rbrace$ then $(X,F)$ is called a Belyi pairs. Two Belyi pairs $(X,F)$ and $(X',F')$ are equivalent if there exists an isomorphism $G:X\rightarrow X'$ such that $F=F'\circ G$. (It is well known that $X$, as an algebraic curve, is defined over $\bar {\Bbb Q}$ if and only if such an $F$ exists. We will not make any explicit usage of this important fact in this paper.)

A graph which is connected, bi-colored (every vertex endowed with white or black color, and every edges connects a white edges with a black edges), and oriented (at every edges a cyclic ordering of the edges contacting with that vertex is given) is called a dessin d'enfants. Two dessins are equivalent if there exists a graph isomorphism between them which preserves the bi-coloring and ordering. 

The Grothendieck correspondence states that there is a bijection between the Belyi pairs modulo equivalence and dessins modulo equivalence. In one direction the graph is given by $F^{-1}([0, 1])$. (Alternatively one may use $[0, \infty]$ or $[1, \infty]$ due to the equivalence.) The two set of vertexes $F^{-1}(0)$ and $F^{-1}(1)$ are labeled by different colors and the orientation is induced from the one on $X$. For the converse side, $(X, F)$ can be constructed through the Riemann existence theorem. In this paper, we only need the case $X=\mathbb{P}^1$ and $F$ is the pullback function in Theorem \ref{Klein}.

\begin{theorem}[\cite{M}, \cite{N}]\label{BM}
Let $PM$ be the projective monodromy group of $L_{n,B,g_2,g_3}$. Then, for each of the following pairs of PM, n there exists a Lam\'e equation:
\begin{enumerate}
\item $PM\cong S_4$ with $n\in \lbrace \pm1/4,\pm1/6\rbrace+\mathbb{Z}$;
\item $PM\cong A_5$ with $ n\in \lbrace \pm1/10,\pm3/10,\pm1/6\rbrace+\mathbb{Z}$.
\end{enumerate}
\end{theorem}

To prove the theorem, it suffices to construct the corresponding dessin. We will use constructions different from the ones used in \cite{N}. Here we give an example and meanwhile fix our notations.

\begin{example}[$M=S_4$, $n= 11/6$] 
We define a convenient notation first. We say the list $[a_1^{p_1}\cdots a_l^{p_l},b_1^{q_1}\cdots b_m^{q_m},c_1^{r_1}\cdots c_n^{r_n} ]$ corresponding to a Belyi function $G$ means the function $G$ has $p_i$ points which maps to $0$ of degree $a_i$, $q_j$ points which maps to $1$ of degree $b_j$ and $r_k$ points which maps to $\infty$ of degree $c_k$. Now, we consider the Lam\'e equation 
\[
L_{11/6,B,g_2,g_3}y\equiv\frac{d^2y}{dx^2}+\frac{1}{2}\sum_{i=1}^{3}\frac{1}{x-e_i}\frac{dy}{dx}-\frac{(187/36)x+B}{p(x)}y=0,
\]
According to the Klein's theorem, if it has $S_4$ as its projective monodromy group, then it is given by the pull-back of $H_{1/2,1/3,1/4}$. The pull-back map is the Belyi function $ F:\mathbb{P^1}\rightarrow \mathbb{P^1}$ with degree$=22$. To compute the degree of $F$, we refer to (\cite{BD}, Lemma 1.5). By checking the local exponent carefully, we get the 2 possible tables:
\begin{table}[h]
\begin{tabular}{|l|l|l|l|l|l|l|}
\hline
 & $e_1$ & $e_2$ & $e_3$ & $\infty$& \\ \hline
0 & 1 & 1 &  &  & 10 points with multiplicity 2\\ \hline
1 &  &  &  & 7 & 5 points with multiplicity 3\\ \hline
$\infty$ &  &  & 2 &  & 5 points with multiplicity 4\\ \hline
\end{tabular}
\end{table}

It follows that the Belyi function $F$ corresponds to list $[1^2 2^{10},3^5 7^1,2^1 4^5]$. From now on, we will just use this notation without drawing the table. According to the Grothendieck correspondence, the existence of the Belyi pair is equivalence to the existence of the corresponding dessin. In this case, we just draw the first list to give the existence. The dessin is given by
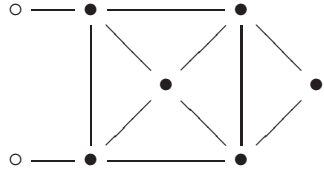
\begin{displaymath}
\xygraph{
!{<0cm,0cm>;<1cm,0cm>:<0cm,1cm>::}
!{(0,0)}*+{\bullet}="1"
!{(1,1)}*+{\bullet}="2"
!{(1,-1)}*+{\bullet}="3"
!{(2,0)}*+{\bullet}="4"
!{(-1,1)}*+{\bullet}="5"
!{(-1,-1)}*+{\bullet}="6"
!{(-2,1)}*+{\circ}="7"
!{(-2,-1)}*+{\circ}="8"
"1"-"2" "1"-"3" "1"-"6" "1"-"5"
"2"-"3" "2"-"4" "2"-"5"
"3"-"4" "3"-"6"
"5"-"7" "6"-"8" "6"-"5"
}
\end{displaymath}
where $\bullet$ represent the points mapped to 1, $\circ$ represent the points mapped to 0. Furthermore, we draw $\bullet - \bullet$ instead of $\bullet -\circ -\bullet$ if no confusion may exist.

\end{example}

\begin{proof}

Since $L_{n,B,g_2,g_3}=L_{-n-1,B,g_2,g_3}$, we assume $n\geq -1/2$. Furthermore, by theorem \ref{BM}, we only need to consider the case $n>0$.

In the proof, we will not give all the possible list. For eash case, we only give the list which has a corresponding dessin. Furthermore, we draw $\bullet - \bullet$  to represent $\bullet -\circ -\bullet$.

\begin{enumerate}[label=\textbf{case \arabic* :}]

\item
($M=S_4,\ n\in\lbrace \pm 1/4\rbrace+\mathbb{Z},\ n>0$)

In this case, $\bullet$ represent the points mapped to $1$; $\circ$ represent the points mapped to $0.$

For $n=(1+2k)/4$, $k\in\mathbb{N}$, the corresponding list is $[1^{3}2^{3k},3^{2k+1},\\(2k+3)^{1}4^{k}]$, and the dessin is given by
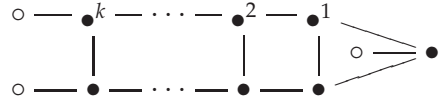
\begin{displaymath}
\xygraph{
!{<0cm,0cm>;<1cm,0cm>:<0cm,1cm>::}
!{(1,0)}*+{\bullet^{2}}="a1"
!{(2,0)}*+{\bullet^{1}}="a2"
!{(3.5,-0.5)}*+{\bullet}="a3"
!{(2,-1)}*+{\bullet}="a4"
!{(1,-1)}*+{\bullet}="a5"
!{(-1,0)}*+{\bullet^{k}}="a6"
!{(-1,-1)}*+{\bullet}="a7"
!{(0,0)}*+{\cdots}="c1"
!{(0,-1)}*+{\cdots}="c2"
!{(2.5,-0.5)}*+{\circ}="b1"
!{(-2,0)}*+{\circ}="b2"
!{(-2,-1)}*+{\circ}="b3"
"c1"-"a1" "c1"-"a6" "c2"-"a5" "c2"-"a7"
"a1"-"a5" "a1"-"a2" "a2"-"a4" "a4"-"a5"
"a2"-"a3" "a3"-"a4" "a3"-"b1"
"a6"-"a7" "a6"-"b2" "a7"-"b3"
}
\end{displaymath}\\
For the case $n=1/4$, the existence is clear.

\item 
($M=S_4,\ n\in\lbrace\pm1/6\rbrace+\mathbb{Z},\ n>0$)

In this case, $\bullet$ represent the points mapped to $\infty$; $\circ$ represent the points mapped to $0.$

For $n=(5+2k)/6$, $k\in\mathbb{N}$, the corresponding list is $[1^{2}2^{2k+4},\\(k+4)^{1}3^{k+2},2^{1}4^{k+2}]$. If $k$ is even, the dessin is given by
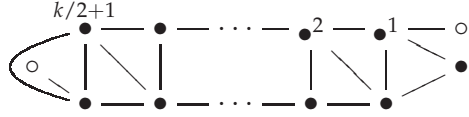
\begin{displaymath}
\xygraph{
!{<0cm,0cm>;<1cm,0cm>:<0cm,1cm>::}
!{(1,0)}*+{\bullet^{2}}="a1"
!{(2,0)}*+{\bullet^{1}}="a2"
!{(3,-0.5)}*+{\bullet}="a3"
!{(2,-1)}*+{\bullet}="a4"
!{(1,-1)}*+{\bullet}="a5"
!{(-1,0)}*+{\bullet}="a6"
!{(-2,0)}*+{\bullet}="a7"
!{(-2,0.2)}*+{{}^{k/2+1}}
!{(0,0)}*+{\cdots}="c1"
!{(0,-1)}*+{\cdots}="c2"
!{(-2,-1)}*+{\bullet}="a9"
!{(-1,-1)}*+{\bullet}="a10"
!{(3,0)}*+{\circ}="b1"
!{(-2.7,-0.5)}*+{\circ}="b2"
"c1"-"a1" "c1"-"a6" "c2"-"a5" "c2"-"a10"
"a1"-"a2" "a5"-"a4" "a1"-"a4" "a1"-"a5"
"a2"-"a4" "a2"-"a3" "a3"-"a4"
"a6"-"a10" "a6"-"a7" "a9"-"a10" "a7"-"a10" 
"a7"-"a9" "a2"-"b1" "a9"-"b2"
"a7"-@/_1cm/"a9"
}
\end{displaymath}\\
If $k$ is odd, the dessin is given by 
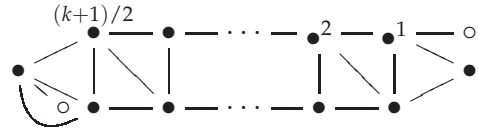
\begin{displaymath}
\xygraph{
!{<0cm,0cm>;<1cm,0cm>:<0cm,1cm>::}
!{(1,0)}*+{\bullet^{2}}="a1"
!{(2,0)}*+{\bullet^{1}}="a2"
!{(3,-0.5)}*+{\bullet}="a3"
!{(2,-1)}*+{\bullet}="a4"
!{(1,-1)}*+{\bullet}="a5"
!{(-1,0)}*+{\bullet}="a6"
!{(-2,0)}*+{\bullet}="a7"
!{(-2,0.2)}*+{{}^{(k+1)/2}}
!{(0,0)}*+{\cdots}="c1"
!{(0,-1)}*+{\cdots}="c2"
!{(-2,-1)}*+{\bullet}="a9"
!{(-1,-1)}*+{\bullet}="a10"
!{(3,0)}*+{\circ}="b1"
!{(-3,-0.5)}*+{\bullet}="a8"
!{(-2.4,-1)}*+{\circ}="b2"
"c1"-"a1" "c1"-"a6" "c2"-"a5" "c2"-"a10"
"a1"-"a2" "a5"-"a4" "a1"-"a4" "a1"-"a5"
"a2"-"a4" "a2"-"a3" "a3"-"a4"
"a6"-"a10" "a6"-"a7" "a9"-"a10" "a7"-"a10" 
"a7"-"a9" "a2"-"b1" "a7"-"a8" "a8"-"a9"
"a8"-"b2" "a8"-@/_0.5cm/"a9"
}
\end{displaymath}\\
The existence of the case $n=1/6,5/6$ is easy.

\item 
($M=A_5,\ n\in\lbrace\pm1/10\pm3/10\rbrace+\mathbb{Z},\ n>0$)

In this case, $\bullet$ represent the points mapped to $1$; $\circ$ represent the points mapped to $0.$

For $n=(7+2k)/10$, $ k\in\mathbb{N}$, the corresponding list is $[1^3 2^{9+3k},3^{7+2k},\\(6+k)^1 5^{3+k}]$. If $k$ is an odd number, the dessin is given by
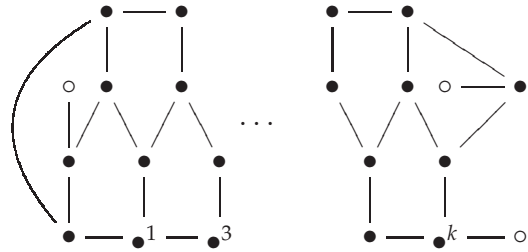
\begin{displaymath}
\xygraph{
!{<0cm,0cm>;<1cm,0cm>:<0cm,1cm>::}
!{(0,-0.5)}*+{\cdots}="c0"
!{(1,0)}*+{\bullet}="a1"
!{(2,0)}*+{\bullet}="a2"
!{(1,1)}*+{\bullet}="a3"
!{(2,1)}*+{\bullet}="a4"
!{(1.5,-1)}*+{\bullet}="a5"
!{(2.5,-1)}*+{\bullet}="a6"
!{(1.5,-2)}*+{\bullet}="a7"
!{(2.5,-2)}*+{\bullet^k}="a8"
!{(-1,0)}*+{\bullet}="a9"
!{(-1,1)}*+{\bullet}="a11"
!{(-2,0)}*+{\bullet}="a10"
!{(-2,1)}*+{\bullet}="a12"
!{(-1.5,-1)}*+{\bullet}="a13"
!{(-2.5,-1)}*+{\bullet}="a14"
!{(-1.5,-2)}*+{\bullet^1}="a15"
!{(-2.5,-2)}*+{\bullet}="a16"
!{(-2.5,0)}*+{\circ}="b1"
!{(3.5,0)}*+{\bullet}="a17"
!{(2.5,0)}*+{\circ}="b2"
!{(3.5,-2)}*+{\circ}="b3"
!{(-0.5,-1)}*+{\bullet}="a18"
!{(-0.5,-2)}*+{\bullet^3}="a19"
"a1"-"a3" "a3"-"a4" "a4"-"a2" "a1"-"a5" "a2"-"a5" "a5"-"a7"
"a7"-"a8" "a8"-"a6" "a6"-"a2" "a15"-"a16" "a14"-"a16" "a10"-"a14"
"a9"-"a11" "a11"-"a12" "a12"-"a10" "a9"-"a13" "a10"-"a13" "a13"-"a15"
"a12"-@/_1cm/"a16" "a14"-"b1" "a4"-"a17" "a6"-"a17" "a17"-"b2"
"a8"-"b3" "a9"-"a18" "a18"-"a19" "a19"-"a15"
}
\end{displaymath}\\
If $k$ is an even number, the dessin is given by
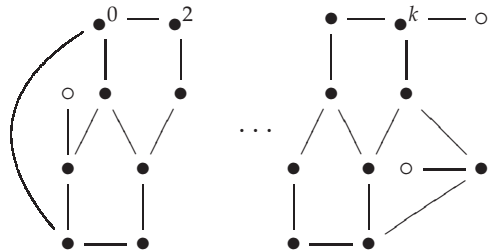
\begin{displaymath}
\xygraph{
!{<0cm,0cm>;<1cm,0cm>:<0cm,1cm>::}
!{(0,-0.5)}*+{\cdots}="c0"
!{(1,0)}*+{\bullet}="a1"
!{(2,0)}*+{\bullet}="a2"
!{(1,1)}*+{\bullet}="a3"
!{(2,1)}*+{\bullet^k}="a4"
!{(3,1)}*+{\circ}="b3"
!{(1.5,-1)}*+{\bullet}="a5"
!{(2,-1)}*+{\circ}="b4"
!{(3,-1)}*+{\bullet}="a20"
!{(1.5,-2)}*+{\bullet}="a7"
!{(-1,0)}*+{\bullet}="a9"
!{(-1,1)}*+{\bullet^2}="a11"
!{(-2,0)}*+{\bullet}="a10"
!{(-2,1)}*+{\bullet^0}="a12"
!{(-1.5,-1)}*+{\bullet}="a13"
!{(-2.5,-1)}*+{\bullet}="a14"
!{(-1.5,-2)}*+{\bullet}="a15"
!{(-2.5,-2)}*+{\bullet}="a16"
!{(-2.5,0)}*+{\circ}="b1"
!{(0.5,-1)}*+{\bullet}="a18"
!{(0.5,-2)}*+{\bullet}="a19"
"a1"-"a3" "a3"-"a4" "a4"-"a2" "a1"-"a5" "a2"-"a5" "a5"-"a7"
 "a15"-"a16" "a14"-"a16" "a10"-"a14"
"a9"-"a11" "a11"-"a12" "a12"-"a10" "a9"-"a13" "a10"-"a13" "a13"-"a15"
"a12"-@/_1cm/"a16" "a14"-"b1" 
 "a1"-"a18" "a18"-"a19" "a19"-"a7" "a4"-"b3"
 "a2"-"a20" "a7"-"a20" "a20"-"b4"
}
\end{displaymath}\\
The existence of the case $1/10,3/10,7/10$ is clear.

\item 
($M=A_5 ,\ n\in\lbrace \pm 1/6\rbrace+\mathbb{Z},\ n>0$)

It is the most complex case. The proof of this case is basic on "repeat" the following picture many times. The picture is given by
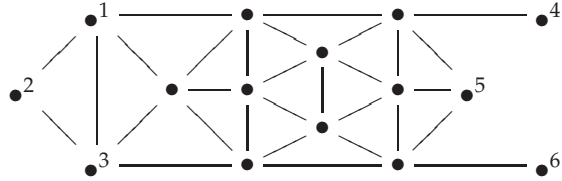
\begin{displaymath}
\xygraph{
!{<0cm,0cm>;<1cm,0cm>:<0cm,1cm>::}
!{(0,0)}*+{\bullet}="a1"
!{(0,1)}*+{\bullet}="a2"
!{(0,-1)}*+{\bullet}="a3"
!{(1,0.5)}*+{\bullet}="a4"
!{(1,-0.5)}*+{\bullet}="a5"
!{(2,0)}*+{\bullet}="a6"
!{(2,1)}*+{\bullet}="a7"
!{(2,-1)}*+{\bullet}="a8"
!{(-1,0)}*+{\bullet}="a9"
!{(-2,1)}*+{\bullet^1}="c1"
!{(-2,-1)}*+{\bullet^3}="c2"
!{(-3,0)}*+{\bullet^2}="c3"
!{(4,1)}*+{\bullet^4}="c4"
!{(3,0)}*+{\bullet^5}="c5"
!{(4,-1)}*+{\bullet^6}="c6"
"c4"-"a7" "c5"-"a7" "c5"-"a6" "c5"-"a8" "c6"-"a8" "a7"-"a2"
"a7"-"a4" "a6"-"a4" "a6"-"a5" "a8"-"a5" "a8"-"a3" "a4"-"a2"
"a4"-"a1" "a4"-"a5" "a5"-"a1" "a5"-"a3" "a6"-"a7" "a6"-"a8"
"a1"-"a2" "a1"-"a3" "a9"-"a2" "a9"-"a1" "a9"-"a3" "c1"-"a9"
"c2"-"a9" "c1"-"c2" "c1"-"c3" "c2"-"c3" "c1"-"a2" "c2"-"a3"
}
\end{displaymath}
Repeat the picture means "glue" the $\bullet^1\bullet^4$, $\bullet^2\bullet^5$ and $\bullet^3\bullet^6$ together.

From now on, we use the notation
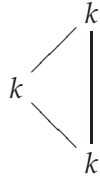
\begin{displaymath}
\xygraph{
!{<0cm,0cm>;<1cm,0cm>:<0cm,1cm>::}
!{(0,0)}*+{k}="a1"
!{(1,1)}*+{k}="a2"
!{(1,-1)}*+{k}="a3"
"a1"-"a2" "a3"-"a1" "a2"-"a3"
}
\end{displaymath}
to denote repeating the picture $k$ times after it.

For example, we draw 
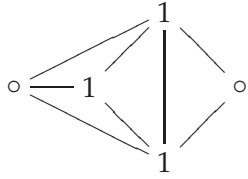
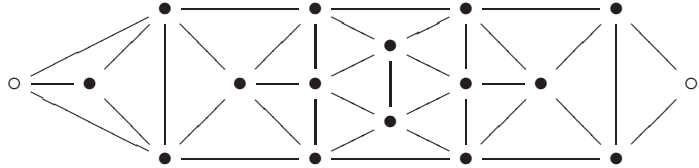
\begin{displaymath}
\xygraph{
!{<0cm,0cm>;<1cm,0cm>:<0cm,1cm>::}
!{(0,0)}*+{1}="a1"
!{(1,1)}*+{1}="a2"
!{(1,-1)}*+{1}="a3"
!{(2,0)}*+{\circ}="a4"
!{(-1,0)}*+{\circ}="a5"
"a1"-"a2" "a3"-"a1" "a2"-"a3" "a2"-"a4" "a3"-"a4"
"a1"-"a5" "a2"-"a5" "a3"-"a5"
}
\end{displaymath} to represent
\begin{displaymath}
\xygraph{
!{<0cm,0cm>;<1cm,0cm>:<0cm,1cm>::}
!{(0,0)}*+{\bullet}="a1"
!{(0,1)}*+{\bullet}="a2"
!{(0,-1)}*+{\bullet}="a3"
!{(1,0.5)}*+{\bullet}="a4"
!{(1,-0.5)}*+{\bullet}="a5"
!{(2,0)}*+{\bullet}="a6"
!{(2,1)}*+{\bullet}="a7"
!{(2,-1)}*+{\bullet}="a8"
!{(-1,0)}*+{\bullet}="a9"
!{(-2,1)}*+{\bullet}="c1"
!{(-2,-1)}*+{\bullet}="c2"
!{(-3,0)}*+{\bullet}="c3"
!{(4,1)}*+{\bullet}="c4"
!{(3,0)}*+{\bullet}="c5"
!{(4,-1)}*+{\bullet}="c6"
!{(-4,0)}*+{\circ}="b1"
!{(5,0)}*+{\circ}="b2"
"c4"-"a7" "c5"-"a7" "c5"-"a6" "c5"-"a8" "c6"-"a8" "a7"-"a2"
"a7"-"a4" "a6"-"a4" "a6"-"a5" "a8"-"a5" "a8"-"a3" "a4"-"a2"
"a4"-"a1" "a4"-"a5" "a5"-"a1" "a5"-"a3" "a6"-"a7" "a6"-"a8"
"a1"-"a2" "a1"-"a3" "a9"-"a2" "a9"-"a1" "a9"-"a3" "c1"-"a9"
"c2"-"a9" "c1"-"c2" "c1"-"c3" "c2"-"c3" "c1"-"a2" "c2"-"a3"
"c4"-"c5" "c5"-"c6" "c4"-"c6" "b2"-"c4"
"b2"-"c6" "b1"-"c1" "b1"-"c2" "b1"-"c3"
}
\end{displaymath}\\

In the following picture, $\bullet$ represent the points mapped to $\infty$; $\circ$ represent the points mapped to $0.$

Now, for $n=(5/6)+2k,\ k\in \mathbb{N}$, the corresponding list is $[1^3 2^{11+30k},\\(4+6k)^1 3^{7+18k},5^{5+12k}]$, and the dessin is given by
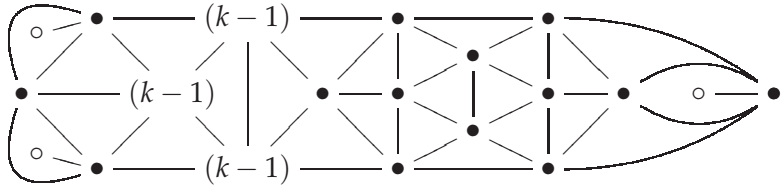
\begin{displaymath}
\xygraph{
!{<0cm,0cm>;<1cm,0cm>:<0cm,1cm>::}
!{(0,0)}*+{\bullet}="a1"
!{(0,1)}*+{\bullet}="a2"
!{(0,-1)}*+{\bullet}="a3"
!{(1,0.5)}*+{\bullet}="a4"
!{(1,-0.5)}*+{\bullet}="a5"
!{(2,0)}*+{\bullet}="a6"
!{(2,1)}*+{\bullet}="a7"
!{(2,-1)}*+{\bullet}="a8"
!{(-1,0)}*+{\bullet}="a9"
!{(-2,1)}*+{(k-1)}="c1"
!{(-2,-1)}*+{(k-1)}="c2"
!{(-3,0)}*+{(k-1)}="c3"
!{(5,0)}*+{\bullet}="c4"
!{(3,0)}*+{\bullet}="c5"
!{(5,0)}*+{\bullet}="c6"
!{(4,0)}*+{\circ}="b1"
!{(-4,1)}*+{\bullet}="a10"
!{(-4,-1)}*+{\bullet}="a11"
!{(-5,0)}*+{\bullet}="a12"
!{(-4.8,0.8)}*+{\circ}="b2"
!{(-4.8,-0.8)}*+{\circ}="b3"
"c4"-@/_0.3cm/"a7" "c5"-"a7" "c5"-"a6" "c5"-"a8" "c6"-@/^0.3cm/"a8" "a7"-"a2"
"a7"-"a4" "a6"-"a4" "a6"-"a5" "a8"-"a5" "a8"-"a3" "a4"-"a2"
"a4"-"a1" "a4"-"a5" "a5"-"a1" "a5"-"a3" "a6"-"a7" "a6"-"a8"
"a1"-"a2" "a1"-"a3" "a9"-"a2" "a9"-"a1" "a9"-"a3" "c1"-"a9"
"c2"-"a9" "c1"-"c2" "c1"-"c3" "c2"-"c3" "c1"-"a2" "c2"-"a3"
"c4"-@/_0.4cm/"c5" "c4"-@/^0.4cm/"c5" "c4"-"b1"
"a10"-"c1" "a10"-"c3" "a11"-"c2" "a11"-"c3" "a12"-"a10" "a12"-"a11"
"a12"-"c3" "a12"-@/^0.8cm/"a10" "a12"-@/_0.8cm/"a11" "a10"-"b2" "a11"-"b3"
}
\end{displaymath}\\

For $n=(7/6)+2k,\ k\in \mathbb{N}$, the corresponding list is 
$[1^3 2^{16+30k},\\(5+6k)3^{10+18k},5^{7+12k}]$, and the dessin is given by
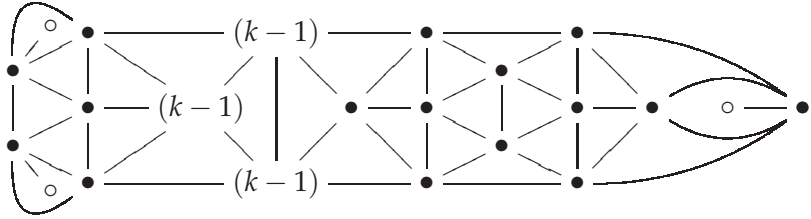
\begin{displaymath}
\xygraph{
!{<0cm,0cm>;<1cm,0cm>:<0cm,1cm>::}
!{(0,0)}*+{\bullet}="a1"
!{(0,1)}*+{\bullet}="a2"
!{(0,-1)}*+{\bullet}="a3"
!{(1,0.5)}*+{\bullet}="a4"
!{(1,-0.5)}*+{\bullet}="a5"
!{(2,0)}*+{\bullet}="a6"
!{(2,1)}*+{\bullet}="a7"
!{(2,-1)}*+{\bullet}="a8"
!{(-1,0)}*+{\bullet}="a9"
!{(-2,1)}*+{(k-1)}="c1"
!{(-2,-1)}*+{(k-1)}="c2"
!{(-3,0)}*+{(k-1)}="c3"
!{(5,0)}*+{\bullet}="c4"
!{(3,0)}*+{\bullet}="c5"
!{(5,0)}*+{\bullet}="c6"
!{(4,0)}*+{\circ}="b1"
!{(-4.5,0)}*+{\bullet}="a10"
!{(-4.5,1)}*+{\bullet}="a11"
!{(-4.5,-1)}*+{\bullet}="a12"
!{(-5.5,0.5)}*+{\bullet}="a13"
!{(-5.5,-0.5)}*+{\bullet}="a14"
!{(-5,1.1)}*+{\circ}="b2"
!{(-5,-1.1)}*+{\circ}="b3"
"c4"-@/_0.3cm/"a7" "c5"-"a7" "c5"-"a6" "c5"-"a8" "c6"-@/^0.3cm/"a8" "a7"-"a2"
"a7"-"a4" "a6"-"a4" "a6"-"a5" "a8"-"a5" "a8"-"a3" "a4"-"a2"
"a4"-"a1" "a4"-"a5" "a5"-"a1" "a5"-"a3" "a6"-"a7" "a6"-"a8"
"a1"-"a2" "a1"-"a3" "a9"-"a2" "a9"-"a1" "a9"-"a3" "c1"-"a9"
"c2"-"a9" "c1"-"c2" "c1"-"c3" "c2"-"c3" "c1"-"a2" "c2"-"a3"
"c4"-@/_0.4cm/"c5" "c4"-@/^0.4cm/"c5" "c4"-"b1"
"a11"-"c1" "a11"-"c3" "a12"-"c2" "a12"-"c3"
"a10"-"a11" "a10"-"a12" "a10"-"c3" "a13"-"a11"
"a13"-"a10" "a14"-"a10" "a14"-"a12" "a13"-"a14"
"a13"-"b2" "a14"-"b3" "a13"-@/^0.7cm/"a11" "a14"-@/_0.7cm/"a12"
}
\end{displaymath}\\

For $n=(11/6)+2k,\ k\in \mathbb{N}$, the corresponding list is 
$[1^3 2^{26+30k},\\(7+6k)3^{16+18k},5^{11+12k}]$, and the dessin is given by
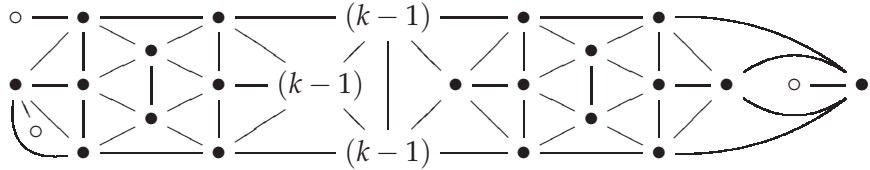
\begin{displaymath}
\xygraph{
!{<0cm,0cm>;<0.9cm,0cm>:<0cm,0.9cm>::}
!{(0,0)}*+{\bullet}="a1"
!{(0,1)}*+{\bullet}="a2"
!{(0,-1)}*+{\bullet}="a3"
!{(1,0.5)}*+{\bullet}="a4"
!{(1,-0.5)}*+{\bullet}="a5"
!{(2,0)}*+{\bullet}="a6"
!{(2,1)}*+{\bullet}="a7"
!{(2,-1)}*+{\bullet}="a8"
!{(-1,0)}*+{\bullet}="a9"
!{(-2,1)}*+{(k-1)}="c1"
!{(-2,-1)}*+{(k-1)}="c2"
!{(-3,0)}*+{(k-1)}="c3"
!{(5,0)}*+{\bullet}="c4"
!{(3,0)}*+{\bullet}="c5"
!{(5,0)}*+{\bullet}="c6"
!{(4,0)}*+{\circ}="b1"
!{(-4.5,0)}*+{\bullet}="a10"
!{(-4.5,1)}*+{\bullet}="a11"
!{(-4.5,-1)}*+{\bullet}="a12"
!{(-5.5,0.5)}*+{\bullet}="a13"
!{(-5.5,-0.5)}*+{\bullet}="a14"
!{(-6.5,1)}*+{\bullet}="a15"
!{(-6.5,0)}*+{\bullet}="a16"
!{(-6.5,-1)}*+{\bullet}="a17"
!{(-7.5,0)}*+{\bullet}="a18"
!{(-7.5,1)}*+{\circ}="b2"
!{(-7.2,-0.7)}*+{\circ}="b3"
"c4"-@/_0.3cm/"a7" "c5"-"a7" "c5"-"a6" "c5"-"a8" "c6"-@/^0.3cm/"a8" "a7"-"a2"
"a7"-"a4" "a6"-"a4" "a6"-"a5" "a8"-"a5" "a8"-"a3" "a4"-"a2"
"a4"-"a1" "a4"-"a5" "a5"-"a1" "a5"-"a3" "a6"-"a7" "a6"-"a8"
"a1"-"a2" "a1"-"a3" "a9"-"a2" "a9"-"a1" "a9"-"a3" "c1"-"a9"
"c2"-"a9" "c1"-"c2" "c1"-"c3" "c2"-"c3" "c1"-"a2" "c2"-"a3"
"c4"-@/_0.4cm/"c5" "c4"-@/^0.4cm/"c5" "c4"-"b1"
"a11"-"c1" "a11"-"c3" "a12"-"c2" "a12"-"c3"
"a10"-"a11" "a10"-"a12" "a10"-"c3" "a13"-"a11"
"a13"-"a10" "a14"-"a10" "a14"-"a12" "a13"-"a14"
"a15"-"a11" "a15"-"a13" "a16"-"a14" "a17"-"a14" "a17"-"a12"
"a15"-"a16" "a16"-"a17" "a18"-"a15" "a18"-"a16" "a18"-"a17"
"a15"-"b2" "a18"-@/_0.5cm/"a17" "a18"-"b3" "a16"-"a13" 
}
\end{displaymath}\\

For $n=(13/6)+2k,\ k\in \mathbb{N}$, the corresponding list is 
$[1^3 2^{31+30k},\\(8+6k)3^{19+18k},5^{13+12k}]$, and the dessin is given by
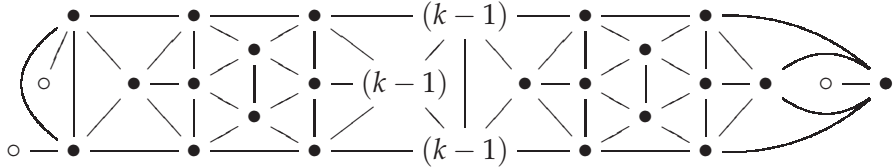
\begin{displaymath}
\xygraph{
!{<0cm,0cm>;<0.8cm,0cm>:<0cm,0.9cm>::}
!{(0,0)}*+{\bullet}="a1"
!{(0,1)}*+{\bullet}="a2"
!{(0,-1)}*+{\bullet}="a3"
!{(1,0.5)}*+{\bullet}="a4"
!{(1,-0.5)}*+{\bullet}="a5"
!{(2,0)}*+{\bullet}="a6"
!{(2,1)}*+{\bullet}="a7"
!{(2,-1)}*+{\bullet}="a8"
!{(-1,0)}*+{\bullet}="a9"
!{(-2,1)}*+{(k-1)}="c1"
!{(-2,-1)}*+{(k-1)}="c2"
!{(-3,0)}*+{(k-1)}="c3"
!{(5,0)}*+{\bullet}="c4"
!{(3,0)}*+{\bullet}="c5"
!{(5,0)}*+{\bullet}="c6"
!{(4,0)}*+{\circ}="b1"
!{(-4.5,0)}*+{\bullet}="a10"
!{(-4.5,1)}*+{\bullet}="a11"
!{(-4.5,-1)}*+{\bullet}="a12"
!{(-5.5,0.5)}*+{\bullet}="a13"
!{(-5.5,-0.5)}*+{\bullet}="a14"
!{(-6.5,1)}*+{\bullet}="a15"
!{(-6.5,0)}*+{\bullet}="a16"
!{(-6.5,-1)}*+{\bullet}="a17"
!{(-7.5,0)}*+{\bullet}="a18"
!{(-8.5,1)}*+{\bullet}="a19"
!{(-8.5,-1)}*+{\bullet}="a20"
!{(-9,0)}*+{\circ}="b2"
!{(-9.5,-1)}*+{\circ}="b3"
"c4"-@/_0.3cm/"a7" "c5"-"a7" "c5"-"a6" "c5"-"a8" "c6"-@/^0.3cm/"a8" "a7"-"a2"
"a7"-"a4" "a6"-"a4" "a6"-"a5" "a8"-"a5" "a8"-"a3" "a4"-"a2"
"a4"-"a1" "a4"-"a5" "a5"-"a1" "a5"-"a3" "a6"-"a7" "a6"-"a8"
"a1"-"a2" "a1"-"a3" "a9"-"a2" "a9"-"a1" "a9"-"a3" "c1"-"a9"
"c2"-"a9" "c1"-"c2" "c1"-"c3" "c2"-"c3" "c1"-"a2" "c2"-"a3"
"c4"-@/_0.4cm/"c5" "c4"-@/^0.4cm/"c5" "c4"-"b1"
"a11"-"c1" "a11"-"c3" "a12"-"c2" "a12"-"c3"
"a10"-"a11" "a10"-"a12" "a10"-"c3" "a13"-"a11"
"a13"-"a10" "a14"-"a10" "a14"-"a12" "a13"-"a14"
"a15"-"a11" "a15"-"a13" "a16"-"a14" "a17"-"a14" "a17"-"a12"
"a15"-"a16" "a16"-"a17" "a18"-"a15" "a18"-"a16" "a18"-"a17"
"a16"-"a13" 
"a19"-"a15" "a19"-"a18" "a20"-"a18" "a20"-"a17"
"a20"-"b3" "a19"-"b2" "a19"-"a20" "a19"-@/_0.7cm/"a20"
}
\end{displaymath}

The existence of $n=1/6,5/6,7/6,11/6,13/6$ is not hard.

\end{enumerate}

This complete the proof.
\end{proof}

\begin{corollary} \label{C1}
Let M be the monodromy group of $L_{n,B,g_2,g_3}y = 0$. Then, for each of the following pairs of M, n there exists a Lam\'e equation:
\begin{enumerate}
\item 
$M\cong G_{12}$ with $n\in \lbrace \pm1/4\rbrace+\mathbb{Z}$;
\item 
$M\cong G_{13}$ with $n\in \lbrace \pm1/6\rbrace+\mathbb{Z}$;
\item 
$M\cong G_{22}$ with $n\in \lbrace \pm1/10,\pm 3/10,\pm1/6\rbrace+\mathbb{Z}$.
\end{enumerate}
\end{corollary}

\begin{proof}
Let $y_1,y_2$ be the two linearly independent solutions of $L_{n,B,g_2,g_3}y = 0$. We have 
$ \displaystyle
\left(\frac{y_1}{y_2}\right)'=\frac{-w}{y_2^2}
$, where $w=\left| \begin{array}{cc}
y_1 & y_2 \\
y_1' &  y_2'
\end{array}\right| $ is the wronskian of $L_{n,B,g_2,g_3}y= 0$. Then, $y_2=\sqrt{\frac{-w}{(y_1/y_2)'}}$ is algebraic if $w$ and $\frac{y_1}{y_2}$ are both algebraic. In our case, $w$ is automatically algebraic, and $(y_1/y_2)'$ is algebraic provided by $L_{n,B,g_2,g_3} y = 0$ has finite projective monodromy. So, if $L_{n,B,g_2,g_3}y=0$ has finite projective monodromy then it has finite monodromy. By Theorem \ref{BD}, we complete the proof.
\end{proof}

\begin{remark}\label{R.R}
According to the case 2 of the proof of Theorem \ref{BM}, we can construct a dessin corresponding to the list $[1^2 2^{14},93^7,24^7]$ by taking $k=5$. However, it doesn't mean that  $L_{5/2,B,g_2,g_3}y=0$ has projective monodromy group isomorphic to $S_4$. Since, for $n\in \lbrace1/2\rbrace+\mathbb{Z}$, $L_{n,B,g_2,g_3}y = 0$ only have $V_4$, the subgroup of $S_4$, as its finite projective monodromy group. In this case, it means that the pullback function is a composition $f=h\circ g$, where $h$ realizes $L_{5/2,B,g_2,g_3}$ as the pullback of $H_{1/2,1/2,1/2}$ and $g$ realized $H_{1/2,1/2,1/2}$ as the pullback of $H_{1/2,1/3,1/4}$. Here we recall an important theorem due to Ritt (c.f.~\cite{R}). It says that a ramified covering with an imprimitive subgroup is a composition of two (or more) coverings. The following example shows how to check whether a Belyi function has primitive monodromy group with its dessin being given.
\end{remark}

\begin{example}
We use the same case as in Remark \ref{R.R}. We construct a dessin corresponding to the list $[1^2 2^{14},93^7,24^7]$. It is given by
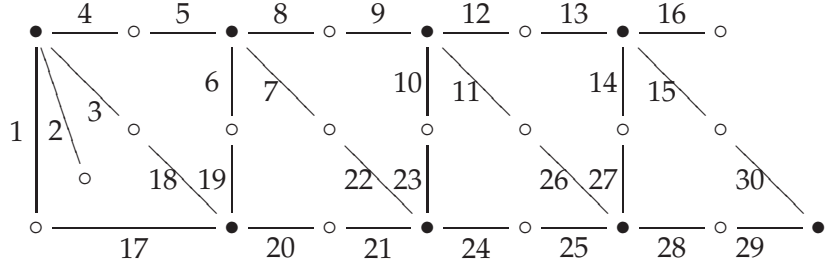
\begin{displaymath}
\xygraph{
!{<0cm,0cm>;<1.3cm,0cm>:<0cm,1.3cm>::}
!{(-4,0)}*+{\bullet}="a1"
!{(-2,0)}*+{\bullet}="a2"
!{(0,0)}*+{\bullet}="a3"
!{(2,0)}*+{\bullet}="a4"
!{(-2,-2)}*+{\bullet}="a5"
!{(0,-2)}*+{\bullet}="a6"
!{(2,-2)}*+{\bullet}="a7"
!{(4,-2)}*+{\bullet}="a8"
!{(-3,0)}*+{\circ}="b1"
!{(-1,0)}*+{\circ}="b2"
!{(1,0)}*+{\circ}="b3"
!{(3,0)}*+{\circ}="b4"
!{(-3.5,-1.5)}*+{\circ}="b5"
!{(-4,-2)}*+{\circ}="b6"
!{(-3,-1)}*+{\circ}="b7"
!{(-2,-1)}*+{\circ}="b8"
!{(-1,-1)}*+{\circ}="b9"
!{(0,-1)}*+{\circ}="b10"
!{(1,-1)}*+{\circ}="b11"
!{(2,-1)}*+{\circ}="b12"
!{(3,-1)}*+{\circ}="b13"
!{(-1,-2)}*+{\circ}="b14"
!{(1,-2)}*+{\circ}="b15"
!{(3,-2)}*+{\circ}="b16"
!{(-4.2,-1)}*+{1}
!{(-3.8,-1)}*+{2}
!{(-3.4,-0.8)}*+{3}
!{(-3.5,0.2)}*+{4}
!{(-2.5,0.2)}*+{5}
!{(-2.2,-0.5)}*+{6}
!{(-1.6,-0.6)}*+{7}
!{(-1.5,0.2)}*+{8}
!{(-0.5,0.2)}*+{9}
!{(-0.2,-0.5)}*+{10}
!{(0.4,-0.6)}*+{11}
!{(0.5,0.2)}*+{12}
!{(1.5,0.2)}*+{13}
!{(1.8,-0.5)}*+{14}
!{(2.4,-0.6)}*+{15}
!{(2.5,0.2)}*+{16}
!{(-3,-2.2)}*+{17}
!{(-2.7,-1.5)}*+{18}
!{(-2.2,-1.5)}*+{19}
!{(-1.5,-2.2)}*+{20}
!{(-0.5,-2.2)}*+{21}
!{(-0.7,-1.5)}*+{22}
!{(-0.2,-1.5)}*+{23}
!{(0.5,-2.2)}*+{24}
!{(1.5,-2.2)}*+{25}
!{(1.3,-1.5)}*+{26}
!{(1.8,-1.5)}*+{27}
!{(2.5,-2.2)}*+{28}
!{(3.3,-2.2)}*+{29}
!{(3.3,-1.5)}*+{30}
"a1"-"b1" "b1"-"a2" "a2"-"b2" "b2"-"a3" 
"a3"-"b3" "b3"-"a4" "a1"-"b5" 
"a1"-"b6" "a1"-"b7" "b7"-"a5" "a2"-"b8"
"b8"-"a5" "a2"-"b9" "b9"-"a6" "a3"-"b11"
"b11"-"a7" "a4"-"b12" "b12"-"a7" "a4"-"b13"
"b13"-"a8" "a2"-"b8" "b8"-"a5" "a4"-"b4" 
"b6"-"a5" "a5"-"b14" "b14"-"a6" "a6"-"b15"
"b15"-"a7" "a7"-"b16" "b16"-"a8"
"a3"-"b10" "b10"-"a6"
}
\end{displaymath}
This dessin corresponds to a Belyi function $f$ with $\sigma_0,\sigma_1,\sigma_{\infty}\in S_{30}$, where $\sigma_0,\sigma_1,\sigma_{\infty} $ be the monodromy of $f$ above $0,1,\infty$ respectively. We have

\begin{eqnarray*}
\begin{split}
\sigma_0=(2)(16)(1\ 17)(3\ 18)&(4\ 5)(6\ 19)(7\ 22)(8\ 9)(10\ 23)(11\ 26)\\
 &(12\ 13)(14\ 27)(15\ 30)(20\ 21)(24\ 25)(28\ 29),
 \end{split}
\end{eqnarray*}

\begin{eqnarray*}
\begin{split}
\sigma_1=(1\ 18\ 2)(3\ 19\ 5)(6\ 20\ 22)(7\ 23\ 9)&(10\ 24\ 26)(11\ 27\ 13)\\
&(14\ 28\ 30)(4\ 8\ 12\ 16\ 15\ 29\ 25\ 21\ 17),
 \end{split}
\end{eqnarray*}

\begin{eqnarray*}
\begin{split}
\sigma_{\infty}=(1\ 2\ 3\ 4)(5\ 6\ 7\ 8)(9\ 10\ 11\ 12)(13\ 14&\ 15\ 16)(17\ 18\ 19\ 20)\\&(21\ 22\ 23\ 24)(25\ 26\ 27\ 28)(29\ 30).
 \end{split}
\end{eqnarray*}
Now, let 

\begin{eqnarray*}
\begin{split}
P_1&=(2,4,5,7,14,16,22,24,25,27),\\
P_2&=(1,3,6,8,13,15,21,23,26,28),\\
P_3&=(9,10,11,12,17,18,19,20,29,30).
 \end{split}
\end{eqnarray*} 
be a partition of $\lbrace 1,2,\dots ,30\rbrace$. Also, let $\Sigma_0=(1)(2\ 3),\Sigma_1=(1\ 2\ 3)$ and $\Sigma_{\infty}=(1\ 2)(3)$ be three elements in $S_{3}$. We observe that $\sigma_{\alpha}(P_{i})=P_{\Sigma_{\alpha}(i)}$. It means that the monoddormy group of $f$ is imprimitive. Moreover, $\Sigma_i$ corresponds to a Belyi function which realizes $H_{1/2,1/2,1/4}$ as a pullback of $H_{1/2,1/3,1/4}$. Hence we have a composition $f=h\circ g$, where $h$ realizes $L_{5/2,B,g_2,g_3}$ as a pullback of $H_{1/2,1/2,1/4}$. In this case, we get a deeper information that $L_{5/2,B,g_2,g_3}y=0$ has finite projective monodromy group which is the subgroup of $D_4$.
\end{example}

\begin{remark}
For given $n,M$ in Corollary \ref{C1}, there is still a problem to compute the number of $L_{n,B,g_2,g_3}y = 0$ modulo equivalence have monodromy group isomorphic to $M$. The equivalence relation of the Lam\'e equation is given by $L_{n,B,g_2,g_3}y = 0\sim L_{n,B/a,g_2/a^2,g_3/a^3}y = 0$ with $a\in\mathbb{C}^{\times}$, which follows from the linear transformation $x\rightarrow ax$ that preserves the monodromy group. This question is equivalent to compute the number of dessin modulo equivalence corresponding to it. 

More precisely, according to the Grothendieck correspondence, for given $n,M$ in Corollary \ref{C1}, the number of dessin modulo equivalence is the same as the number of the pullback map modulo homography ($\displaystyle x\rightarrow\frac{ax+b}{cx+d},\\ ad-bc\neq 0$), the automorphism of $\mathbb{P}^1$. Furthermore, according to the Klein's theorem and the fact that $L_{n,B,g_2,g_3}$ has projective dihedral group if and only if $n\in \mathbb{Z}/2$, the pullback equation must have $S_4$, or $A_5$ as its projective monodromy group. Notice that we can further determine its monodromy group type by Corollary \ref{C1}. This lead to the conclusion that the number of Lam\'e equation with $n,M$ being given in corollary \ref{C1} modulo equivalence equals the number of dessins corresponding to it modulo equivalence. This counting problem is unsolved now. No explicit formula have been given except for the case $n\in\mathbb{Z}$. For the case $n\in\mathbb{Z}$, Dahmen completes the counting for the projective monodromy group level (see \cite{D} for the detail).
\end{remark}

\begin{example}
Here we give an explicit computation for $n=5/6$. 

In this case, if $L_{5/6,B,g_2,g_3}y=0$ has finite monodromy then it is a pullback of $H_{1/2,1/3,1/4}y=0$. The corresponding list is $[1^2 2^4,3^2 4^1,2^14^2]$, and all the possible dessins are given by 
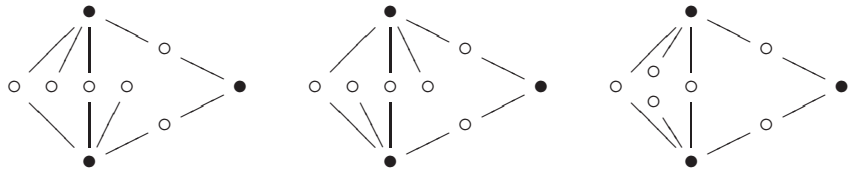
\begin{displaymath}
\xygraph{
!{<0cm,0cm>;<1cm,0cm>:<0cm,1cm>::}
!{(0,1)}*+{\bullet}="a1"
!{(0,-1)}*+{\bullet}="a2"
!{(2,0)}*+{\bullet}="a3"
!{(-1,0)}*+{\circ}="b1"
!{(-0.5,0)}*+{\circ}="b2"
!{(0,0)}*+{\circ}="b3"
!{(0.5,0)}*+{\circ}="b4"
!{(1,0.5)}*+{\circ}="b5"
!{(1,-0.5)}*+{\circ}="b6"
!{(4,1)}*+{\bullet}="aa1"
!{(4,-1)}*+{\bullet}="aa2"
!{(6,0)}*+{\bullet}="aa3"
!{(3,0)}*+{\circ}="bb1"
!{(3.5,0)}*+{\circ}="bb2"
!{(4,0)}*+{\circ}="bb3"
!{(4.5,0)}*+{\circ}="bb4"
!{(5,0.5)}*+{\circ}="bb5"
!{(5,-0.5)}*+{\circ}="bb6"
!{(8,1)}*+{\bullet}="aaa1"
!{(8,-1)}*+{\bullet}="aaa2"
!{(10,0)}*+{\bullet}="aaa3"
!{(7,0)}*+{\circ}="bbb1"
!{(7.5,0.2)}*+{\circ}="bbb4"
!{(8,0)}*+{\circ}="bbb3"
!{(7.5,-0.2)}*+{\circ}="bbb2"
!{(9,0.5)}*+{\circ}="bbb5"
!{(9,-0.5)}*+{\circ}="bbb6"
"a1"-"b1" "a2"-"b1" "a1"-"b2" "a1"-"b3"
"a2"-"b3" "a2"-"b4" "a1"-"b5" "a3"-"b5"
"a2"-"b6" "a3"-"b6"
"aa1"-"bb1" "aa2"-"bb1" "aa1"-"bb4" "aa1"-"bb3"
"aa2"-"bb3" "aa2"-"bb2" "aa1"-"bb5" "aa3"-"bb5"
"aa2"-"bb6" "aa3"-"bb6"
"aaa1"-"bbb1" "aaa2"-"bbb1" "aaa1"-"bbb4" "aaa1"-"bbb3"
"aaa2"-"bbb3" "aaa2"-"bbb2" "aaa1"-"bbb5" "aaa3"-"bbb5"
"aaa2"-"bbb6" "aaa3"-"bbb6"
}
\end{displaymath}
where $\bullet$ represent the points mapped to $\infty$, $\circ$ represent the points mapped to 0. It means that if we fix $B=0$ or $1$, there are only three tori such that $L_{5/6,B,g_2,g_3}$ has finite monodromy group isomorphic to $G_{13}$. This is consistent with the result obtained in \cite[p.134]{Waall} based on computer algorithms.

To construct the Belyi functions corresponding to the dessins given above, we can use the method of undetermined coefficients. Let $f$ be the corresponding Belyi function. Then we have 
\[
\begin{split}
f(x) &= K\frac{(x^4+px^3+qx^2+rx+s)^2(x-m)(x-n)}{(x^2+Ax+B)^4(x-C)^3},\\
f(x)-1 &= K\frac{(x^2+ax+b)^3(x-c)^4}{(x^2+Ax+B)^4(x-C)^2}.
\end{split}
\]
By eliminating $f$, we obtain a system of $10$ algebraic relations with $13$ unknowns $a,b,c,m,n,p,q,r,s,K,A,B,C.$
The 3 remaining degrees of freedom correspond to the probability  of making a m\"obius transformation of $x$. There are many ways to fix the freedom. For example, we can assume that $A=0,B=1$ and $C=0$. Hence we can solve the system of $10$ algebraic relations with $10$ unknowns $a,b,c,m,n,p,q,r,s,K$ theoretically. However, it needs a large amount of calculation and takes hours even using a computer software like \emph{Mathematica}. Instead, we give another method that enable us to construct one of them in a more efficient manner. 

We first recall the fact that $L_{5/6,0,1,0}$ is a pullback of $H_{1/2,1/3,1/4}$ by $x^2+1$ (c.f.~\cite{M}, Lemma 3.3). Now, let $F$ be the rational function which realizes $H_{1/2,2/3,1/4}$ as a pullback of $H_{1/2,1/3,1/4}$. Then the corresponding list is $[1^1 2^2, 2^1 3^1, 1^1 4^1]$, and the only possible dessin is given by 
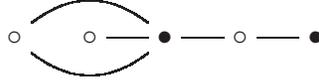
\begin{displaymath}
\xygraph{
!{<0cm,0cm>;<1cm,0cm>:<0cm,1cm>::}
!{(3,1)}*+{\bullet}="a1"
!{(5,1)}*+{\bullet}="a2"
!{(1,1)}*+{\circ}="b1"
!{(2,1)}*+{\circ}="b2"
!{(4,1)}*+{\circ}="b3"
"b2"-"a1" "a1"-"b3" "b3"-"a2" "a1"-@/_5mm/"b1" "a1"-@/^5mm/"b1"
}
\end{displaymath}

Also, we can use the method of undetermined coefficients to construct $F$. Then we have 
\[
\begin{split}
F(x) &= k \frac{(x-p)(x^2+qx+r)^2}{(x-A)(x-B)^4};\\
F(x)-1 &= k\frac{(x-s)^2(x-t)^3}{(x-A)(x-B)^4}.
\end{split}
\]
For convenience, we take $p=0$, $s=1$ and $B=\infty$. By direct computations, we find
\[
F(x)=\frac{-2x(2x^2-20x+45)^2}{27(5x-32)}.
\]
Hence we get the conclusion that $G(x)=F(x^2+1)$ is a rational function which realizes $L_{5/6,0,1,0}$ as a pullback of $H_{1/2,1/3,1/4}$. This completes the construction.
\end{example}

\section{Finite monodromy problem for the generalized Lam\'e equations with 2 singular points}

The main goal of this section is to study the finite monodromy problem of equation (\ref{Alg-2sing-Lame}). We start with a general result due to Chai, Lin, and Wang:

\begin{theorem}[\cite{Boss3}]
Let $PM$ be the projective monodromy group of the generalized Lam\'e equation (\ref{Wei-gene-Lame}). Suppose $n_i\in\mathbb{Z}/2$ and the solutions are all free from logarithm. Then we have 
\begin{enumerate}
\item 
$PM\cong K_4$, if $\sum n_i\in \lbrace 1/2\rbrace+\mathbb{Z}$,
\item 
$PM\cong K_4 $ or $C_k$ , if $\sum n_i\in \mathbb{Z}$ and $PM$ is finite.
\end{enumerate}
\end{theorem}

\begin{proof}
We sketch the proof here and refer to \cite{Boss3} for the details. Let $\sigma_{p_i}$ be the local monodromy matrix at $p_i$. We have $\sigma_{p_i}=(-1)^{2n_i}I_2$, provided by the log-free assumption and $n_i\in\mathbb{Z}/2$. Notice that this matrix is identity in $PGL(2,\mathbb{C})$ and is independent of the choice of the basis. Let $S_1$, $S_2$ be the monodromy representation for the standard homology basis $w_1$ and $w_2$ respectively. Then the homotypo relation gives that 
$$ S_2^{-1}S_1^{-1}S_2S_1=\prod_{i=1}^{s}\sigma_{p_i}=(-1)^{2\sum n_i}I_2 $$.

If $\sum n_i\in\mathbb{Z}$, we have $S_2S_1=S_1S_2$. In this case, the projective monodromy group of the equation (\ref{Wei-gene-Lame}) is abelian. The conclusion follows from the fact that $K_4$ and $C_k$ are the only possible abelian groups in basic Schwarz list.

If $\sum n_i\in \lbrace 1/2\rbrace+\mathbb{Z}$, we have $S_2S_1=-S_1S_2$. Now, let $v$ be an eigenvector of $S_1$ with $S_1v=\lambda_1 v$ and $\lambda\neq 0$. Then 
$$S_1(S_2 v)=-S_2S_1v=-\lambda_1(S_2v).$$
It means that $S_2v$ is also an eigenvector of $S_1$ with eigenvalue $-\lambda_1\neq \lambda_1$. We have $S_1^2=\lambda_1^2 I_2$. With the same computation, we have $S_2^2=\lambda_2^2 I_2$ for some $\lambda_2\neq 0$. It is clear now that $PM\cong K_4$ generated by $S_1$ and $S_2$.
\end{proof}

\begin{corollary}\label{C-PM-Coro}
Let $PM'$ be the projective monodromy group of equation (\ref{C-Alg-Lame}). Suppose $n_i\in\mathbb{Z}/2$ and $PM'$ is finite. Then 
\begin{enumerate}
\item $PM'\cong K_4$, if $n_0\in \lbrace 1/2\rbrace+ \mathbb{Z}$;
\item $PM'\cong D_N$ for some $N\in\mathbb{N}_{\geq 2}$, if $n_0 \in \mathbb{Z}$.
\end{enumerate}
\end{corollary}

\begin{proof}
Let $\gamma_1,\gamma_2,\gamma_3$ and $\gamma_{\infty} $ be the local monodromy matrices at $z=e_1,e_2,e_3$ and $\infty $ respectively. The double cover $z=\wp(a)$ corresponds to the relation $S_1=\gamma_1\gamma_{\infty}^{-1}$ and $S_2=\gamma_2\gamma_{\infty}^{-1}$. Let $\sigma_{\wp(p_i)}$ be the local monodromy at $\wp(a)$. We have $\sigma_{\wp(p_i)}=(-1)^{4n_i}I_2$, by checking its degree. Hence, $\sigma_{\wp(p_i)}=I_2$, provided by the assumption $n_i\in\mathbb{Z}/2$. Now, we have $PM'$ is generated by $PM$ and $\gamma_{\infty}$.

If $n_{0}\in \lbrace 1/2\rbrace+\mathbb{Z} $, $\gamma_{\infty}=I_2$ in $PGL(2,\mathbb{C})$. In this case, $PM'\cong PM\cong K_4$.

If $n_{0}\in\mathbb{Z}$, then $\gamma_{\infty}$ is an element with order $2$ in $PGL(2,\mathbb{C})$. By using the fact that $PM'$ is never finite cyclic, we have $PM\cong D_N$ for some $N\geq 2$
\end{proof}

\begin{theorem} \label{A-Lame-Cyc}
The projective monodromy group of equation (\ref{C-Alg-Lame}) is never finite cyclic. 
\end{theorem}

\begin{proof}
The method of the proof is similar with the method in original Lam\'e equation (c.f. \cite{BD}). We just sketch the proof. 

Suppose the monodromy group $M$ of the equation (\ref{C-Alg-Lame}) acts completely reducibly. We can assume the two solutions is of the form
$$ 
y_1=\prod_{i=1}^3 (x-e_i)^{\epsilon_i}\prod_{j=1}^s (x-\wp(p_j))^{\mu_j}P_1(x),
$$
$$ 
y_2=\prod_{i=1}^3 (x-e_i)^{1/2-\epsilon_i}\prod_{j=1}^s (x-\wp(p_j))^{1-\mu_j}P_2(x)
$$ 
with $\epsilon_i\in\lbrace 0,1/2\rbrace $, $\mu_j\in\lbrace -n_j,n_j+1\rbrace$. $P_1(x),P_2(x)$ are fixed under the action of $M$ and hence they are contained in $\mathbb{C}(x)$. Moreover, $P_1(x),P_2(x)\in\mathbb{C}[x]$, since they have no finite poles. Now, we have
$$-\sum_{i=1}^3 e_i-\sum_{j=1}^s \mu_j-m_1=-n_0/2,$$
$$\sum_{i=1}^3 e_i+\sum_{j=1}^s \mu_j-3/2-s-m_2=(n_0+1)/2.$$ 
where $m_1$ and $m_2$ be the degree of $P_1$ and $P_2$ respectively. Addition of the two equation above leads to a contradiction that $-3/2-s-m_1-m_2=1/2$. Hence $M$ is not completely reducibly.

Now, if $PM\cong C_k$ with $k<\infty$, the group $M$ would be reducible. $M$ is also finite since equation (\ref{C-Alg-Lame}) has algebraic Wronskian. According to Maschke's theorem, $M$ acts completely reducibly on the solution space. It contradicts with the above conclusion.
\end{proof}

We come back to the monodromy problem of equation (\ref{Alg-2sing-Lame}), we divide it into the following two cases:

\begin{enumerate}
\item $n\notin \mathbb{Z}/2$. This is the case that the pullback function in Theorem \ref{Klein} is a Belyi function. That is, the pullback function is unramified on $\mathbb{P}^1\setminus \lbrace 0,1,\infty\rbrace$.
\item $n\in\mathbb{Z}/2$. The exponent difference of the equation (\ref{Alg-2sing-Lame}) at $\wp(a)$ is an integer in this case. So, the pullback function in Theorem \ref{Klein} may not be a Belyi function.
\end{enumerate} 

According to the Corollary \ref{C-PM-Coro}, if $n\in\mathbb{Z}/2$, then the finite projective monodromy group of equation (\ref{Alg-2sing-Lame}) is $D_N$. The remaining problem is to determine all the possible $N$ with given $n$ and further compute the number of the equation modulo linear equivalence. This problem will relate to the modulo form problem. We will study it in the following paper.

The remaining content of this section will only focus on the case 1. From now on, we assume $n\notin \mathbb{Z}/2$. The following theorem restricts the possible cases of $n$.

\begin{theorem} \label{BM1}
Let M be the monodromy group of equation (\ref{Alg-2sing-Lame}) and $n \notin\mathbb{N}/2$. Then 
\begin{enumerate}
\item 
$M \cong G(N,N/2,2)\Rightarrow n\in\lbrace\pm 1/4\rbrace+\mathbb{Z}$, $N$ is even and $N\geq 4$;
\item
 $M \cong G_{12} \Rightarrow n\in\lbrace \pm 1/6,\pm 1/4,\pm1/3\rbrace+\mathbb{Z}$;
\item 
$M\cong G_{13}\Rightarrow n\in\lbrace\pm 1/8,\pm3/8,\pm 1/6,\pm 1/4,\pm 1/3\rbrace +\mathbb{Z}$ ;
\item 
$M\cong G_{22}\Rightarrow n\in(\mathbb{Z}/6)\cup(\mathbb{Z}/10) \setminus(\mathbb{Z}/2) $.
\end{enumerate}
Moreover, if $n\geq -1/2$, then $n\geq 1/6$.
\end{theorem}

\begin{proof}
The proof is similar to the proof of Theorem \ref{BD} (\cite{BA}, Theorem 4.4). From now on, we assume $n\geq -1/2$, since substituting $n$ to $-n-1$ does not change the equation.
\begin{enumerate}

\item 
This group has invariant polynomial with degree 4. Hence, the fourth symmetric product of the equation (\ref{Alg-2sing-Lame}) has solution invariant under the action of monodromy i.e.~the rational function. This leads to one of the local exponents $-4n$, $-2n+1$, $2$, $2n+3$, $4n+4$ at $z=\wp(a)$ is an integer. Furthermore, the solution is not a polynomial, since the local exponents are $0$, $1/2$, $1$, $3/2$, $2$ at $z=\infty$. We only need one of $-4n$, $-2n+1$ is an integer. However, $n\notin\mathbb{Z}/2$, or the group would has invariant polynomial with degree 2. This leads to a contradiction. 

\item 
Let $\sigma$ be the sign of the permutation of $M$ modulo scalars. $\delta$ be the determinant of the element of $M$. We will prove that $\sigma\delta (g)=1$, $\forall\ g\in G_{12}$. First, if $g$ has 3-cycle as projection, then $g^3$ is a scalar multiple of Id. Hence, $\sigma(g)=1$. Second, at least one element $h\in G_{12}$ has 2-cycle as projection has determinant $-1$, since $G_{12}$ is generated by $\pm$Id and the 2-cycle in $S_{4}$. At last, $S_{4}$ is generated by $A_{4}$ and $h$. We get the conclusion that $\sigma\delta(g)=1$, $\forall g\in G_{12}$. Note that the local monodromy  $\gamma_1,\gamma_2,\gamma_3,\gamma_{\infty}$ have determinant $-1$ and hence have order $2$. So, the local monodromy at $\wp(a)$ have determinant $1$ and hence has even permutation. It has order 2 or 3. This implies that $2n+1$, the exponent difference at $\wp(a)$, is contained in either $\lbrace 1/2\rbrace +\mathbb{Z}$ or $\lbrace\pm 1/3\rbrace+\mathbb{Z}$. Thus we have $n\in\lbrace \pm1/6,\pm1/4,\pm1/3\rbrace +\mathbb{Z}$.

\item 
This group is generated by the element of $G_{12}$ and $i\cdot$Id. Let $\sigma,\delta$ be the same as (2). Notice that $\sigma\delta(i\cdot{\rm Id})=-1$. Since $\gamma_1,\gamma_2,\gamma_3,\gamma_{\infty}$ generate $G_{13}$, at least one of them must have value $-1$ when taking $\sigma\delta$. Also, if $\gamma_i$ has even permutation, its projection must be the product of 2-cycles. It is because $\gamma_i$ has eigenvalues $-1$ and $1$. Note that $S_4$ cannot be generated by its subgroup $V_4$ and a 2-cycle. This arises two cases:
\begin{enumerate}

\item 
Two of $\gamma_i$ has the product of 2-cycles as  its projection. In this case, the local monodromy at $\wp(a)$ has even permutation. This leads to the condition $2n+1\in\lbrace 1/2,\pm1/3\rbrace+\mathbb{Z}$ and hence $n\in\lbrace \pm1/6,\pm1/4,\pm1/3\rbrace+\mathbb{Z}$.

\item 
One of $\gamma_i$ has the product of 2-cycles as  its projection. In this case, the local monodromy at $\wp(a)$ has odd permutation. This leads to the conditon  $2n+1\in\lbrace 1/2,\pm1/4\rbrace$ and hence $n\in\lbrace \pm1/4,\pm1/8,\pm3/8\rbrace$.
\end{enumerate}
We complete all the cases.

\item 
In this case, we check the exponent difference. We need $2n+1$, the exponent difference at $\wp(a)$, contained in $\lbrace 1/2,\pm1/3,\pm1/5,\pm2/5\rbrace+\mathbb{Z}$, since it is the pullback of the equation $H_{1/2,1/3,1/5}y=0$ which has exponent difference $1/2,1/3,1/5$ at $0,1,\infty$ respectively. This leads to $n\in\lbrace \pm1/10,\pm3/10,\pm1/5,\pm2/5,\pm1/6,\pm1/4,\pm1/3\rbrace+\mathbb{Z}$.
\end{enumerate}
To finish the last statement of the proof, we need to exclude the possible that $-1/2\leq n<1/6$. For the case $-1/2\leq n<0$, we observe that the degree of the pullback function in Klein's theorem is less than zero. It is impossible. Now we only need to exclude the cases $n=1/8,1/10$. In these two cases, we can easily observe that there exists a point whose ramification index will greater than the degree of the pullback function in Klein's theorem. We get the contradiction. 
\end{proof}

\begin{theorem} \label{BM2}
Let M be the projective monodromy group of the generalized Lam\'e equation. Then, for each of the following pairs of M, n there exists a generalized Lam\'e equation (\ref{Alg-2sing-Lame}):
\begin{enumerate}
\item 
$M\cong G_{12}$ with $n\in\lbrace \pm1/6,\pm1/4,\pm1/3\rbrace+\mathbb{Z}$;
\item 
$M\cong G_{13}$ with $n\in\lbrace\pm1/8,\pm3/8,\pm1/6,\pm1/4,\pm1/3\rbrace +\mathbb{Z}$;
\item 
$M\cong G_{22}$ with $n\in(\mathbb{Z}/6)\cup(\mathbb{Z}/10) \setminus(\mathbb{Z}/2) $.
\end{enumerate}
\end{theorem}

\begin{proof}
According to Theorem \ref{BM1}, we only need to consider the case $n\geq1/6$.
The proof of this theorem has two differences from the proof of Theorem \ref{BM}. 

For $n\in \lbrace \pm1/6,\pm1/4,\pm1/3\rbrace$, it is not enough to show the existence of the equation has $S_4$ as its projective monodromy group and determine its monodromy group as in Corollary \ref{C1}, since it has $G_{12},G_{13}$ two possibilities. To solve this problem, we check whether $e_1,e_2,e_3,\infty $ all map to $0$. If it is, then the local monodromy $\gamma_1,\gamma_2,\gamma_3,\gamma_{\infty}$ are all two cycles. It is the case 2 in the proof of Theorem \ref{BM1}. Otherwise, at least one of $\gamma_1,\gamma_2,\gamma_3,\gamma_{\infty} $ will be the product of 2 cycles. It is because at least one of $e_1,e_2,e_3,\infty $ will map to $\infty$ with ramification index equals $2$. It relates to the case 3 in theorem \ref{BM1}.

For $n\in \lbrace\pm1/4\rbrace+\mathbb{Z}$, we need to construct the dessin with primitive monodromy group, since in this case, it may have dihedral as its projective monodromy group.

\begin{enumerate}[label=\textbf{case \arabic* :}]

\item 
($M=G_{12},\ n\in\lbrace \pm 1/4\rbrace+\mathbb{Z},\ n>0$)

For $n=(1+2k)/4,k\in\mathbb{N}$, the corresponding list is $[1^{4}2^{6k+1},3^{4k+2},\\(4k+6)^{1}4^{2k}]$. If $k$ is even, the dessin is given by 
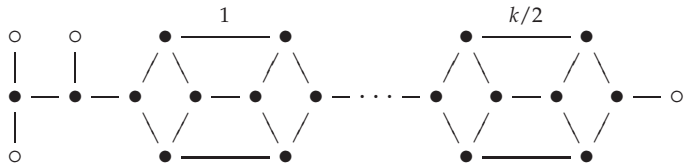
\begin{displaymath}
\xygraph{
!{<0cm,0cm>;<0.8cm,0cm>:<0cm,0.8cm>::}
!{(0,0)}*+{\cdots}="c1"
!{(1,0)}*+{\bullet}="a1"
!{(2,0)}*+{\bullet}="a2"
!{(3,0)}*+{\bullet}="a3"
!{(4,0)}*+{\bullet}="a4"
!{(1.5,1)}*+{\bullet}="a5"
!{(3.5,1)}*+{\bullet}="a6"
!{(1.5,-1)}*+{\bullet}="a7"
!{(3.5,-1)}*+{\bullet}="a8"
!{(-1,0)}*+{\bullet}="a9"
!{(-2,0)}*+{\bullet}="a10"
!{(-3,0)}*+{\bullet}="a11"
!{(-4,0)}*+{\bullet}="a12"
!{(-1.5,1)}*+{\bullet}="a13"
!{(-3.5,1)}*+{\bullet}="a14"
!{(-1.5,-1)}*+{\bullet}="a15"
!{(-3.5,-1)}*+{\bullet}="a16"
!{(-5,0)}*+{\bullet}="a17"
!{(-6,0)}*+{\bullet}="a18"
!{(5,0)}*+{\circ}="b1"
!{(-5,1)}*+{\circ}="b2"
!{(-6,1)}*+{\circ}="b3"
!{(-6,-1)}*+{\circ}="b4"
!{(2.5,1.3)}*+{^{k/2}}
!{(-2.5,1.3)}*+{^1}
"c1"-"a1" "a1"-"a5" "a1"-"a7" "a5"-"a2"
"a7"-"a2" "a5"-"a6" "a2"-"a3" "a7"-"a8"
"a3"-"a6" "a3"-"a8" "a6"-"a4" "a8"-"a4"
"c1"-"a9" "a9"-"a13" "a9"-"a15" "a13"-"a10"
"a10"-"a15" "a13"-"a14" "a10"-"a11" "a15"-"a16"
"a11"-"a14" "a11"-"a16" "a14"-"a12" "a16"-"a12"
"a12"-"a17" "a17"-"a18" "a4"-"b1" "a17"-"b2"
"a18"-"b3" "a18"-"b4" 
}
\end{displaymath}\\
If $k$ is odd, the dessin is given by
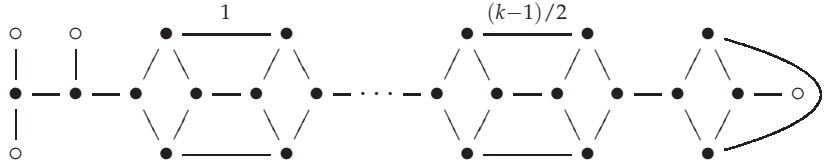
\begin{displaymath}
\xygraph{
!{<0cm,0cm>;<0.8cm,0cm>:<0cm,0.8cm>::}
!{(0,0)}*+{\cdots}="c1"
!{(1,0)}*+{\bullet}="a1"
!{(2,0)}*+{\bullet}="a2"
!{(3,0)}*+{\bullet}="a3"
!{(4,0)}*+{\bullet}="a4"
!{(1.5,1)}*+{\bullet}="a5"
!{(3.5,1)}*+{\bullet}="a6"
!{(1.5,-1)}*+{\bullet}="a7"
!{(3.5,-1)}*+{\bullet}="a8"
!{(5,0)}*+{\bullet}="a19"
!{(6,0)}*+{\bullet}="a20"
!{(5.5,1)}*+{\bullet}="a21"
!{(5.5,-1)}*+{\bullet}="a22"
!{(-1,0)}*+{\bullet}="a9"
!{(-2,0)}*+{\bullet}="a10"
!{(-3,0)}*+{\bullet}="a11"
!{(-4,0)}*+{\bullet}="a12"
!{(-1.5,1)}*+{\bullet}="a13"
!{(-3.5,1)}*+{\bullet}="a14"
!{(-1.5,-1)}*+{\bullet}="a15"
!{(-3.5,-1)}*+{\bullet}="a16"
!{(-5,0)}*+{\bullet}="a17"
!{(-6,0)}*+{\bullet}="a18"
!{(7,0)}*+{\circ}="b1"
!{(-5,1)}*+{\circ}="b2"
!{(-6,1)}*+{\circ}="b3"
!{(-6,-1)}*+{\circ}="b4"
!{(2.5,1.3)}*+{^{(k-1)/2}}
!{(-2.5,1.3)}*+{^1}
"c1"-"a1" "a1"-"a5" "a1"-"a7" "a5"-"a2"
"a7"-"a2" "a5"-"a6" "a2"-"a3" "a7"-"a8"
"a3"-"a6" "a3"-"a8" "a6"-"a4" "a8"-"a4"
"c1"-"a9" "a9"-"a13" "a9"-"a15" "a13"-"a10"
"a10"-"a15" "a13"-"a14" "a10"-"a11" "a15"-"a16"
"a11"-"a14" "a11"-"a16" "a14"-"a12" "a16"-"a12"
"a12"-"a17" "a17"-"a18" "a17"-"b2"
"a18"-"b3" "a18"-"b4" 
"a4"-"a19" "a19"-"a21" "a19"-"a22"
"a21"-"a20" "a22"-"a20" "a20"-"b1"
"a21"-@/^1.5cm/"a22"
}
\end{displaymath}\\
where $\bullet$ represent the points mapped to 1; $\circ$ represent the points mapped to 0.

Note that the dessin given above has primitive monodromy group. See Remark \ref{R.k/4} for the detail.

\item 
($M=G_{12},n\in\lbrace\pm1/6,\pm1/3\rbrace+\mathbb{Z}$)

For $n=k/6,k\in\mathbb{N}$, the corresponding list is $[1^{4}2^{2k-2},(k+3)^1 3^{k-1},\\4^{k}]$. If $k$ is odd, the dessin is given by
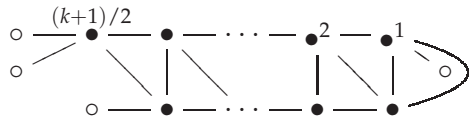
\begin{displaymath}
\xygraph{
!{<0cm,0cm>;<1cm,0cm>:<0cm,1cm>::}
!{(1,0)}*+{\bullet^2}="a1"
!{(2,0)}*+{\bullet^1}="a2"
!{(2,-1)}*+{\bullet}="a4"
!{(1,-1)}*+{\bullet}="a5"
!{(-1,0)}*+{\bullet}="a6"
!{(-2,0)}*+{\bullet}="a7"
!{(-2,0.2)}*+{^{(k+1)/2}}
!{(-1,-1)}*+{\bullet}="a9"
!{(0,0)}*+{\cdots}="c1"
!{(-2,-1)}*+{\circ}="c2"
!{(-3,-0.5)}*+{\circ}="c3"
!{(0,-1)}*+{\cdots}="c4"
!{(-3,0)}*+{\circ}="b2"
!{(2.7,-0.5)}*+{\circ}="b7"
"c1"-"a1" "c1"-"a6" 
"a1"-"a5" "a1"-"a2" "a2"-"a4" "a4"-"a5"
"a6"-"a7" "a6"-"a7" "b2"-"a7" 
"a6"-"a9" "a7"-"a9" "a1"-"a4" "a2"-@/^1cm/"a4"
"a2"-"b7" "a9"-"c2" "a7"-"c3"
"c4"-"a9" "c4"-"a6" "c4"-"a5"
}
\end{displaymath}\\
If $k$ is even, the dessin is given by
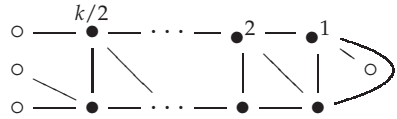
\begin{displaymath}
\xygraph{
!{<0cm,0cm>;<1cm,0cm>:<0cm,1cm>::}
!{(1,0)}*+{\bullet^2}="a1"
!{(2,0)}*+{\bullet^1}="a2"
!{(2,-1)}*+{\bullet}="a4"
!{(1,-1)}*+{\bullet}="a5"
!{(-1,0)}*+{\bullet}="a6"
!{(-1,0.2)}*+{^{k/2}}
!{(-1,-1)}*+{\bullet}="a9"
!{(0,0)}*+{\cdots}="c1"
!{(-2,-1)}*+{\circ}="c2"
!{(-2,-0.5)}*+{\circ}="c3"
!{(0,-1)}*+{\cdots}="c4"
!{(-2,0)}*+{\circ}="b2"
!{(2.7,-0.5)}*+{\circ}="b7"
"c1"-"a1" "c1"-"a6" 
"a1"-"a5" "a1"-"a2" "a2"-"a4" "a4"-"a5"
"a6"-"a9"  "a1"-"a4" "a2"-@/^1cm/"a4"
"a2"-"b7" "a9"-"c2" "a6"-"b2" "a9"-"c3"
"a9"-"c4" "c4"-"a6" "c4"-"a5"
}
\end{displaymath}\\
where $\bullet$ represent the points mapped to $\infty$; $\circ$ represent the points mapped to 0.

\item 
($M=G_{13},\ n\in\lbrace \pm1/4\rbrace+\mathbb{Z}$)

For $n=(1+2k)/4,k\in\mathbb{N}$, the corresponding list is $[1^3 (2k+3^1)\\2^{5k},3^{4k+2},2^14^{2+3k}]$. To construct this dessin, we recall the dessin given in case 1. If $k$ is even, we construct the dessin corresponding to the list $[1^4 2^{6k+1},3^{4k+2},(4k+6)^1 4^{2k}]$. It is given by
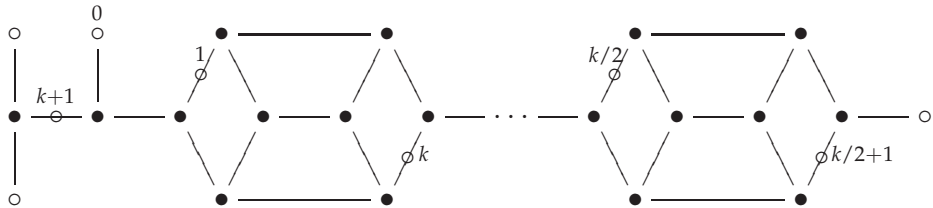
\begin{displaymath}
\xygraph{
!{<0cm,0cm>;<1.1cm,0cm>:<0cm,1.1cm>::}
!{(0,0)}*+{\cdots}="c1"
!{(1,0)}*+{\bullet}="a1"
!{(2,0)}*+{\bullet}="a2"
!{(1.25,0.5)}*+{\circ}="p1"
!{(1.15,0.7)}*+{^{k/2}}
!{(3.75,-0.5)}*+{\circ}="p2"
!{(4.25,-0.5)}*+{^{k/2+1}}
!{(-1.25,-0.5)}*+{\circ}="p3"
!{(-1.05,-0.5)}*+{^{k}}
!{(-3.75,0.5)}*+{\circ}="p4"
!{(-3.75,0.7)}*+{^{1}}
!{(-5.5,0)}*+{\circ}="p5"
!{(-5.5,0.2)}*+{^{k+1}}
!{(3,0)}*+{\bullet}="a3"
!{(4,0)}*+{\bullet}="a4"
!{(1.5,1)}*+{\bullet}="a5"
!{(3.5,1)}*+{\bullet}="a6"
!{(1.5,-1)}*+{\bullet}="a7"
!{(3.5,-1)}*+{\bullet}="a8"
!{(-1,0)}*+{\bullet}="a9"
!{(-2,0)}*+{\bullet}="a10"
!{(-3,0)}*+{\bullet}="a11"
!{(-4,0)}*+{\bullet}="a12"
!{(-1.5,1)}*+{\bullet}="a13"
!{(-3.5,1)}*+{\bullet}="a14"
!{(-1.5,-1)}*+{\bullet}="a15"
!{(-3.5,-1)}*+{\bullet}="a16"
!{(-5,0)}*+{\bullet}="a17"
!{(-6,0)}*+{\bullet}="a18"
!{(5,0)}*+{\circ}="b1"
!{(-5,1)}*+{\circ}="b2"
!{(-5,1.2)}*+{^0}
!{(-6,1)}*+{\circ}="b3"
!{(-6,-1)}*+{\circ}="b4"
"c1"-"a1" "a1"-"a5" "a1"-"a7" "a5"-"a2"
"a7"-"a2" "a5"-"a6" "a2"-"a3" "a7"-"a8"
"a3"-"a6" "a3"-"a8" "a6"-"a4" "a8"-"a4"
"c1"-"a9" "a9"-"a13" "a9"-"a15" "a13"-"a10"
"a10"-"a15" "a13"-"a14" "a10"-"a11" "a15"-"a16"
"a11"-"a14" "a11"-"a16" "a14"-"a12" "a16"-"a12"
"a12"-"a17" "a17"-"a18" "a4"-"b1" "a17"-"b2"
"a18"-"b3" "a18"-"b4" 
}
\end{displaymath}
Now, by gluing the points $\circ^0,\circ^1,\dots \circ^{k+1} $ counterclockwise outside in order, we get the dessin as desired. For the case $k$ is odd, the construction is the same. The method to prove the monodromy group corresponding to this dessin is primitive is similar to the method in case 1. See Remark \ref{R.k/4} for the detail.

\item 
($M=G_{13},\ n\in\lbrace \pm1/6,\pm1/3\rbrace+\mathbb{Z}$)

For $n=(1+2k)/6,n\in\mathbb{N}$, the corresponding list is $[1^2 2^{4k+2},\\(2k+4)^1 3^{2k}, 2^2 4^{2k}]$, and the dessin is given by 
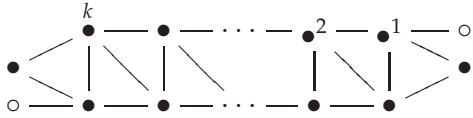
\begin{displaymath}
\xygraph{
!{<0cm,0cm>;<1cm,0cm>:<0cm,1cm>::}
!{(1,0)}*+{\bullet^2}="a1"
!{(2,0)}*+{\bullet^1}="a2"
!{(2,-1)}*+{\bullet}="a4"
!{(1,-1)}*+{\bullet}="a5"
!{(-1,0)}*+{\bullet}="a6"
!{(-2,0)}*+{\bullet}="a7"
!{(-2,0.2)}*+{^{k}}
!{(-1,-1)}*+{\bullet}="a9"
!{(0,0)}*+{\cdots}="c1"
!{(0,-1)}*+{\cdots}="c4"
!{(3,-0.5)}*+{\bullet}="a11"
!{(-2,-1)}*+{\bullet}="a10"
!{(-3,-0.5)}*+{\bullet}="a12"
!{(3,0)}*+{\circ}="b1"
!{(-3,-1)}*+{\circ}="b2"
"c1"-"a1" "c1"-"a6" 
"a1"-"a5" "a1"-"a2" "a2"-"a4" "a4"-"a5"
"a6"-"a7" "a6"-"a7" 
"a6"-"a9" "a7"-"a9" "a1"-"a4" 
"a9"-"a10" "b2"-"a10" "b1"-"a2"
"a10"-"a7" "a11"-"a2" "a11"-"a4"
"a12"-"a7" "a12"-"a10"
"c4"-"a5" "c4"-"a6" "c4"-"a9"
}
\end{displaymath}\\

For $n=k/3,k\in\mathbb{N}$, the corresponding list is $[1^2 2^{4k-1}, (2k+3)^1\\ 3^{2k-1}, 2^2 4^{2k-1}]$, and the dessin is given by
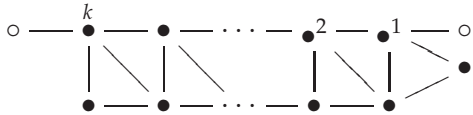
\begin{displaymath}
\xygraph{
!{<0cm,0cm>;<1cm,0cm>:<0cm,1cm>::}
!{(1,0)}*+{\bullet^2}="a1"
!{(2,0)}*+{\bullet^1}="a2"
!{(2,-1)}*+{\bullet}="a4"
!{(1,-1)}*+{\bullet}="a5"
!{(-1,0)}*+{\bullet}="a6"
!{(-2,0)}*+{\bullet}="a7"
!{(-2,0.2)}*+{^{k}}
!{(-1,-1)}*+{\bullet}="a9"
!{(0,0)}*+{\cdots}="c1"
!{(0,-1)}*+{\cdots}="c4"
!{(3,-0.5)}*+{\bullet}="a11"
!{(-2,-1)}*+{\bullet}="a10"
!{(3,0)}*+{\circ}="b1"
!{(-3,0)}*+{\circ}="b2"
"c1"-"a1" "c1"-"a6" 
"a1"-"a5" "a1"-"a2" "a2"-"a4" "a4"-"a5"
"a6"-"a7" "a6"-"a7" 
"a6"-"a9" "a7"-"a9" "a1"-"a4" 
"a9"-"a10"  "b1"-"a2"
"a10"-"a7" "a11"-"a2" "a11"-"a4"
"b2"-"a7" "c4"-"a5" "c4"-"a6" "c4"-"a9"
}
\end{displaymath}
where $\bullet$ represent the points mapped to 1; $\circ$ represent the points mapped to 0.
\item 
($M=G_{13},\ n\in\lbrace\pm1/8,\pm3/8\rbrace+\mathbb{Z}$)

For $n=(1+2k)/8,k\in\mathbb{N}$, the corresponding list is $[1^{3}2^{3k},3^{2k+1},\\2^1(2k+5)^{1}4^{k-1}]$, and the dessin is given by
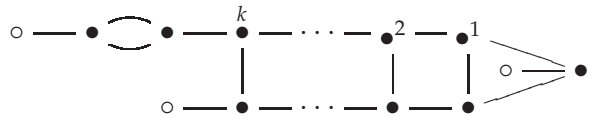
\begin{displaymath}
\xygraph{
!{<0cm,0cm>;<1cm,0cm>:<0cm,1cm>::}
!{(1,0)}*+{\bullet^2}="a1"
!{(2,0)}*+{\bullet^1}="a2"
!{(3.5,-0.5)}*+{\bullet}="a3"
!{(2,-1)}*+{\bullet}="a4"
!{(1,-1)}*+{\bullet}="a5"
!{(-1,0)}*+{\bullet}="a6"
!{(-1,0.2)}*+{^{k}}
!{(-1,-1)}*+{\bullet}="a7"
!{(0,0)}*+{\cdots}="c1"
!{(0,-1)}*+{\cdots}="c2"
!{(2.5,-0.5)}*+{\circ}="b1"
!{(-2,0)}*+{\bullet}="b2"
!{(-2,-1)}*+{\circ}="b3"
!{(-3,0)}*+{\bullet}="b4"
!{(-4,0)}*+{\circ}="b5"
"c1"-"a1" "c1"-"a6" "c2"-"a5" "c2"-"a7"
"a1"-"a5" "a1"-"a2" "a2"-"a4" "a4"-"a5"
"a2"-"a3" "a3"-"a4" "a3"-"b1"
"a6"-"a7" "a6"-"b2" "a7"-"b3" "b2"-@/_/"b4"
"b2"-@/^/"b4" "b4"-"b5"
}
\end{displaymath}
where $\bullet$ represent the points mapped to 1; $\circ$ represent the points mapped to 0.

\item 
$(M=G_{22},n\in\lbrace \pm1/6,\pm1/3\rbrace+\mathbb{Z},n>0)$

The proof of this case is similar to the case 4 of the proof of Theorem \ref{BM}. We will use the notation 
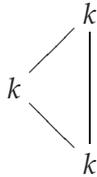
\begin{displaymath}
\xygraph{
!{<0cm,0cm>;<1cm,0cm>:<0cm,1cm>::}
!{(0,0)}*+{k}="a1"
!{(1,1)}*+{k}="a2"
!{(1,-1)}*+{k}="a3"
"a1"-"a2" "a3"-"a1" "a2"-"a3"
}
\end{displaymath}
to denote repeating the following picture $k$ times after it.
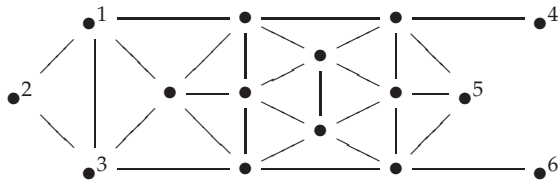
\begin{displaymath}
\xygraph{
!{<0cm,0cm>;<1cm,0cm>:<0cm,1cm>::}
!{(0,0)}*+{\bullet}="a1"
!{(0,1)}*+{\bullet}="a2"
!{(0,-1)}*+{\bullet}="a3"
!{(1,0.5)}*+{\bullet}="a4"
!{(1,-0.5)}*+{\bullet}="a5"
!{(2,0)}*+{\bullet}="a6"
!{(2,1)}*+{\bullet}="a7"
!{(2,-1)}*+{\bullet}="a8"
!{(-1,0)}*+{\bullet}="a9"
!{(-2,1)}*+{\bullet^1}="c1"
!{(-2,-1)}*+{\bullet^3}="c2"
!{(-3,0)}*+{\bullet^2}="c3"
!{(4,1)}*+{\bullet^4}="c4"
!{(3,0)}*+{\bullet^5}="c5"
!{(4,-1)}*+{\bullet^6}="c6"
"c4"-"a7" "c5"-"a7" "c5"-"a6" "c5"-"a8" "c6"-"a8" "a7"-"a2"
"a7"-"a4" "a6"-"a4" "a6"-"a5" "a8"-"a5" "a8"-"a3" "a4"-"a2"
"a4"-"a1" "a4"-"a5" "a5"-"a1" "a5"-"a3" "a6"-"a7" "a6"-"a8"
"a1"-"a2" "a1"-"a3" "a9"-"a2" "a9"-"a1" "a9"-"a3" "c1"-"a9"
"c2"-"a9" "c1"-"c2" "c1"-"c3" "c2"-"c3" "c1"-"a2" "c2"-"a3"
}
\end{displaymath}

For $n=7/6+k,k\in\mathbb{Z}_{\geq 0}$, the corresponding list is $[1^4 2^{33+30k},\\(10+6k)^1 3^{20+18k},5^{14+12k}]$, and the dessin is given by
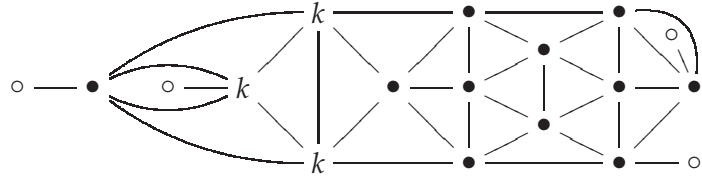
\begin{displaymath}
\xygraph{
!{<0cm,0cm>;<1cm,0cm>:<0cm,1cm>::}
!{(0,0)}*+{\bullet}="a1"
!{(0,1)}*+{\bullet}="a2"
!{(0,-1)}*+{\bullet}="a3"
!{(1,0.5)}*+{\bullet}="a4"
!{(1,-0.5)}*+{\bullet}="a5"
!{(2,0)}*+{\bullet}="a6"
!{(2,1)}*+{\bullet}="a7"
!{(2,-1)}*+{\bullet}="a8"
!{(-1,0)}*+{\bullet}="a9"
!{(-2,1)}*+{k}="c1"
!{(-2,-1)}*+{k}="c2"
!{(-3,0)}*+{k}="c3"
!{(-4,0)}*+{\circ}="b2"
!{(-5,0)}*+{\bullet}="a11"
!{(-6,0)}*+{\circ}="b1"
!{(3,0)}*+{\bullet}="c5"
!{(2.7,0.7)}*+{\circ}="b3"
!{(3,-1)}*+{\circ}="b4"
"c5"-"a7" "c5"-"a6" "c5"-"a8" "a7"-"a2"
"a7"-"a4" "a6"-"a4" "a6"-"a5" "a8"-"a5" "a8"-"a3" "a4"-"a2"
"a4"-"a1" "a4"-"a5" "a5"-"a1" "a5"-"a3" "a6"-"a7" "a6"-"a8"
"a1"-"a2" "a1"-"a3" "a9"-"a2" "a9"-"a1" "a9"-"a3" "c1"-"a9"
"c2"-"a9" "c1"-"c2" "c1"-"c3" "c2"-"c3" "c1"-"a2" "c2"-"a3"
"b2"-"c3" "a11"-@/^3mm/"c1" "a11"-@/^3mm/"c3" "a11"-"b1"
"a11"-@/_3mm/"c2" "a11"-@/_3mm/"c3"
"c5"-"b3" "b4"-"a8" "a7"-@/^0.5cm/"c5"
}
\end{displaymath}\\

For $n=4/3+k,k\in\mathbb{Z}_{\geq 0}$, the corresponding list is $[1^4 2^{43+30k},\\(11+6k)3^{23+18k},5^{16+12k}]$, and the dessin is given by 
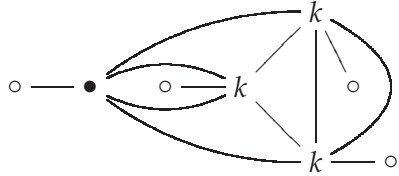
\begin{displaymath}
\xygraph{
!{<0cm,0cm>;<1cm,0cm>:<0cm,1cm>::}
!{(-2,1)}*+{k}="c1"
!{(-2,-1)}*+{k}="c2"
!{(-3,0)}*+{k}="c3"
!{(-4,0)}*+{\circ}="b2"
!{(-5,0)}*+{\bullet}="a11"
!{(-6,0)}*+{\circ}="b1"
!{(-1.5,0)}*+{\circ}="b3"
!{(-1,-1)}*+{\circ}="b4"
"c1"-"c2" "c1"-"c3" "c2"-"c3"
"b2"-"c3" "a11"-@/^3mm/"c1" "a11"-@/^3mm/"c3" "a11"-"b1"
"a11"-@/_3mm/"c2" "a11"-@/_3mm/"c3"
"c1"-"b3" "c2"-"b4" "c1"-@/^1cm/"c2"
}
\end{displaymath}\\

For $n=5/3+k,k\in\mathbb{Z}_{\geq 0}$, the corresponding list is $[1^4 2^{48+30k},\\(13+6k)^1 3^{29+18k},5^{20+12k}]$, and the dessin is given by 
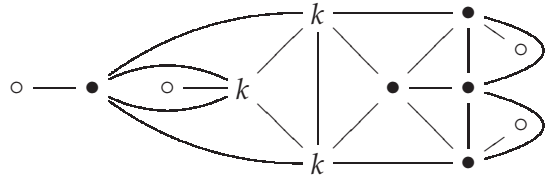
\begin{displaymath}
\xygraph{
!{<0cm,0cm>;<1cm,0cm>:<0cm,1cm>::}
!{(-2,1)}*+{k}="c1"
!{(-2,-1)}*+{k}="c2"
!{(-3,0)}*+{k}="c3"
!{(-4,0)}*+{\circ}="b2"
!{(-5,0)}*+{\bullet}="a11"
!{(-6,0)}*+{\circ}="b1"
!{(-1,0)}*+{\bullet}="a1"
!{(0,1)}*+{\bullet}="a2"
!{(0,0)}*+{\bullet}="a3"
!{(0,-1)}*+{\bullet}="a4"
!{(0.7,0.5)}*+{\circ}="b3"
!{(0.7,-0.5)}*+{\circ}="b4"
"c1"-"c2" "c1"-"c3" "c2"-"c3"
"b2"-"c3" "a11"-@/^3mm/"c1" "a11"-@/^3mm/"c3" "a11"-"b1"
"a11"-@/_3mm/"c2" "a11"-@/_3mm/"c3"
"c1"-"a1" "c2"-"a1" "c1"-"a2" "c2"-"a4"
"a1"-"a2" "a1"-"a3" "a1"-"a4" "a2"-"a3"
"a3"-"a4" "a2"-"b3" "a4"-"b4"
"a2"-@/^1cm/"a3" "a3"-@/^1cm/"a4"
}
\end{displaymath}\\

For $n=11/6+k,k\in\mathbb{Z}_{\geq 0}$, the corresponding list is $[1^4 2^{53+30k},\\(14+6k)^1 3^{32+18k},5^{22+12k}]$, and the dessin is given by 
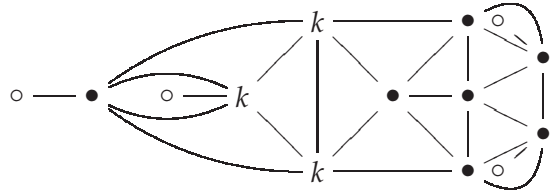
\begin{displaymath}
\xygraph{
!{<0cm,0cm>;<1cm,0cm>:<0cm,1cm>::}
!{(-2,1)}*+{k}="c1"
!{(-2,-1)}*+{k}="c2"
!{(-3,0)}*+{k}="c3"
!{(-4,0)}*+{\circ}="b2"
!{(-5,0)}*+{\bullet}="a11"
!{(-6,0)}*+{\circ}="b1"
!{(-1,0)}*+{\bullet}="a1"
!{(0,1)}*+{\bullet}="a2"
!{(0,0)}*+{\bullet}="a3"
!{(0,-1)}*+{\bullet}="a4"
!{(1,0.5)}*+{\bullet}="a5"
!{(1,-0.5)}*+{\bullet}="a6"
!{(0.4,1)}*+{\circ}="b3"
!{(0.4,-1)}*+{\circ}="b4"
"c1"-"c2" "c1"-"c3" "c2"-"c3"
"b2"-"c3" "a11"-@/^3mm/"c1" "a11"-@/^3mm/"c3" "a11"-"b1"
"a11"-@/_3mm/"c2" "a11"-@/_3mm/"c3"
"c1"-"a1" "c2"-"a1" "c1"-"a2" "c2"-"a4"
"a1"-"a2" "a1"-"a3" "a1"-"a4" "a2"-"a3"
"a3"-"a4" "a5"-"a2" "a5"-"a3" "a6"-"a3"
"a6"-"a4" "a5"-"b3" "a6"-"b4" "a5"-"a6"
"a2"-@/^5mm/"a5" "a6"-@/^5mm/"a4"
}
\end{displaymath}
where $\bullet$ represent the points mapped to $\infty$; $\circ$ represent the points mapped to 0.

\item 
($M=G_{22}$, $n\in\lbrace \pm{1/4}\rbrace+\mathbb{Z}$, $n>0$)

We use the same notation as in case 6.

For $n=1/4+k,k\in\mathbb{Z}_{\geq 0}$, the corresponding list is $[1^4 \\(2k+3)^1 2^{29k+5},3^{20k+5},5^{12k+3}]$. We first construct the dessin corresponding to the list $[1^5 2^{30k+5},(3k+6)^1 3^{19k+3},5^{12k+3}]$ and label $\circ$ which bound the outer face counterclockwisely. It is given by 
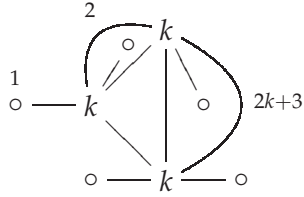
\begin{displaymath}
\xygraph{
!{<0cm,0cm>;<1cm,0cm>:<0cm,1cm>::}
!{(-1,0)}*+{k}="c1"
!{(-2,0.3)}*+{^1}
!{(-1,1.2)}*+{^2}
!{(1.5,0)}*+{^{2k+3}}
!{(0,1)}*+{k}="c2"
!{(0,-1)}*+{k}="c3"
!{(-0.5,0.8)}*+{\circ}="b1"
!{(-2,0)}*+{\circ}="b2"
!{(-1,-1)}*+{\circ}="b3"
!{(1,-1)}*+{\circ}="b4"
!{(0.5,0)}*+{\circ}="b5"
"c1"-"c2" "c2"-"c3" "c1"-"c3" "c1"-"b2"
"c3"-"b3" "c3"-"b4" "c2"-"b5" "c1"-@/^6mm/"c2"
"c2"-@/^1cm/"c3" "c1"-"b1" 
}
\end{displaymath}
Then gluing $\circ^1, \circ^4,\cdots,\circ^{3k+4}$ together can get our desired dessin. 

For $n=3/4+k,k\in\mathbb{Z}_{\geq 0}$, the corresponding list is $[1^5 2^{20+30k},\\(3k+9)^1 3^{19k+12},5^{12k+9}]$. We first construct the dessin corresponding to the list $[1^5 2^{30k+20},(3k+9)^1 3^{19k+12},5^{12k+9}]$ and label $\circ$ which bound the outer face counterclockwisely. It is given by 
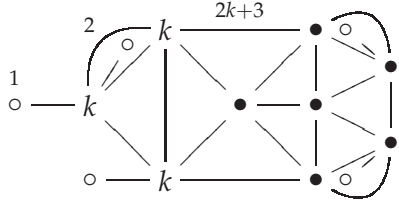
\begin{displaymath}
\xygraph{
!{<0cm,0cm>;<1cm,0cm>:<0cm,1cm>::}
!{(-2,1)}*+{k}="c1"
!{(-2,-1)}*+{k}="c2"
!{(-3,0)}*+{k}="c3"
!{(-4,0.3)}*+{^1}
!{(-3,1)}*+{^2}
!{(-1,1.2)}*+{^{2k+3}}
!{(-2.5,0.8)}*+{\circ}="b0"
!{(-4,0)}*+{\circ}="b1"
!{(-3,-1)}*+{\circ}="b2"
!{(-1,0)}*+{\bullet}="a1"
!{(0,1)}*+{\bullet}="a2"
!{(0,0)}*+{\bullet}="a3"
!{(0,-1)}*+{\bullet}="a4"
!{(1,0.5)}*+{\bullet}="a5"
!{(1,-0.5)}*+{\bullet}="a6"
!{(0.4,1)}*+{\circ}="b3"
!{(0.4,-1)}*+{\circ}="b4"
"c1"-"c2" "c1"-"c3" "c2"-"c3"
"c1"-"a1" "c2"-"a1" "c1"-"a2" "c2"-"a4"
"a1"-"a2" "a1"-"a3" "a1"-"a4" "a2"-"a3"
"a3"-"a4" "a5"-"a2" "a5"-"a3" "a6"-"a3"
"a6"-"a4" "a5"-"b3" "a6"-"b4" "a5"-"a6"
"a2"-@/^5mm/"a5" "a6"-@/^5mm/"a4"
"c3"-@/^5mm/"c1" "c3"-"b0" "c3"-"b1" "c2"-"b2"
}
\end{displaymath}
Then gluing $\circ^1,\circ^4, \cdots \circ^{3k+7}$ together can get our desired dessin.

\item 
$(M=G_{22},\ n\in\lbrace\pm1/10,\pm3/10,\pm1/5,\pm2/5\rbrace+\mathbb{Z},n>0)$

For $n=(1+2k)/10,k\in\mathbb{N}$, the corresponding list is $[1^4 2^{6k+1},3^{4k+2},\\(6+2k)^1 5^{2k}]$, and the dessin is given by 
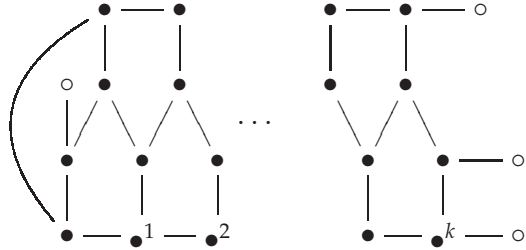
\begin{displaymath}
\xygraph{
!{<0cm,0cm>;<1cm,0cm>:<0cm,1cm>::}
!{(0,-0.5)}*+{\cdots}="c0"
!{(1,0)}*+{\bullet}="a1"
!{(2,0)}*+{\bullet}="a2"
!{(1,1)}*+{\bullet}="a3"
!{(2,1)}*+{\bullet}="a4"
!{(1.5,-1)}*+{\bullet}="a5"
!{(2.5,-1)}*+{\bullet}="a6"
!{(1.5,-2)}*+{\bullet}="a7"
!{(2.5,-2)}*+{\bullet^k}="a8"
!{(-1,0)}*+{\bullet}="a9"
!{(-1,1)}*+{\bullet}="a11"
!{(-2,0)}*+{\bullet}="a10"
!{(-2,1)}*+{\bullet}="a12"
!{(-1.5,-1)}*+{\bullet}="a13"
!{(-2.5,-1)}*+{\bullet}="a14"
!{(-1.5,-2)}*+{\bullet^1}="a15"
!{(-2.5,-2)}*+{\bullet}="a16"
!{(-2.5,0)}*+{\circ}="b1"
!{(3.5,-2)}*+{\circ}="b3"
!{(3,1)}*+{\circ}="b4"
!{(3.5,-1)}*+{\circ}="b5"
!{(-0.5,-1)}*+{\bullet}="a18"
!{(-0.5,-2)}*+{\bullet^2}="a19"
"a1"-"a3" "a3"-"a4" "a4"-"a2" "a1"-"a5" "a2"-"a5" "a5"-"a7"
"a7"-"a8" "a8"-"a6" "a6"-"a2" "a15"-"a16" "a14"-"a16" "a10"-"a14"
"a9"-"a11" "a11"-"a12" "a12"-"a10" "a9"-"a13" "a10"-"a13" "a13"-"a15"
"a12"-@/_1cm/"a16" "a14"-"b1"
"a8"-"b3" "a9"-"a18" "a18"-"a19" "a19"-"a15" "b4"-"a4" "b5"-"a6"
}
\end{displaymath}\\

For $n=(1+k)/5,k\in\mathbb{N}$, the corresponding list is $[1^4 2^{6k+4},3^{4k+4},\\(7+2k)^1 5^{2k+1}]$, and the dessin is given by 
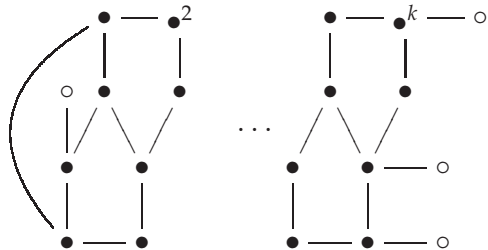
\begin{displaymath}
\xygraph{
!{<0cm,0cm>;<1cm,0cm>:<0cm,1cm>::}
!{(0,-0.5)}*+{\cdots}="c0"
!{(1,0)}*+{\bullet}="a1"
!{(2,0)}*+{\bullet}="a2"
!{(1,1)}*+{\bullet}="a3"
!{(2,1)}*+{\bullet^k}="a4"
!{(3,1)}*+{\circ}="b3"
!{(1.5,-1)}*+{\bullet}="a5"
!{(1.5,-2)}*+{\bullet}="a7"
!{(-1,0)}*+{\bullet}="a9"
!{(-1,1)}*+{\bullet^2}="a11"
!{(-2,0)}*+{\bullet}="a10"
!{(-2,1)}*+{\bullet}="a12"
!{(-1.5,-1)}*+{\bullet}="a13"
!{(-2.5,-1)}*+{\bullet}="a14"
!{(-1.5,-2)}*+{\bullet}="a15"
!{(-2.5,-2)}*+{\bullet}="a16"
!{(-2.5,0)}*+{\circ}="b1"
!{(2.5,-1)}*+{\circ}="b4"
!{(2.5,-2)}*+{\circ}="b5"
!{(0.5,-1)}*+{\bullet}="a18"
!{(0.5,-2)}*+{\bullet}="a19"
"a1"-"a3" "a3"-"a4" "a4"-"a2" "a1"-"a5" "a2"-"a5" "a5"-"a7"
 "a15"-"a16" "a14"-"a16" "a10"-"a14"
"a9"-"a11" "a11"-"a12" "a12"-"a10" "a9"-"a13" "a10"-"a13" "a13"-"a15"
"a12"-@/_1cm/"a16" "a14"-"b1" 
 "a1"-"a18" "a18"-"a19" "a19"-"a7" "a4"-"b3" "a5"-"b4" "a7"-"b5"
}
\end{displaymath}
The case $n=1/5,1/10$ is easy.
\end{enumerate}
We complete all the cases.
\end{proof}

\begin{remark} \label{R. Ex}
From now on, it seems that if $n$ being given and the corresponding table exists, then there exists dessin corresponding to it. Here, we give an example which has a table but no dessin corresponds to it. 

For $n=5/8$, by checking the local exponent carefully, we get the two possible tables:
\begin{table}[h]
\begin{tabular}{|l|l|c|c|c|c|c|l|}
\hline
 & $e_1$ & $e_2$ & $e_3$ & $\wp(a)$&$\infty$& \\ \hline
0 & 1 &1 & 1 &  &  & 6 points with multiplicity 2\\ \hline
1 &  & &  &  &  & 5 points with multiplicity 3\\ \hline
$\infty$ &&  &  & 9 & 2 & 1 point with multiplicity 4\\ \hline
\end{tabular}
\end{table}
\begin{table}[h]
\begin{tabular}{|l|l|c|c|c|c|c|l|}
\hline
 & $e_1$ & $e_2$ & $e_3$ & $\wp(a)$&$\infty$& \\ \hline
0 & 1 & &  &  &  & 7 points with multiplicity 2\\ \hline
1 &  & &  &  &  & 5 points with multiplicity 3\\ \hline
$\infty$ &&  2& 2 & 9 & 2 & no point with multiplicity 4\\ \hline
\end{tabular}
\end{table}\\
There exists dessin corresponding to the first table. It is given in the case 5 of the proof of the Theorem \ref{BM1}. However, no dessin corresponds to the second table. According to the case 3 of the proof of the Theorem \ref{BM1}, only one of $\gamma_i$ is the product of $2$-cycles as its projection in $S_4$. If there exists a dessin corresponding to the second table, three of the $\gamma_i$ will be the product of $2$-cycles, a contradiction. 
\end{remark}

\begin{remark}\label{R.k/4}
According to the case 1 of the proof of Theorem \ref{BM2}, we can construct a dessin corresponding to the list $[1^4 2^{6k+1},3^{4k+2},(4k+6)^1 4^{2k}]$. Here we give another dessin corresponding to the list
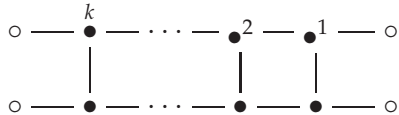
\begin{displaymath}
\xygraph{
!{<0cm,0cm>;<1cm,0cm>:<0cm,1cm>::}
!{(1,0)}*+{\bullet^2}="a1"
!{(2,0)}*+{\bullet^1}="a2"
!{(2,-1)}*+{\bullet}="a4"
!{(1,-1)}*+{\bullet}="a5"
!{(-1,0)}*+{\bullet}="a6"
!{(-1,0.2)}*+{^{k}}
!{(-1,-1)}*+{\bullet}="a7"
!{(0,0)}*+{\cdots}="c1"
!{(0,-1)}*+{\cdots}="c2"
!{(-2,0)}*+{\circ}="b2"
!{(-2,-1)}*+{\circ}="b3"
!{(3,0)}*+{\circ}="b4"
!{(3,-1)}*+{\circ}="b5"
"c1"-"a1" "c1"-"a6" "c2"-"a5" "c2"-"a7"
"a1"-"a5" "a1"-"a2" "a2"-"a4" "a4"-"a5"
"a6"-"a7" "a6"-"b2" "a7"-"b3" "a2"-"b4"
"a4"-"b5"
}
\end{displaymath}
The difference between the two dessins is that the one given in the proof has primitive monodromy group and the one given here hasn't. It is easy to show the monodromy group of the dessin given here is imprimitive. We skip the proof. It remains to show that the monodromy group of the dessin given in the proof is primitive. To show it, we only need to focus on the left part of the dessin 
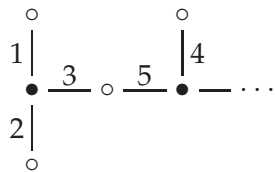
\begin{displaymath}
\xygraph{
!{<0cm,0cm>;<1cm,0cm>:<0cm,1cm>::}
!{(0,0)}*+{\bullet}="a1"
!{(0,1)}*+{\circ}="b1"
!{(0,-1)}*+{\circ}="b2"
!{(1,0)}*+{\circ}="bb"
!{(2,0)}*+{\bullet}="a2"
!{(2,1)}*+{\circ}="b3"
!{(3,0)}*+{\cdots}="c1"
!{(-0.2,0.5)}*+{1}
!{(-0.2,-0.5)}*+{2}
!{(0.5,0.2)}*+{3}
!{(1.5,0.2)}*+{5}
!{(2.2,0.5)}*+{4}
"a1"-"b1" "a1"-"b2" "a1"-"bb"
"bb"-"a2" "a2"-"b3" "a2"-"c1"
}
\end{displaymath}
In this case, we can observe that 
$$\sigma_0=(1)(2)(3\ 5)\dots,\ \sigma_1=(123)(45\dots)\dots.$$ 
If this monodromy group is imprimitive, then let $\Sigma_0=(1)(23),\Sigma_1=(123),\\\Sigma_{\infty}=(12)(3)$ which realizes $H_{1/2,1/3,1/4}$ as a pullback of $H_{1/2,1/2,1/4}$. Also, let $P_1,P_2,P_3$ be a partition such that $\sigma_{\alpha}(P_i)=P_{\Sigma_{\alpha}(i)}$. We claim that the partition does not exists. It means the monodromy group is primitive. It followed by some observations. First, the number 1, 2 and 4 must lies in the same partition say $P_1$, since the are fixed by $\sigma_0$. Second, 3 lies in $P_1$, or we will have $\sigma_1(P_1)\neq P_1$ and $\sigma_1(P_1)\cap P_1 \neq \lbrace \emptyset\rbrace$. Third, 5 lies in $P_1$ since $P_1$ is fixed under $\sigma_0$. Now, we have $P_1$ is fixed under $\sigma_0,\sigma_1,\sigma_{\infty}$. It is impossible and proves the claim. 
\end{remark}

\section{Finite monodromy problem for the generalized Lam\'e euqations with 3 singular points}

In this section, we focus on the finite monodromy problem of equation (\ref{C-ale-3pts}). For convenience, we denote it $L_{n_0,n_1,g_2,g_3,A,B}y=0$. We assume $n_0,n_1\geq -1/2$ as before. The following theorem gives all the possible pairs $(n_0, n_1)$ such that $L_{n_0,n_1,g_2,g_3,A,B}y=0$ may have dihedral projective monodromy group.

\begin{theorem}\label{n0n1}
If the equation $L_{n_0,n_1,g_2,g_3,A,B}y=0$ has finite projective monodromy group isomorphic to $D_N$ for some $N\geq 2$, then the pairs $(n_0,n_1)$ lies in one of the following lists:
\begin{enumerate}[label=(\alph*)]
\item 
$n_0\in\mathbb{Z}/2$, $n_0\geq 1$ and $n_1\in\mathbb{Q}$,
\item 
$n_0\in\mathbb{Q}$, $n_1\in \lbrace \pm1/4\rbrace+\mathbb{Z}$ and $n_1\geq 3/4$.
\end{enumerate}
Conversely, for given $n_0,n_1$ above, there exists equation $L_{n_0,n_1,g_2,g_3,A,B}y=0$ which has finite projective monodromy group isomorphic to $D_N$ for some $N\geq 2$.
\end{theorem}

\begin{proof}

The cases
\begin{enumerate}
\item 
$n_0\in\mathbb{Z}/2$ and $n_1\in\mathbb{Z}/2$;
\item 
$n_0\in\lbrace 1/2\rbrace+\mathbb{Z}$ and $n_1\in\lbrace \pm1/4\rbrace+\mathbb{Z}$.
\end{enumerate}
will be treated in Theorem \ref{n0n1-1} and Theorem \ref{n0n1-2}. Then in the following proof we exclude the 2 cases.

Recall that the local exponent of the equation $L_{n_0,n_1,g_2,g_3,A,B}y=0$ is given by
$$\left\lbrace
\begin{array}{ccccccc}
e_1 & e_2 & e_3 & \wp(a_1)&  \infty\\
0 & 0 & 0 & -n_1  & -n_0/2 \\
\frac{1}{2} & \frac{1}{2} & \frac{1}{2} & n_1+1  & (n_0+1)/2
\end{array}
\right\rbrace.$$
Notice that $n_0,n_1$ in lists ($a$) and ($b$) correspond to the same local exponent by exchanging $\wp(a_1)$ with $\infty$. Hence, it is enough to consider (a). 

To prove the necessary condition of $n_0$ and $n_1$, let $n_0=p/q,n_1=p_1/q_1$ and $PM\cong D_N$. By checking the local exponent carefully, we can get the following table:

\begin{table}[h]
\begin{tabular}{|c|c|c|c|c|c|c|c|}
\hline
 & $e_1$ & $e_2$ & $e_3$ & $\wp(a)$&$\infty$& \\ \hline
0 &  & &  &  &  & points with degree 2\\ \hline
1 &  & &  & $2p/q+1$ &  & points with degree 2\\ \hline
$\infty$ &  &  &  &  & $Np_1/q_1+N$ & points with degree $N$\\ \hline
\end{tabular}
\end{table}

Notice that at least one of the points $\wp(a),\infty$ is mapped to $\infty$. Also, $\wp(a)$ and $\infty$ can not be mapped to $\infty$ simultaneously, since otherwise the sum of the ramification index will greater than the degree of the pullback function which leads to a contradiction. Furthermore, we can assume that $\infty$ is mapped to $\infty$ by exchanging $\wp(a)$ and $\infty$ if necessary. Now, the ramification index $2p/q+1\in\mathbb{Z}$ gives the condition $n_0\in\mathbb{Z}/2$. At last, the degree of the pullback function is greater than or equal to the degree of the ramification index gives the equation: $N(n_0+2n_1)\geq Nn_2+N$. It implies $n_0\geq 1$. This completes the necessary part. 

Conversely, let $n_0=p/2,n_1=k/N$ with g.c.d.$(k,N)=1$ and $N\geq 3$ (since we have excluded the case that not all $n_0,n_1$ contain in $\mathbb{Z}/2$). We claim that in this case $L_{n_0,n_1,g_2,g_3,A,B}$ is the pullback of the hypergeometric operator $H_{1/2,1/2,1/N'}$ whenever $4N\mid N'$. The method is simply to construct the corresponding dessin. Let $N'=mN,\ 4\mid m$. By checking the degree carefully, we get the following tables. 

For $p$ is even

\begin{table}[h]
\begin{tabular}{|c|c|c|c|c|c|c|c|}
\hline
 & $e_1$ & $e_2$ & $e_3$ & $\wp(a)$&$\infty$& \\ \hline
0 & 0,1 & 0,1 & 0,1 &  &  & points with degree 2\\ \hline
1 & 0,1 & 0,1 & 0,1 & $p+1$ &  & points with degree 2\\ \hline
$\infty$ &  &  &  &  & $m(2k+N)$ & $p/2-1$ points with degree $mN$\\ \hline
\end{tabular}
\end{table}

For $p$ is odd

\begin{table}[h]
\begin{tabular}{|c|c|c|c|c|c|c|c|}
\hline
 & $e_1$ & $e_2$ & $e_3$ & $\wp(a)$&$\infty$& \\ \hline
0 & 0,1 & 0,1 &  &  &  & points with degree 2\\ \hline
1 & 0,1 & 0,1 &  & $p+1$ &  & points with degree 2\\ \hline
$\infty$ &  &  & $mN/2$ &  & $m(2k+N)$ & $(p-3)/2$ points with degree $mN$\\ \hline
\end{tabular}
\end{table}

In all cases, there exist dessin of the form given below corresponding to the table. 
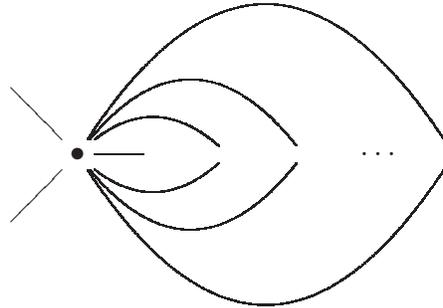
\begin{displaymath}
\xygraph{
!{<0cm,0cm>;<1cm,0cm>:<0cm,1cm>::}
!{(0,0)}*+{\bullet}="a1"
!{(1,0)}*+{}="l1"
!{(-1,1)}*+{}="l2"
!{(-1,-1)}*+{}="l3"
!{(2,0)}*+{}="a2"
!{(3,0)}*+{}="a3"
!{(4,0)}*+{\cdots}
!{(5,0)}*+{}="a4"
"a1"-@/_0.5cm/"a2" "a1"-@/^0.5cm/"a2"
"a1"-@/_1cm/"a3" "a1"-@/^1cm/"a3"
"a1"-@/_2cm/"a4" "a1"-@/^2cm/"a4"
"a1"-"l1" "a1"-"l2" "a1"-"l3"
}
\end{displaymath}
The dessin contains $p/2-1({\rm resp.}/ (p-1)/2)$ cycles and $3$ (resp.~$2$) lines if $n$ is even (resp.~odd). It is easy to show that the dessin of this form exists, hence we complete the proof.
\end{proof}

\begin{theorem} \label{n0n1-1}
Suppose $L_{n_0,n_1,g_2,g_3,A,B}y=0$ has finite projective monodromy group $PM$. Let $n_0\in\lbrace 1/2\rbrace+\mathbb{Z}$ and $n_1\in\mathbb{Z}/2$. Then $PM\cong K_4$.
\end{theorem}
This is just the special case of the Corollary \ref{C-PM-Coro}. We skip the proof.

\begin{theorem}\label{n0n1-2}
As the same notation in Theorem \ref{n0n1-1}. Let $n_0,n_1$ lie in one of the following list:
\begin{enumerate}
\item $n_0\in\mathbb{Z}$ and $n_1\in\mathbb{Z}/2$;
\item $n_0\in\lbrace 1/2\rbrace+\mathbb{Z}$ and $n_1\in\lbrace\pm1/4\rbrace+\mathbb{Z}$.
\end{enumerate}
Then we have $PM\cong D_N$ for some $N\in\mathbb{N}_{\geq 2}$.
\end{theorem}
\begin{proof}
The list (1) is the special case of the Corollary \ref{C-PM-Coro}. For the list (2), notice that they have the same local exponent by exchanging the point $\wp(a_1)$ with $\infty$. We finish the proof.
\end{proof}

In the following content, we will consider the cases that $L_{n_0,n_1}$ has projective monodromy which is isomorphic to $A_4,S_4$ or $A_5$. Hence, we assume $n_0,n_1$ do not lie in the list given in Theorem \ref{n0n1-1} and Theorem \ref{n0n1-2}. Since, in this case, $L_{n_0,n_1,g_2,g_3,A,B}y=0$ only has dihedral group as its projective monodromy.

\begin{theorem} \label{n0n1-3}
Suppose $L_{n_0,n_1,g_2,g_3,A,B}y=0$ has finite projective monodromy group $PM$.
\begin{enumerate}
\item If $PM\cong A_4$, then $n_0=(1+2k)/6$, $n_1=l/6$ with $k\equiv 0,2 \pmod{3}$, $k+l\geq 4$ and $k+l\equiv 1 \pmod{3}$;
\item If $PM\cong S_4$, then $n_0\in (\mathbb{Z}/4)\cup(\mathbb{Z}/6)\setminus (\mathbb{Z}/3) $ and $n_1\in (\mathbb{Z}/8)\cup(\mathbb{Z}/6)$;
\item If $PM\cong A_5$, then $n_0\in (\mathbb{Z}/6)\cup(\mathbb{Z}/10)\setminus ((\mathbb{Z}/3)\cup (\mathbb{Z}/5) )$ and $n_1\in(\mathbb{Z}/4)\cup(\mathbb{Z}/6)\cup(\mathbb{Z}/10)$.
\end{enumerate}
\end{theorem}

\begin{proof}
The proof is again by checking the degree and show that there exists a table corresponding to it. 
\end{proof}

The degree checking can always give the necessary condition of $n_i$ of equation (\ref{C-Alg-Lame}). We say that $n_i$ satisfy the degree checking if the corresponding table exists. Hence the important problem is to show the list is sufficient. The next theorem completes this problem for the list given in Theorem \ref{n0n1-3}.

\begin{theorem} \label{3-pts M}
Suppose $n_0$ and $n_1$ do not lie in the list given in Theorem \ref{n0n1-1} and Theorem \ref{n0n1-2}. Let $n_0,n_1$ in Theorem \ref{n0n1-3} being given and satisfy the degree checking, then there exists equation $L_{n_0,n_1,g_2,g_3,A,B}y=0$ with corresponding projective monodromy group.
\end{theorem}

\begin{proof}
We divide the problem into several cases. Here are part the cases.
\begin{enumerate}
\item 
$PM\cong A_4$ with $n_0\in\lbrace \pm1/6\rbrace+\mathbb{Z}$ and $n_1\in\lbrace\pm1/3,\pm1/6\rbrace+\mathbb{Z}$;
\item 
$PM\cong S_4$ with $n_0\in\mathbb{Z}$ and $n_1\in\lbrace\pm1/4\rbrace+\mathbb{Z}$;
\item 
$PM\cong S_4$ with $n_0\in\mathbb{Z}$ and $n_1\in\lbrace\pm1/8,\pm3/8\rbrace+\mathbb{Z}$;
\item 
$PM\cong S_4$ with $n_0\in\mathbb{Z}$ and $n_1\in\lbrace\pm1/6,\pm1/3\rbrace+\mathbb{Z}$;
\item 
$PM\cong S_4$ with $n_0\in\lbrace 1/2\rbrace+\mathbb{Z}$ and $n_1\in\lbrace\pm1/8,\pm3/8\rbrace+\mathbb{Z}$;
\item 
$PM\cong S_4$ with $n_0\in\lbrace 1/2\rbrace+\mathbb{Z}$ and $n_1\in\lbrace\pm1/6,\pm1/3\rbrace+\mathbb{Z}$;
\item 
$PM\cong S_4$ with $n_0\in\lbrace \pm1/4\rbrace+\mathbb{Z}$ and $n_1\in\lbrace\pm1/8,\pm3/8\rbrace+\mathbb{Z}$;
\item 
$PM\cong S_4$ with $n_0\in\lbrace \pm1/4\rbrace+\mathbb{Z}$ and $n_1\in\lbrace\pm1/6,\pm1/3\rbrace+\mathbb{Z}$;
\item 
$PM\cong A_5$ with $n_0\in\mathbb{Z}$ and $n_1\in\lbrace\pm1/4\rbrace+\mathbb{Z}$;
\item $PM\cong A_5$ with $n_0\in\mathbb{Z}$ and $n_1\in\lbrace\pm1/3,\pm1/6\rbrace+\mathbb{Z}$;
\item 
$PM\cong A_5$ with $n_0\in\mathbb{Z}$ and $n_1\in (\mathbb{Z}/10)\setminus(\mathbb{Z}/2)$;
\item 
$PM\cong A_5$ with $n_0\in\lbrace 1/2\rbrace+\mathbb{Z}$ and $n_1\in\lbrace\pm1/3,\pm1/6\rbrace+\mathbb{Z}$;
\item 
$PM\cong A_5$ with $n_0\in\lbrace 1/2\rbrace+\mathbb{Z}$ and $n_1\in (\mathbb{Z}/10)\setminus(\mathbb{Z}/2)$;
\item 
$PM\cong A_5$ with $n_0\in\lbrace \pm1/6\rbrace+\mathbb{Z}$ and $n_1\in\lbrace\pm1/3,\pm1/6\rbrace+\mathbb{Z}$;
\item 
$PM\cong A_5$ with $n_0\in\lbrace\pm1/6\rbrace+\mathbb{Z}$ and $n_1\in (\mathbb{Z}/10)\setminus(\mathbb{Z}/2)$;
\item 
$PM\cong A_5$ with $n_0\in\lbrace\pm1/10,3/10\rbrace+\mathbb{Z}$ and $n_1\in (\mathbb{Z}/10)\setminus(\mathbb{Z}/2)$.
\end{enumerate}
We do not write down all the cases, since we can cover all the remaining cases by exchanging the point $\wp(a_1)$ with $\infty$.

Notice that cases (2), (3), (4), (5), (9), (10) and (11) may have dihedral group as its projective monodromy group. In these cases, we need to construct the dessin and check whether the corresponding monodromy group is primitive. For other cases, showing the existence of dessin is enough.
\begin{enumerate}[label=\textbf{case (\arabic*) :}]

\item 
($PM\cong A_4$ with $n_0\in\lbrace \pm1/6\rbrace+\mathbb{Z}$, $n_1\in\lbrace\pm1/3,\pm1/6\rbrace+\mathbb{Z}$)

Let $n_0=(2k+1)/6$ and $n_2=l/6$ with $k+l\geq 4$, $k+l\equiv 1({\rm mod}3)$ and $k\equiv 0,2({\rm mod}3)$. 

For $k=-1,l=5+3m$, the corresponding list is $[1^3 2^{3+3m},1^1 \\(8+3m)^1 3^m,3^{3+2m}]$. If $k$ is even, the dessin is given by
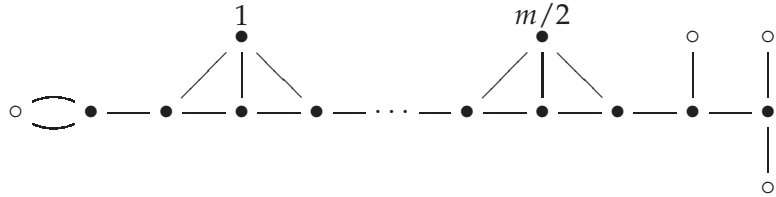
\begin{displaymath}
\xygraph{
!{<0cm,0cm>;<1cm,0cm>:<0cm,1cm>::}
!{(0,0)}*+{\cdots}="c"
!{(1,0)}*+{\bullet}="a1"
!{(2,0)}*+{\bullet}="a2"
!{(3,0)}*+{\bullet}="a3"
!{(2,1)}*+{\bullet}="a4"
!{(2,1.3)}*+{m/2}
!{(-1,0)}*+{\bullet}="a5"
!{(-2,0)}*+{\bullet}="a6"
!{(-3,0)}*+{\bullet}="a7"
!{(-2,1)}*+{\bullet}="a8"
!{(-2,1.3)}*+{1}
!{(4,0)}*+{\bullet}="a9"
!{(5,0)}*+{\bullet}="a10"
!{(-4,0)}*+{\bullet}="a11"
!{(4,1)}*+{\circ}="b1"
!{(5,1)}*+{\circ}="b2"
!{(5,-1)}*+{\circ}="b3"
!{(-5,0)}*+{\circ}="b4"
"a1"-"a2" "a2"-"a3" "a4"-"a1" "a4"-"a2"
"a4"-"a3" "a1"-"c" "c"-"a5" "a5"-"a6"
"a6"-"a7" "a8"-"a5" "a8"-"a6" "a8"-"a7"
"a3"-"a9" "a9"-"a10" "a7"-"a11" "a9"-"b1"
"a10"-"b2" "a10"-"b3" "a11"-@/_/"b4" "a11"-@/^/"b4"
}
\end{displaymath}\\
If $k$ is odd, the dessin is given by
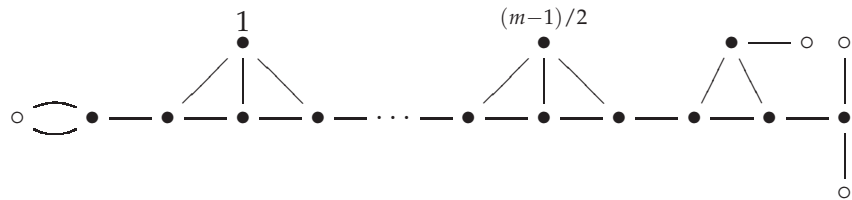
\begin{displaymath}
\xygraph{
!{<0cm,0cm>;<1cm,0cm>:<0cm,1cm>::}
!{(0,0)}*+{\cdots}="c"
!{(1,0)}*+{\bullet}="a1"
!{(2,0)}*+{\bullet}="a2"
!{(3,0)}*+{\bullet}="a3"
!{(2,1)}*+{\bullet}="a4"
!{(2,1.3)}*+{^{(m-1)/2}}
!{(-1,0)}*+{\bullet}="a5"
!{(-2,0)}*+{\bullet}="a6"
!{(-3,0)}*+{\bullet}="a7"
!{(-2,1)}*+{\bullet}="a8"
!{(-2,1.3)}*+{1}
!{(4,0)}*+{\bullet}="a9"
!{(5,0)}*+{\bullet}="a10"
!{(6,0)}*+{\bullet}="A1"
!{(4.5,1)}*+{\bullet}="A2"
!{(-4,0)}*+{\bullet}="a11"
!{(5.5,1)}*+{\circ}="b1"
!{(6,1)}*+{\circ}="b2"
!{(6,-1)}*+{\circ}="b3"
!{(-5,0)}*+{\circ}="b4"
"a1"-"a2" "a2"-"a3" "a4"-"a1" "a4"-"a2"
"a4"-"a3" "a1"-"c" "c"-"a5" "a5"-"a6"
"a6"-"a7" "a8"-"a5" "a8"-"a6" "a8"-"a7"
"a3"-"a9" "a9"-"a10" "a7"-"a11" "A2"-"b1"
"A1"-"b2" "A1"-"b3" "a11"-@/_/"b4" "a11"-@/^/"b4"
"A2"-"a9" "A2"-"a10" "A1"-"a10"
}
\end{displaymath}

Now, without loss of generality, we assume $l\geq k\geq 2$. For $k=2+3m$, we first construct the dessin
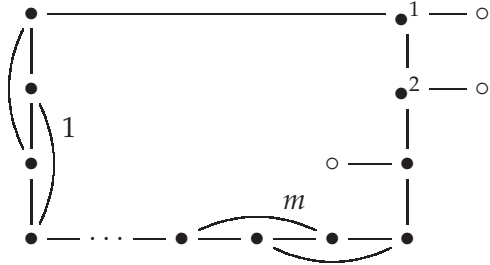
\begin{displaymath}
\xygraph{
!{<0cm,0cm>;<1cm,0cm>:<0cm,1cm>::}
!{(0,0)}*+{\cdots}="c"
!{(1,0)}*+{\bullet}="a1"
!{(2,0)}*+{\bullet}="a2"
!{(3,0)}*+{\bullet}="a3"
!{(4,0)}*+{\bullet}="a4"
!{(-0.5,1.5)}*+{1}
!{(2.5,0.5)}*+{m}
!{(4,1)}*+{\bullet}="a9"
!{(4,2)}*+{\bullet^2}="a10"
!{(4,3)}*+{\bullet^1}="a11"
!{(-1,0)}*+{\bullet}="a5"
!{(-1,1)}*+{\bullet}="a6"
!{(-1,2)}*+{\bullet}="a7"
!{(-1,3)}*+{\bullet}="a8"
!{(3,1)}*+{\circ}="b1"
!{(5,2)}*+{\circ}="b2"
!{(5,3)}*+{\circ}="b3"
"c"-"a1" "c"-"a5" "a1"-"a2" "a2"-"a3"
"a3"-"a4" "a5"-"a6" "a6"-"a7" "a7"-"a8"
"a4"-"a9" "a9"-"a10" "a10"-"a11" "a8"-"a11"
"a9"-"b1" "a10"-"b2" "a11"-"b3" "a1"-@/^3mm/"a3"
"a2"-@/_3mm/"a4" "a8"-@/_3mm/"a6" "a7"-@/^3mm/"a5"
}
\end{displaymath}
Then replace the right part 
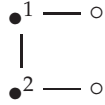
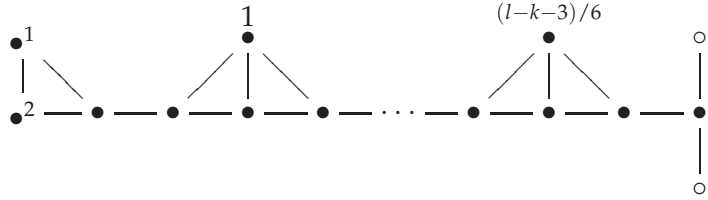
\begin{displaymath}
\xygraph{
!{<0cm,0cm>;<1cm,0cm>:<0cm,1cm>::}
!{(0,0)}*+{\bullet^2}="a1"
!{(0,1)}*+{\bullet^1}="a2"
!{(1,0)}*+{\circ}="b1"
!{(1,1)}*+{\circ}="b2"
"a1"-"a2" "a1"-"b1" "a2"-"b2"
}
\end{displaymath}
to
\begin{displaymath}
\xygraph{
!{<0cm,0cm>;<1cm,0cm>:<0cm,1cm>::}
!{(0,0)}*+{\cdots}="c"
!{(1,0)}*+{\bullet}="a1"
!{(2,0)}*+{\bullet}="a2"
!{(3,0)}*+{\bullet}="a3"
!{(2,1)}*+{\bullet}="a4"
!{(2,1.3)}*+{^{(l-k-3)/6}}
!{(-1,0)}*+{\bullet}="a5"
!{(-2,0)}*+{\bullet}="a6"
!{(-3,0)}*+{\bullet}="a7"
!{(-2,1)}*+{\bullet}="a8"
!{(-2,1.3)}*+{1}
!{(4,0)}*+{\bullet}="a9"
!{(-4,0)}*+{\bullet}="a11"
!{(-5,0)}*+{\bullet^2}="a12"
!{(-5,1)}*+{\bullet^1}="a13"
!{(4,1)}*+{\circ}="b1"
!{(4,-1)}*+{\circ}="b2"
"a1"-"a2" "a2"-"a3" "a4"-"a1" "a4"-"a2"
"a4"-"a3" "a1"-"c" "c"-"a5" "a5"-"a6"
"a6"-"a7" "a8"-"a5" "a8"-"a6" "a8"-"a7"
"a3"-"a9" "a7"-"a11" "a9"-"b1"
"a9"-"b2" "a11"-"a12" "a12"-"a13" "a11"-"a13"
}
\end{displaymath}
,when $6|(l-k-3)$; to
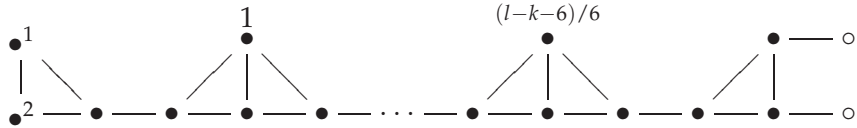
\begin{displaymath}
\xygraph{
!{<0cm,0cm>;<1cm,0cm>:<0cm,1cm>::}
!{(0,0)}*+{\cdots}="c"
!{(1,0)}*+{\bullet}="a1"
!{(2,0)}*+{\bullet}="a2"
!{(3,0)}*+{\bullet}="a3"
!{(2,1)}*+{\bullet}="a4"
!{(2,1.3)}*+{^{(l-k-6)/6}}
!{(-1,0)}*+{\bullet}="a5"
!{(-2,0)}*+{\bullet}="a6"
!{(-3,0)}*+{\bullet}="a7"
!{(-2,1)}*+{\bullet}="a8"
!{(-2,1.3)}*+{1}
!{(4,0)}*+{\bullet}="a9"
!{(5,0)}*+{\bullet}="a14"
!{(5,1)}*+{\bullet}="a15"
!{(-4,0)}*+{\bullet}="a11"
!{(-5,0)}*+{\bullet^2}="a12"
!{(-5,1)}*+{\bullet^1}="a13"
!{(6,0)}*+{\circ}="b1"
!{(6,1)}*+{\circ}="b2"
"a1"-"a2" "a2"-"a3" "a4"-"a1" "a4"-"a2"
"a4"-"a3" "a1"-"c" "c"-"a5" "a5"-"a6"
"a6"-"a7" "a8"-"a5" "a8"-"a6" "a8"-"a7"
"a3"-"a9" "a7"-"a11"
"a11"-"a12" "a12"-"a13" "a11"-"a13"
"a14"-"b1" "a15"-"b2" "a9"-"a14" "a9"-"a15" "a14"-"a15"
}
\end{displaymath}
,when $6|(l-k)$. For $k=3m$, the construction is similar. 

\item 
($PM\cong S_4$ with $n_0\in\mathbb{Z}$, $n_1\in\lbrace\pm1/4\rbrace+\mathbb{Z}$)

For $n_0=k$ and $n_1=(2l+1)/4$, the corresponding list is $[1^3 \\(2k+1)^1 2^{5k+6l+1},3^{2+4k+4l},(6+4l)^1 4^{2l+3k}]$. We first construct the dessin corresponding to the list $[1^4 2^{6k+6l+1},3^{2+4k+4l},(6+4k+4l)^1 4^{2k+2l}]$. It is given by the case 1 of the proof of Theorem \ref{BM2}. We construct the dessin for $k+l$ is even and label the $\circ$ which bound the outer face counterclockwise. It is given by
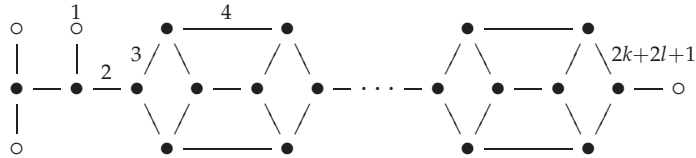
\begin{displaymath}
\xygraph{
!{<0cm,0cm>;<0.8cm,0cm>:<0cm,0.8cm>::}
!{(0,0)}*+{\cdots}="c1"
!{(1,0)}*+{\bullet}="a1"
!{(2,0)}*+{\bullet}="a2"
!{(3,0)}*+{\bullet}="a3"
!{(4,0)}*+{\bullet}="a4"
!{(1.5,1)}*+{\bullet}="a5"
!{(3.5,1)}*+{\bullet}="a6"
!{(1.5,-1)}*+{\bullet}="a7"
!{(3.5,-1)}*+{\bullet}="a8"
!{(-1,0)}*+{\bullet}="a9"
!{(-2,0)}*+{\bullet}="a10"
!{(-3,0)}*+{\bullet}="a11"
!{(-4,0)}*+{\bullet}="a12"
!{(-1.5,1)}*+{\bullet}="a13"
!{(-3.5,1)}*+{\bullet}="a14"
!{(-1.5,-1)}*+{\bullet}="a15"
!{(-3.5,-1)}*+{\bullet}="a16"
!{(-5,0)}*+{\bullet}="a17"
!{(-6,0)}*+{\bullet}="a18"
!{(5,0)}*+{\circ}="b1"
!{(-5,1)}*+{\circ}="b2"
!{(-6,1)}*+{\circ}="b3"
!{(-6,-1)}*+{\circ}="b4"
!{(-5,1.2)}*+{^1}
!{(-4.5,0.2)}*+{^2}
!{(-4,0.5)}*+{^3}
!{(-2.5,1.2)}*+{^4}
!{(4.6,0.5)}*+{^{2k+2l+1}}
"c1"-"a1" "a1"-"a5" "a1"-"a7" "a5"-"a2"
"a7"-"a2" "a5"-"a6" "a2"-"a3" "a7"-"a8"
"a3"-"a6" "a3"-"a8" "a6"-"a4" "a8"-"a4"
"c1"-"a9" "a9"-"a13" "a9"-"a15" "a13"-"a10"
"a10"-"a15" "a13"-"a14" "a10"-"a11" "a15"-"a16"
"a11"-"a14" "a11"-"a16" "a14"-"a12" "a16"-"a12"
"a12"-"a17" "a17"-"a18" "a4"-"b1" "a17"-"b2"
"a18"-"b3" "a18"-"b4" 
}
\end{displaymath}
Now, it is easy too see that by gluing $\circ^1,\circ^5,\dots,\circ^{4k+1}$ together, we get the dessin as desired. For $k+l$ is odd, the construction is similar. The method to prove whether the monodromy group corresponding to this dessin is primitive is the same as the method given in Remark \ref{R.k/4}.

\item 
($PM\cong S_4$ with $n_0\in\mathbb{Z}$, $n_1\in\lbrace\pm1/8,\pm3/8\rbrace+\mathbb{Z}$)

We first prove the dessin given in case 5 of the proof of Theorem \ref{BM2} has primitive monodromy group. The dessin is given by 
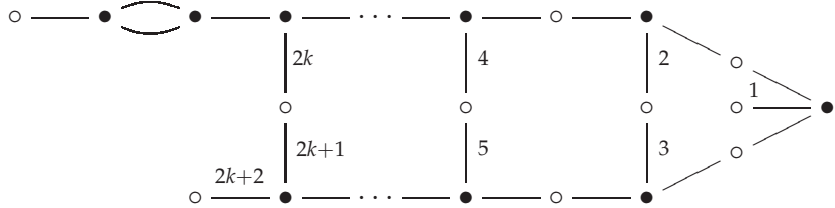
\begin{displaymath}
\xygraph{
!{<0cm,0cm>;<1.2cm,0cm>:<0cm,1.2cm>::}
!{(0,0)}*+{\cdots}="c1"
!{(0,2)}*+{\cdots}="c2"
!{(-1,0)}*+{\bullet}="a1"
!{(1,0)}*+{\bullet}="a2"
!{(3,0)}*+{\bullet}="a3"
!{(5,1)}*+{\bullet}="a4"
!{(3,2)}*+{\bullet}="a5"
!{(1,2)}*+{\bullet}="a6"
!{(-1,2)}*+{\bullet}="a7"
!{(-2,2)}*+{\bullet}="a8"
!{(-3,2)}*+{\bullet}="a9"
!{(4,1.5)}*+{\circ}="b1"
!{(4,1)}*+{\circ}="b2"
!{(4,0.5)}*+{\circ}="b3"
!{(3,1)}*+{\circ}="b4"
!{(2,2)}*+{\circ}="b5"
!{(2,0)}*+{\circ}="b6"
!{(1,1)}*+{\circ}="b7"
!{(-1,1)}*+{\circ}="b8"
!{(-2,0)}*+{\circ}="b9"
!{(-4,2)}*+{\circ}="b10"
!{(4.2,1.15)}*+{^1}
!{(3.2,1.5)}*+{^2}
!{(3.2,0.5)}*+{^3}
!{(1.2,1.5)}*+{^4}
!{(1.2,0.5)}*+{^5}
!{(-0.8,1.5)}*+{^{2k}}
!{(-0.6,0.5)}*+{^{2k+1}}
!{(-1.5,0.2)}*+{^{2k+2}}
"a4"-"b1" "a4"-"b2" "a4"-"b3" "b1"-"a5"
"b3"-"a3" "a5"-"b4" "b4"-"a3" "a5"-"b5"
"a6"-"b5" "a3"-"b6" "b6"-"a2" "a6"-"b7"
"b7"-"a2" "a2"-"c1" "a6"-"c2" "c1"-"a1"
"c2"-"a7" "a7"-"b8" "b8"-"a1" "a1"-"b9"
"a7"-"a8" "a8"-@/^/"a9" "a8"-@/_/"a9" "a9"-"b10"
}
\end{displaymath}
If this monodromy group is imprimitive, then let $\Sigma_0=(1)(23)$, $\Sigma_1=(123)$ and $\Sigma_{\infty}=(12)(3)$ which realizes $H_{1/2,1/3,1/4}$ as a pullback of $H_{1/2,1/2,1/4}$. Also, let $P_1,P_2,P_3$ be a partition such that $\sigma_{\alpha}(P_i)=P_{\Sigma_{\alpha}(i)}$, where $\sigma_i$ be the monodromy of the dessin above $i$. To prove such partition does not exists, it followed by some observations. First, $1$ and $2k+2$ lie in the same partition, since they are fixed by $\sigma_0$. Second, $1$, $2$ and $3$ lie in same partition say $P_1$. Also, $4$ and $5$ lie in $P_1$ for the same reason. By continuing the process, we have $2k$ and $2k+1$ lie in $P_1$. Now, $2k+1$ and $2k+2$ lie in $P_1$, a contradiction.

Now, we label the $\circ$ which bound the outer face counterclockwise. 
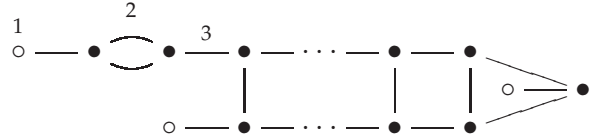
\begin{displaymath}
\xygraph{
!{<0cm,0cm>;<1cm,0cm>:<0cm,1cm>::}
!{(1,0)}*+{\bullet}="a1"
!{(2,0)}*+{\bullet}="a2"
!{(3.5,-0.5)}*+{\bullet}="a3"
!{(2,-1)}*+{\bullet}="a4"
!{(1,-1)}*+{\bullet}="a5"
!{(-1,0)}*+{\bullet}="a6"
!{(-1,-1)}*+{\bullet}="a7"
!{(0,0)}*+{\cdots}="c1"
!{(0,-1)}*+{\cdots}="c2"
!{(2.5,-0.5)}*+{\circ}="b1"
!{(-2,0)}*+{\bullet}="b2"
!{(-2,-1)}*+{\circ}="b3"
!{(-3,0)}*+{\bullet}="b4"
!{(-4,0)}*+{\circ}="b5"
!{(-4,0.3)}*+{^1}
!{(-2.5,0.5)}*+{^2}
!{(-1.5,0.2)}*+{^3}
"c1"-"a1" "c1"-"a6" "c2"-"a5" "c2"-"a7"
"a1"-"a5" "a1"-"a2" "a2"-"a4" "a4"-"a5"
"a2"-"a3" "a3"-"a4" "a3"-"b1"
"a6"-"a7" "a6"-"b2" "a7"-"b3" "b2"-@/_/"b4"
"b2"-@/^/"b4" "b4"-"b5"
}
\end{displaymath} 
Then by gluing $\circ^1,\circ^5,\cdots,\circ^{4n_0+1}$ together, we get the dessin as desired. The method to show this dessin has primitive monodromy is the same as the method given above.

\item 
($PM\cong S_4$ with $n_0\in\mathbb{Z}$, $n_1\in\lbrace\pm1/6,\pm1/3\rbrace+\mathbb{Z}$)

The method is similar with the method used in case (3). We just sketch the proof. We first recall the dessin given by the case 2 of the proof of Theorem \ref{BM2}.
\begin{displaymath}
\xygraph{
!{<0cm,0cm>;<1cm,0cm>:<0cm,1cm>::}
!{(1,0)}*+{\bullet}="a1"
!{(2,0)}*+{\bullet}="a2"
!{(2,-1)}*+{\bullet}="a4"
!{(1,-1)}*+{\bullet}="a5"
!{(-1,0)}*+{\bullet}="a6"
!{(-2,0)}*+{\bullet}="a7"
!{(-1,-1)}*+{\bullet}="a9"
!{(0,0)}*+{\cdots}="c1"
!{(-2,-1)}*+{\circ}="c2"
!{(-3,-0.5)}*+{\circ}="c3"
!{(0,-1)}*+{\cdots}="c4"
!{(-3,0)}*+{\circ^1}="b2"
!{(2.7,-0.5)}*+{\circ}="b7"
"c1"-"a1" "c1"-"a6" 
"a1"-"a5" "a1"-"a2" "a2"-"a4" "a4"-"a5"
"a6"-"a7" "a6"-"a7" "b2"-"a7" 
"a6"-"a9" "a7"-"a9" "a1"-"a4" "a2"-@/^1cm/"a4"
"a2"-"b7" "a9"-"c2" "a7"-"c3"
"c4"-"a9" "c4"-"a6" "c4"-"a5"
}
\end{displaymath}
Then we begin with $\circ^1$ and label $\circ$ which bound the outer face counterclockwise. By gluing $\circ^1,\circ^4,\cdots,\circ^{3n_0+1}$ together, we get the dessin as desired. It is not hard to show this dessin has primitive monodromy group.

\item 
($PM\cong S_4$ with $n_0\in\lbrace 1/2\rbrace+\mathbb{Z}$, $n_1\in\lbrace\pm1/8,\pm3/8\rbrace+\mathbb{Z}$)

The construction also uses the gluing method. We sketch the proof. We recall the dessin given in case 5 of the proof of Theorem \ref{BM2}.
\begin{displaymath}
\xygraph{
!{<0cm,0cm>;<1cm,0cm>:<0cm,1cm>::}
!{(1,0)}*+{\bullet}="a1"
!{(2,0)}*+{\bullet}="a2"
!{(3.5,-0.5)}*+{\bullet}="a3"
!{(2,-1)}*+{\bullet}="a4"
!{(1,-1)}*+{\bullet}="a5"
!{(-1,0)}*+{\bullet}="a6"
!{(-1,-1)}*+{\bullet}="a7"
!{(0,0)}*+{\cdots}="c1"
!{(0,-1)}*+{\cdots}="c2"
!{(2.5,-0.5)}*+{\circ}="b1"
!{(-2,0)}*+{\bullet}="b2"
!{(-2,-1)}*+{\circ}="b3"
!{(-3,0)}*+{\bullet}="b4"
!{(-4,0)}*+{\circ}="b5"
!{(-2.5,0.5)}*+{^1}
!{(-1.5,0.2)}*+{^2}
"c1"-"a1" "c1"-"a6" "c2"-"a5" "c2"-"a7"
"a1"-"a5" "a1"-"a2" "a2"-"a4" "a4"-"a5"
"a2"-"a3" "a3"-"a4" "a3"-"b1"
"a6"-"a7" "a6"-"b2" "a7"-"b3" "b2"-@/_/"b4"
"b2"-@/^/"b4" "b4"-"b5"
}
\end{displaymath}
Then we begin with $\circ^1$ and label $\circ$ which bound the outer face counterclockwise. By gluing $\circ^1,\circ^5,\cdots,\circ^{4n_0-1}$ together, we get the dessin as desired. The method to show that this dessin has primitive monodromy group is the same as the method given in case (3).

\item 
($PM\cong S_4$ with $n_0\in\lbrace 1/2\rbrace+\mathbb{Z}$, $n_1\in\lbrace\pm1/6,\pm1/3\rbrace+\mathbb{Z}$)

In this case, we don't need to worried about whether the dessin has primitive monodromy group, since $PM$ is never isomorphic to $D_2,D_3$ or $D_4$, the subgroup of $S_4$. The construction also uses the gluing method. We recall the dessin given in case 2 of the proof of Theorem \ref{BM2}.
\begin{displaymath}
\xygraph{
!{<0cm,0cm>;<1cm,0cm>:<0cm,1cm>::}
!{(1,0)}*+{\bullet}="a1"
!{(2,0)}*+{\bullet}="a2"
!{(2,-1)}*+{\bullet}="a4"
!{(1,-1)}*+{\bullet}="a5"
!{(-1,0)}*+{\bullet}="a6"
!{(-2,0)}*+{\bullet}="a7"
!{(-1,-1)}*+{\bullet}="a9"
!{(0,0)}*+{\cdots}="c1"
!{(-2,-1)}*+{\circ}="c2"
!{(-3,-0.5)}*+{\circ}="c3"
!{(0,-1)}*+{\cdots}="c4"
!{(-3,0)}*+{\circ}="b2"
!{(2.7,-0.5)}*+{\circ}="b7"
"c1"-"a1" "c1"-"a6" 
"a1"-"a5" "a1"-"a2" "a2"-"a4" "a4"-"a5"
"a6"-"a7" "a6"-"a7" "b2"-"a7" 
"a6"-"a9" "a7"-"a9" "a1"-"a4" "a2"-@/^1cm/"a4"
"a2"-"b7" "a9"-"c2" "a7"-"c3"
"c4"-"a9" "c4"-"a6" "c4"-"a5"
}
\end{displaymath}
Then gluing $\circ$ or $\bullet$ together can both give our desired dessin.

\item 
($PM\cong S_4$ with $n_0\in\lbrace \pm1/4\rbrace+\mathbb{Z}$, $n_1\in\lbrace\pm1/8,\pm3/8\rbrace+\mathbb{Z}$)

The dessin we need has two points with large odd degree over $\infty$.

We first construct the dessin 
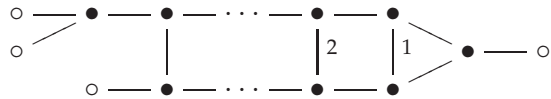
\begin{displaymath}
\xygraph{
!{<0cm,0cm>;<1cm,0cm>:<0cm,1cm>::}
!{(1,0)}*+{\bullet}="a1"
!{(2,0)}*+{\bullet}="a2"
!{(4,-0.5)}*+{\circ}="b1"
!{(2,-1)}*+{\bullet}="a4"
!{(1,-1)}*+{\bullet}="a5"
!{(-1,0)}*+{\bullet}="a6"
!{(-1,-1)}*+{\bullet}="a7"
!{(2.2,-0.5)}*+{^1}
!{(1.2,-0.5)}*+{^2}
!{(0,0)}*+{\cdots}="c1"
!{(0,-1)}*+{\cdots}="c2"
!{(3,-0.5)}*+{\bullet}="a3"
!{(-2,0)}*+{\bullet}="b2"
!{(-2,-1)}*+{\circ}="b3"
!{(-3,0)}*+{\circ}="b4"
!{(-3,-0.5)}*+{\circ}="b5"
"c1"-"a1" "c1"-"a6" "c2"-"a5" "c2"-"a7"
"a1"-"a5" "a1"-"a2" "a2"-"a4" "a4"-"a5"
"a2"-"a3" "a3"-"a4" "a3"-"b1"
"a6"-"a7" "a6"-"b2" "a7"-"b3"
"b2"-"b4" "b2"-"b5"
}
\end{displaymath}
Notice that this dessin has two faces which have odd degree. One of them has degree $3$. If we delete the line $1$ and draw it outside, we get the following dessin 
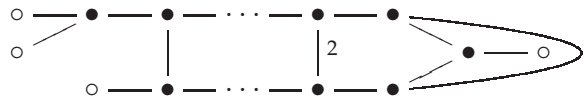
\begin{displaymath}
\xygraph{
!{<0cm,0cm>;<1cm,0cm>:<0cm,1cm>::}
!{(1,0)}*+{\bullet}="a1"
!{(2,0)}*+{\bullet}="a2"
!{(4,-0.5)}*+{\circ}="b1"
!{(2,-1)}*+{\bullet}="a4"
!{(1,-1)}*+{\bullet}="a5"
!{(-1,0)}*+{\bullet}="a6"
!{(-1,-1)}*+{\bullet}="a7"
!{(1.2,-0.5)}*+{^2}
!{(0,0)}*+{\cdots}="c1"
!{(0,-1)}*+{\cdots}="c2"
!{(3,-0.5)}*+{\bullet}="a3"
!{(-2,0)}*+{\bullet}="b2"
!{(-2,-1)}*+{\circ}="b3"
!{(-3,0)}*+{\circ}="b4"
!{(-3,-0.5)}*+{\circ}="b5"
"c1"-"a1" "c1"-"a6" "c2"-"a5" "c2"-"a7"
"a1"-"a5" "a1"-"a2" "a4"-"a5"
"a2"-"a3" "a3"-"a4" "a3"-"b1"
"a6"-"a7" "a6"-"b2" "a7"-"b3"
"b2"-"b4" "b2"-"b5" "a2"-@/^2.5cm/"a4"
}
\end{displaymath}
This dessin also has two faces which have odd degree. One of them has degree $5$. Notice that we can do the same process on the line $2$. Hence we can construct a dessin with two faces have large odd degree. This finish the construction.

\item 
($PM\cong S_4$ with $n_0\in\lbrace \pm1/4\rbrace+\mathbb{Z}$, $n_1\in\lbrace\pm1/8,\pm3/8\rbrace+\mathbb{Z}$)

The dessin we need has two points with large degree over $\infty$. One of them has odd degree and the other has even.

We first construct the dessin
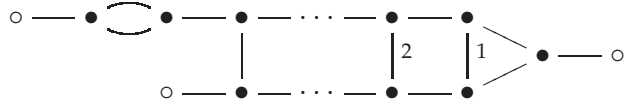
\begin{displaymath}
\xygraph{
!{<0cm,0cm>;<1cm,0cm>:<0cm,1cm>::}
!{(1,0)}*+{\bullet}="a1"
!{(2,0)}*+{\bullet}="a2"
!{(4,-0.5)}*+{\circ}="b1"
!{(2,-1)}*+{\bullet}="a4"
!{(1,-1)}*+{\bullet}="a5"
!{(-1,0)}*+{\bullet}="a6"
!{(-1,-1)}*+{\bullet}="a7"
!{(0,0)}*+{\cdots}="c1"
!{(0,-1)}*+{\cdots}="c2"
!{(3,-0.5)}*+{\bullet}="a3"
!{(-2,0)}*+{\bullet}="b2"
!{(-2,-1)}*+{\circ}="b3"
!{(-3,0)}*+{\bullet}="b4"
!{(-4,0)}*+{\circ}="b5"
!{(2.2,-0.5)}*+{^1}
!{(1.2,-0.5)}*+{^2}
"c1"-"a1" "c1"-"a6" "c2"-"a5" "c2"-"a7"
"a1"-"a5" "a1"-"a2" "a2"-"a4" "a4"-"a5"
"a2"-"a3" "a3"-"a4" "a3"-"b1"
"a6"-"a7" "a6"-"b2" "a7"-"b3" "b2"-@/_/"b4"
"b2"-@/^/"b4" "b4"-"b5"
}
\end{displaymath} 
This dessin has a face with degree 3 and a face with large even degree. We delete the line $1$ and draw it outside as in case (7). We get a dessin with a degree $3$ face. We can do the same process on the line 2. Hence, we can construct a dessin with a face with large odd degree and a face with large even degree. We complete the construction.
\end{enumerate}

By using the gluing method, it is not hard to construct the dessin corresponding to case (9), (10), (11), (12) and (13). We skip the detail.
 
\begin{enumerate}[label=\textbf{case (14) :}]
\item 
($PM\cong A_5$ with $n_0\in\lbrace \pm1/6\rbrace+\mathbb{Z}$, $n_1\in\lbrace\pm1/3,\pm1/6\rbrace+\mathbb{Z}$)

The dessin we need has two points with large ramification index over $1$. We recall the the dessin given in case 3 of the proof of Theorem \ref{BM}. Then label the $\circ$ which bound the outer face clockwise.
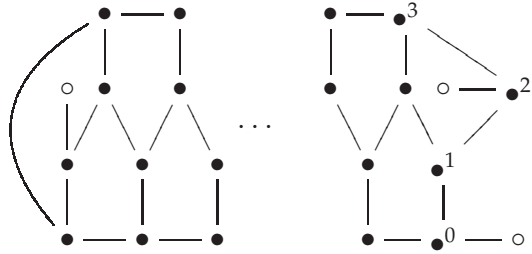
\begin{displaymath}
\xygraph{
!{<0cm,0cm>;<1cm,0cm>:<0cm,1cm>::}
!{(0,-0.5)}*+{\cdots}="c0"
!{(1,0)}*+{\bullet}="a1"
!{(2,0)}*+{\bullet}="a2"
!{(1,1)}*+{\bullet}="a3"
!{(2,1)}*+{\bullet^3}="a4"
!{(1.5,-1)}*+{\bullet}="a5"
!{(2.5,-1)}*+{\bullet^1}="a6"
!{(1.5,-2)}*+{\bullet}="a7"
!{(2.5,-2)}*+{\bullet^0}="a8"
!{(-1,0)}*+{\bullet}="a9"
!{(-1,1)}*+{\bullet}="a11"
!{(-2,0)}*+{\bullet}="a10"
!{(-2,1)}*+{\bullet}="a12"
!{(-1.5,-1)}*+{\bullet}="a13"
!{(-2.5,-1)}*+{\bullet}="a14"
!{(-1.5,-2)}*+{\bullet}="a15"
!{(-2.5,-2)}*+{\bullet}="a16"
!{(-2.5,0)}*+{\circ}="b1"
!{(3.5,0)}*+{\bullet^2}="a17"
!{(2.5,0)}*+{\circ}="b2"
!{(3.5,-2)}*+{\circ}="b3"
!{(-0.5,-1)}*+{\bullet}="a18"
!{(-0.5,-2)}*+{\bullet}="a19"
"a1"-"a3" "a3"-"a4" "a4"-"a2" "a1"-"a5" "a2"-"a5" "a5"-"a7"
"a7"-"a8" "a8"-"a6" "a6"-"a2" "a15"-"a16" "a14"-"a16" "a10"-"a14"
"a9"-"a11" "a11"-"a12" "a12"-"a10" "a9"-"a13" "a10"-"a13" "a13"-"a15"
"a12"-@/_1cm/"a16" "a14"-"b1" "a4"-"a17" "a6"-"a17" "a17"-"b2"
"a8"-"b3" "a9"-"a18" "a18"-"a19" "a19"-"a15"
}
\end{displaymath}
Then by gluing $\circ^1,\circ^6,\cdots$ together, we get a point with high ramification index divided by $3$. However, in our case, we need to construct the point with large ramification index not divided by $3$. To solve this problem, we replace the term 
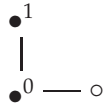
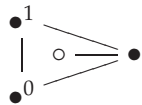
\begin{displaymath}
\xygraph{
!{<0cm,0cm>;<1cm,0cm>:<0cm,1cm>::}
!{(0,1)}*+{\bullet^1}="a1"
!{(0,0)}*+{\bullet^0}="a2"
!{(1,0)}*+{\circ}="b1"
"a2"-"b1" "a1"-"a2"
} 
\end{displaymath} 
by
\begin{displaymath}
\xygraph{
!{<0cm,0cm>;<1cm,0cm>:<0cm,1cm>::}
!{(0,1)}*+{\bullet^1}="a1"
!{(0,0)}*+{\bullet^0}="a2"
!{(1.5,0.5)}*+{\bullet}="a3"
!{(0.5,0.5)}*+{\circ}="b1"
"a2"-"a1" "a3"-"a1" "a3"-"a2" "a3"-"b1"
} 
\end{displaymath}
After doing this process, the degree of $\bullet^1$ plus $1$. Hence we can continue doing this process to get the point which has large degree above $1$.
\end{enumerate}

\begin{enumerate}[label=\textbf{case (15) :}]

\item 
($PM\cong A_5$ with $n_0\in\lbrace\pm1/6\rbrace+\mathbb{Z}$, $n_1\in (\mathbb{Z}/10)\setminus(\mathbb{Z}/2)$)

The dessin we need has one point with large degree above $1$. The method given in case (14) works and the construction is the same.
\end{enumerate}

\begin{enumerate}[label=\textbf{case (16) :}]

\item 
($PM\cong A_5$ with $n_0\in\lbrace\pm1/10,3/10\rbrace+\mathbb{Z}$, $n_1\in (\mathbb{Z}/10)\setminus(\mathbb{Z}/2)$)

This dessin we need has two points with large degree above $\infty$. We recall the dessin given in case 3 in the proof of Theorem \ref{BM}
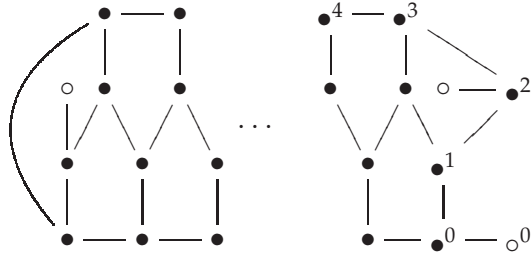
\begin{displaymath}
\xygraph{
!{<0cm,0cm>;<1cm,0cm>:<0cm,1cm>::}
!{(0,-0.5)}*+{\cdots}="c0"
!{(1,0)}*+{\bullet}="a1"
!{(2,0)}*+{\bullet}="a2"
!{(1,1)}*+{\bullet^4}="a3"
!{(2,1)}*+{\bullet^3}="a4"
!{(1.5,-1)}*+{\bullet}="a5"
!{(2.5,-1)}*+{\bullet^1}="a6"
!{(1.5,-2)}*+{\bullet}="a7"
!{(2.5,-2)}*+{\bullet^0}="a8"
!{(-1,0)}*+{\bullet}="a9"
!{(-1,1)}*+{\bullet}="a11"
!{(-2,0)}*+{\bullet}="a10"
!{(-2,1)}*+{\bullet}="a12"
!{(-1.5,-1)}*+{\bullet}="a13"
!{(-2.5,-1)}*+{\bullet}="a14"
!{(-1.5,-2)}*+{\bullet}="a15"
!{(-2.5,-2)}*+{\bullet}="a16"
!{(-2.5,0)}*+{\circ}="b1"
!{(3.5,0)}*+{\bullet^2}="a17"
!{(2.5,0)}*+{\circ}="b2"
!{(3.5,-2)}*+{\circ^0}="b3"
!{(-0.5,-1)}*+{\bullet}="a18"
!{(-0.5,-2)}*+{\bullet}="a19"
"a1"-"a3" "a3"-"a4" "a4"-"a2" "a1"-"a5" "a2"-"a5" "a5"-"a7"
"a7"-"a8" "a8"-"a6" "a6"-"a2" "a15"-"a16" "a14"-"a16" "a10"-"a14"
"a9"-"a11" "a11"-"a12" "a12"-"a10" "a9"-"a13" "a10"-"a13" "a13"-"a15"
"a12"-@/_1cm/"a16" "a14"-"b1" "a4"-"a17" "a6"-"a17" "a17"-"b2"
"a8"-"b3" "a9"-"a18" "a18"-"a19" "a19"-"a15"
}
\end{displaymath}
To construct a face with large degree, we delete $\circ^0$ and the line connect $\bullet^k\bullet^{k+1}$. Then we connect $\bullet^0$ with $\bullet^{k+1}$ in a line and connect $\circ^0$ with $\bullet^k$. For example, if we choose $k=2$, we get the dessin 
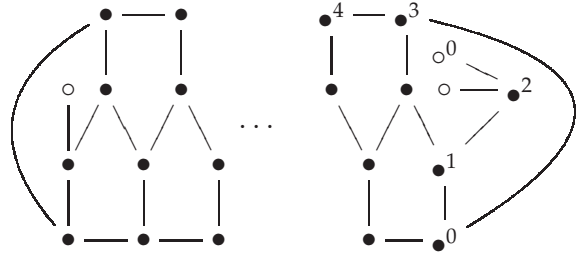
\begin{displaymath}
\xygraph{
!{<0cm,0cm>;<1cm,0cm>:<0cm,1cm>::}
!{(0,-0.5)}*+{\cdots}="c0"
!{(1,0)}*+{\bullet}="a1"
!{(2,0)}*+{\bullet}="a2"
!{(1,1)}*+{\bullet^4}="a3"
!{(2,1)}*+{\bullet^3}="a4"
!{(1.5,-1)}*+{\bullet}="a5"
!{(2.5,-1)}*+{\bullet^1}="a6"
!{(1.5,-2)}*+{\bullet}="a7"
!{(2.5,-2)}*+{\bullet^0}="a8"
!{(-1,0)}*+{\bullet}="a9"
!{(-1,1)}*+{\bullet}="a11"
!{(-2,0)}*+{\bullet}="a10"
!{(-2,1)}*+{\bullet}="a12"
!{(-1.5,-1)}*+{\bullet}="a13"
!{(-2.5,-1)}*+{\bullet}="a14"
!{(-1.5,-2)}*+{\bullet}="a15"
!{(-2.5,-2)}*+{\bullet}="a16"
!{(-2.5,0)}*+{\circ}="b1"
!{(3.5,0)}*+{\bullet^2}="a17"
!{(2.5,0)}*+{\circ}="b2"
!{(2.5,0.5)}*+{\circ^0}="b3"
!{(-0.5,-1)}*+{\bullet}="a18"
!{(-0.5,-2)}*+{\bullet}="a19"
"a1"-"a3" "a4"-"a2" "a1"-"a5" "a2"-"a5" "a5"-"a7"
"a7"-"a8" "a8"-"a6" "a6"-"a2" "a15"-"a16" "a14"-"a16" "a10"-"a14"
"a9"-"a11" "a11"-"a12" "a12"-"a10" "a9"-"a13" "a10"-"a13" "a13"-"a15"
"a12"-@/_1cm/"a16" "a14"-"b1" "a6"-"a17" "a17"-"b2" "a3"-"a4" "a8"-@/_2cm/"a4" "b3"-"a17"
"a9"-"a18" "a18"-"a19" "a19"-"a15"
}
\end{displaymath}
This construction can construct a dessin with 2 faces which has large degree. We finish the construction. 
\end{enumerate}
We finish all the cases.
\end{proof}

\begin{remark} \label{R. Gene}
In the proof of Theorem \ref{3-pts M}, we only give the construction for the points with large degree. In other words, $n_0,n_1 \geq 1/6$. In this case, the gluing method gives a systematic way to construct the corresponding dessin. For the case $n_0$ or $n_1<1/6$, we can just list all the possible as in the case (1) of the proof of Theorem \ref{3-pts M}. The existence of the corresponding dessin of all the cases are not hard. We can further generalize the problem to the equation (\ref{C-Alg-Lame}) with $2s+1$ regular singular points. The gluing method can also be generalized and give a systematic construction for $n_i\geq 1/6$. The problem raises when some $n_i<1/6$. The gluing method seems not work. It is still unsolved now.
\end{remark}

\section*{Acknowledgements}
I am grateful to the Taida Institute for Mathematical Sciences (TIMS) for organizing the very helpful \emph{TIMS Undergraduate Seminar} from September 2013 to February 2014. It gave me the chance to study various topics related to this research and I had benefited a lot from discussions with those participants. I am particularly grateful to my \emph{Undergraduate Thesis} adviser Professor Chin-Lung Wang. He guided me to study further recent results and gave me many useful suggestions and comments throughout my research on the finite monodromy problem. I am also grateful to Professor Ching-Li Chai for discussions on the Grothendieck correspondence during the preparation of this paper.

\end{document}